\documentclass[11pt]{amsart}
\usepackage[utf8]{inputenc}
\usepackage{amsmath,amssymb,amsthm,mathtools}
\usepackage[margin=2.7cm]{geometry}
\usepackage{tikz}
\usepackage{hyperref}

\newtheorem{theorem}{Theorem}[section]
\newtheorem{proposition}[theorem]{Proposition}
\newtheorem{lemma}[theorem]{Lemma}
\newtheorem{corollary}[theorem]{Corollary}

\newtheorem*{theorem*}{Theorem}
\theoremstyle{definition}
\newtheorem{definition}[theorem]{Definition}

\theoremstyle{remark}
\newtheorem{remark}[theorem]{Remark}
\newtheorem{example}[theorem]{Example}

\newcommand{\K}{\mathbb{K}}
\newcommand{\C}{\mathbb{C}}
\newcommand{\R}{\mathbb{R}}
\newcommand{\N}{\mathbb{N}}
\newcommand{\Z}{\mathbb{Z}}
\newcommand{\Q}{\mathbb{Q}}
\newcommand{\cO}{\mathcal{O}}
\newcommand{\cF}{\mathcal{F}}
\newcommand{\cP}{\mathcal{P}}
\newcommand{\br}[1]{\langle #1\rangle}

\title[Hill's lunar equation as a leaf]{Hill's Lunar Equation as a Leaf: Riccati Reduction and Frobenius Theory for Second-Order PDEs}
\author{V\'ictor Le\'on}
\address{ILACVN, Universidade Federal da Integra\c{c}\~ao Latino-Americana (UNILA), Foz do Igua\c{c}u, PR, Brazil}
\email{victor.leon@unila.edu.br}

\author{Bruno Sc\'ardua}
\address{Instituto de Matem\'atica, Universidade Federal do Rio de Janeiro (UFRJ), CP 68530, Rio de Janeiro, RJ, 21945-970, Brazil}
\email{bruno.scardua@gmail.com}


\numberwithin{equation}{section}
\begin{document}

\begin{abstract}
The equations of Hill and Adams for the lunar perigee and node, Mathieu's equation and Lam\'e's equation arise as leafwise restrictions of
one class of linear partial differential equations, the \emph{square-type foliated equations}
$L^2u+A\,Lu+B\,u=0$, where $L$ is an analytic vector field and $A,B$ are analytic functions. Such an equation has a regular singularity
exactly when the foliation defined by $L$
is analytically linearizable. At
rational saddles, where every point of the indicial conic is resonant, we obtain convergent solutions whose
exponents are holomorphic functions of a first integral; their variation describes the monodromy along the
leaves and decides when classical Frobenius solutions exist.

For a Hill equation, the characteristic exponent becomes the value of such a holomorphic function, given by a
convergent series. This yields rigorous computer-assisted proofs of
Hill's constant for the motion of the lunar perigee and of Adams' constant for the motion of the node: Hill's
value is confirmed within his own error estimate and its fifteenth decimal is determined, his value of the
motion of the perigee is shown to be correct in all its digits, and Adams' value is confirmed to twelve
decimals. For Mathieu's equation, the curves of constant characteristic exponent that govern ion trapping
appear as classical Frobenius loci. For Lam\'e's equation, we prove a product formula for the logarithmic
obstructions, which coincide with the leading coefficients of the eigenvalue splittings computed by Volkmer and
whose vanishing reflects the finite-gap property.
\end{abstract}

\maketitle
\tableofcontents

\part*{Introduction}

Linear second order equations with periodic coefficients are among the oldest objects of mathematical
physics. Hill introduced his equation in 1877 to compute the motion of the lunar perigee
\cite{Hill1886}, and Adams had independently arrived at the same equation in his study of the motion of
the node \cite{Adams1877}. Mathieu's equation arose in 1868 from the vibrations of an elliptic membrane
\cite{Mathieu1868}; today it governs the motion of ions in radio-frequency (Paul) traps \cite{Paul1990}.
Lam\'e's equation, which goes back to Lam\'e's 1837 work on the equilibrium of heat in ellipsoids
\cite{Lame1837}, reappears as the equation of parametric resonance in models of preheating after
inflation \cite{GKLS1997}. In this paper these equations appear in a geometric way: each of them is the
restriction, to the leaves of a holomorphic foliation near a saddle singularity, of one linear partial
differential equation of the form \eqref{eq:intro} below, which we call a square-type foliated equation. The idea of placing Hill's equation inside a
holomorphic foliation goes back to \cite{ReisScardua2019}, where the complex Hill equation is associated
with an integrable complex analytic one-form in dimension three. Here this point of view is developed
for square-type foliated equations as a whole; it makes the characteristic exponents of these equations
holomorphic functions of a first integral, and leads to convergent and rigorously certified computations
of classical constants (Part~\ref{part:applications}).

The method of Frobenius is a classical tool for studying linear differential equations near regular
singular points. Its extension to partial differential equations requires a suitable notion of regular
singularity, an indicial equation and conditions under which the resulting formal solutions converge.
Fuchsian partial differential equations go back to Baouendi and Goulaouic \cite{BG73}, and the general
Fuchsian PDE setting has also been investigated by Mandai \cite{Mandai2000}. Formal, holomorphic and singular
solutions of singular nonlinear first order PDEs, including totally characteristic equations of the form
$(ax\partial_x+by\partial_y)v=F(x,y,v)$ and Briot--Bouquet type theorems for them, are treated systematically by
G\'erard and Tahara \cite{GT96}.
In \cite{LS21,LS23}, Frobenius type methods were developed for analytic linear PDEs of third and second order,
respectively. In particular, the second order theory in \cite{LS23} is formulated in terms of an Euler model,
an indicial conic and resonance conditions on its lattice points.

The present work is concerned with a particular class of second order PDEs, namely those obtained by applying
twice an analytic vector field, together with lower order terms. The special form of the operator makes it
possible to relate the Frobenius problem to a Riccati equation and to the local geometry of the vector field.
Our purpose is to determine how this structure affects regular singularities, convergence and resonance, and
to study in detail the case of rational saddles, where the standard non-resonance arguments are not available.

The Riccati viewpoint has also appeared in our earlier work on linear ordinary differential equations
\cite{LSERA}. Here we use it in the setting of PDEs, without requiring the global holonomy constructions of
that paper. The main new feature is the appearance, at rational saddles, of Frobenius exponents that are
holomorphic functions of a first integral rather than constants. Restricted to a regular leaf, which near the
saddle is an annulus, our solutions take the Floquet form $s^{\lambda(h_0)}\times(\text{Laurent series in }s)$, so
$\lambda(h_0)$ is a Floquet exponent of the leafwise equation \cite{Floquet1883,MagnusWinkler1966}, \cite[Chapter~3]{CL}. The point is that this exponent extends
holomorphically across the singular fiber $h=0$, where the annuli degenerate, and that the Floquet factors glue
together into a single holomorphic germ at the singular point. This leads to a precise criterion for
classical Frobenius solutions and to a locally holomorphic parametrization of the corresponding coefficients.
Unlike a parametrization of solutions for a fixed operator, this describes the locus of operators admitting
a prescribed classical Frobenius behavior; its finite compatibility truncations have explicitly determined codimension.

\subsection*{Square-type foliated equations}
Let $L=a\,\partial_x+b\,\partial_y$ be an analytic vector field and $A,B$ analytic functions. We call
\begin{equation}\label{eq:intro}
L^2u+A\,Lu+B\,u=0
\end{equation}
a \emph{square-type foliated equation} (STFE). The operator is a polynomial of degree two in the single vector
field $L$, so \eqref{eq:intro} is a family of linear second order ODEs, one along each leaf of the foliation
$\cF_L$ defined by $L$, depending analytically on the leaf. Written out, it reads
\[a^2u_{xx}+2ab\,u_{xy}+b^2u_{yy}+(L(a)+aA)u_x+(L(b)+bA)u_y+Bu=0.\]These are exactly the parabolic
equations whose principal part is the square of a vector field $L$ and whose first order part, corrected by
$L(a)\partial_x+L(b)\partial_y$, is a multiple of $L$. Two simple observations organize the paper. The
substitution $v=Lu/u$ transforms \eqref{eq:intro} into the \emph{Riccati PDE} $Lv=-v^2-Av-B$
(Proposition~\ref{prop:riccati}). On the other hand, at a regular point of $L$ the flow-box theorem
reduces \eqref{eq:intro} to a linear ODE with a parameter (Section~\ref{sec:riccati}), so the phenomena
studied here are concentrated near the \emph{singular points of $L$}, where the geometry of holomorphic
foliations becomes essential.
Square-type foliated equations are preserved by changes of coordinates and by multiplication of $L$ by a
unit (Lemma~\ref{lem:invariance}), so the notion depends only on the foliation $\cF_L$. The restriction of
\eqref{eq:intro} to a single leaf will be called a \emph{leafwise equation}.

\subsection*{Main results}
Let $0$ be an isolated singular point of $L$, with eigenvalues $a,b$ of $DL(0)$, and let $\lambda_1,\lambda_2$ be the
roots of $\lambda^2+A(0)\lambda+B(0)$, with $\mu=\lambda_2-\lambda_1$.

Following \cite{LS23}, we say that \eqref{eq:intro} has a \emph{regular singularity} at $0$ if, after a holomorphic
change of coordinates fixing $0$ and division by a holomorphic function, the operator takes the Euler-type form
\[
\alpha\,x^2\partial_x^2+\beta\,xy\,\partial_x\partial_y+\gamma\,y^2\partial_y^2+x\,a_1\,\partial_x+y\,b_1\,\partial_y+c_1,
\]
with constant principal coefficients $\alpha,\beta,\gamma$ and holomorphic $a_1,b_1,c_1$ (Definition~\ref{def:regular}).
In \cite{LS23} the principal constants are required to be real; here they may be complex, and we speak of
\emph{real principal constants} when they can be chosen real. The \emph{indicial conic} is the set of exponents
$(r,s)$ for which $x^ry^s$ solves the Euler-type model (Section~\ref{sec:euler}).

\begin{theorem*}[A --- regularity is linearizability]
The square-type foliated equation \eqref{eq:intro} has a regular singularity at $0$ if and only if the foliation $\cF_L$ is analytically linearizable with non-degenerate diagonal linear
part. The regular form can be chosen with real principal constants if and only if, in addition, the ratio of the
eigenvalues of $DL(0)$ is real. In that case the indicial
conic is the pair of parallel lines $\{ar+bs=\lambda_1\}\cup\{ar+bs=\lambda_2\}$. These lines are the preimages of
the equilibria of the Riccati PDE.
\end{theorem*}

This is Theorem~\ref{thm:characterization}. Its proof is short: for square-type foliated equations the principal
symbol is $L\otimes L$, and regularity is decided by this symbol alone. Its interest lies in the consequences.
Combined with the classical linearization theory (Poincar\'e,
Dulac, Brjuno, Mattei--Moussu, P\'erez-Marco--Yoccoz), it relates the local linearization properties of $\cF_L$ to the regular-singular Frobenius theory of \eqref{eq:intro} (Corollary~\ref{cor:dictionary}). Resonance of a
point of the indicial conic depends only on the line containing it, and it is decided by the arithmetic of $a/b$
and $\mu$ (Proposition~\ref{prop:resonance}). In particular, when $a/b\in\Q_{<0}$ (rational saddle) every point of
the conic is resonant.

\begin{theorem*}[B --- Poincar\'e domain]
Let $L$ be a holomorphic vector field near $0\in\C^2$ with $L(0)=0$. Assume that $DL(0)$ has eigenvalues $a,b$ with
$a/b\notin\R_{\le0}$, and that the foliation defined by $L$ is analytically linearizable, that is, there are holomorphic
coordinates $(x,y)$ centred at $0$ and a holomorphic function $g$ with $g(0)=1$ such that
$L=g\cdot(a\,x\partial_x+b\,y\partial_y)$. Let $A,B$ be holomorphic near $0$, let $\lambda_1,\lambda_2$ be the roots of
$\lambda^2+A(0)\lambda+B(0)$, put $\mu=\lambda_2-\lambda_1$, and let
$\Gamma=\{aq_1+bq_2:(q_1,q_2)\in\N^2\setminus\{0\}\}$. Then:
\begin{enumerate}
\item $\mu$ and $-\mu$ do not both belong to $\Gamma$.
\item If $\mu\notin\Gamma$, there is a unique formal power series $F_1(x,y)$ with $F_1(0)=1$ such that
$u=x^ry^sF_1(x,y)$ solves equation~\eqref{eq:intro} for every $(r,s)$ with $ar+bs=\lambda_1$. This series converges near $0$
and does not depend on $(r,s)$. Consequently $u=f(x^by^{-a})\,x^{\lambda_1/a}F_1$ is a solution, locally away from
the axes, for every holomorphic function $f$ of one variable. The same holds with $\lambda_1$ and $\mu$ replaced by
$\lambda_2$ and $-\mu$.
\item Every holomorphic germ $v$ solving the Riccati equation $Lv=-v^2-Av-B$ satisfies $v(0)\in\{\lambda_1,\lambda_2\}$.
If $\mu\notin\Gamma$, there is exactly one such germ with $v(0)=\lambda_1$, namely $v=Lu_1/u_1$ with
$u_1=x^{\lambda_1/a}F_1$. The same holds for $\lambda_2$ when $-\mu\notin\Gamma$.
\end{enumerate}
\end{theorem*}

In the language of Section~\ref{sec:euler}, the condition $\mu\notin\Gamma$ means that the indicial line
$\{ar+bs=\lambda_1\}$ is non-resonant, so item~(1) of the theorem says that at least one of the two lines is non-resonant. Item~(3)
is a Briot--Bouquet type theorem for the Riccati PDE. The theorem combines Theorem~\ref{thm:poincare} and
Corollaries~\ref{cor:always} and~\ref{cor:riccati-unique} with Lemma~\ref{lem:invariance}, which reduces
\eqref{eq:intro} to the linear model; since $g(0)=1$, this reduction does not change $\lambda_1$, $\lambda_2$ and $\mu$.
The linearizability hypothesis holds automatically when $DL(0)$ is diagonalizable and neither $a/b$ nor $b/a$ is an
integer $\ge2$, by Poincar\'e's theorem (Corollary~\ref{cor:dictionary}).

\begin{theorem*}[C --- rational saddles]
Let $L=p\,x\partial_x-q\,y\partial_y$, with $p,q$ positive coprime integers, and let $h=x^qy^p$, a first integral of
$L$. Let $A,B$ be holomorphic near $0$, let $\lambda_1,\lambda_2$ be the roots of $\lambda^2+A(0)\lambda+B(0)$, and
assume $\mu=\lambda_2-\lambda_1\notin\Z$. Then equation~\eqref{eq:intro} has convergent solutions
\[
u_j=x^{\lambda_j(h)/p}\,G_j(x,y),\qquad \lambda_j\in\C\{h\},\ \lambda_j(0)=\lambda_j,\ G_j\in\C\{x,y\},\ G_j(0)=1,
\]
($j=1,2$), and, locally on $\{xy\neq0\}$ near $0$, every solution is $f(h)u_1+g(h)u_2$. The exponent $\lambda_j(h)$
varies from leaf to leaf. For small $h_0\neq0$ and a small polydisc $\Delta_\rho$, the local leaf $\{h=h_0\}\cap\Delta_\rho$ is
an annulus, and analytic continuation of $u_j$ around it multiplies $u_j$ by $e^{2\pi i\lambda_j(h_0)}$. Finally,
equation~\eqref{eq:intro} has a classical Frobenius solution $x^ry^sF(x,y)$ with $pr-qs=\lambda_j$, $F$ holomorphic and $F(0)=1$,
if and only if $\lambda_j(h)$ is constant, that is, if and only if this leafwise monodromy eigenvalue is the same on
every leaf.
\end{theorem*}

This is Theorem~\ref{thm:rational} with Corollary~\ref{cor:rational}. By Lemma~\ref{lem:invariance}, the theorem
applies equally when the foliation of $L$ is analytically linearizable with this linear part, that is, when
$L=g\cdot(p\,x\partial_x-q\,y\partial_y)$ in suitable coordinates with $g(0)=1$, since the reduction does not change
$\lambda_1$ and $\lambda_2$. Two further results complete the picture
at rational saddles. First, the Riccati PDE has exactly one holomorphic solution through each equilibrium, as in the
Poincar\'e domain (Theorem~\ref{thm:riccati-saddle}). Second, classical Frobenius solutions are exceptional. For a fixed indicial line $\ell_1$,
write $B=B'+\beta(h)$, where $\beta$ collects the pure powers of $h$. The equations with a classical Frobenius
solution on $\ell_1$ are exactly those with $\beta=\Phi(A,B')$, for a map $\Phi$ which is holomorphic on
neighbourhoods in weighted $\ell^1$ spaces of power series, with a possible loss of convergence radius (Theorem~\ref{thm:isomonodromic}). The associated finite truncations define submanifolds of arbitrarily large finite codimension (Corollary~\ref{cor:generic}). The first obstruction is an explicit
coefficient $\lambda'(0)$ (Proposition~\ref{prop:first-coefficient}). The integral case $\mu\in\Z$ is treated for
one-sided coefficients in Section~\ref{sec:integral} (Theorem~\ref{thm:integral}). As in the classical theory, logarithms then appear
(Example~\ref{ex:bessel}, where $I_1$ and $K_1$ occur). If the classical ansatz with a constant exponent is
nevertheless used when $\lambda(h)$ is not constant, the expansion of $x^{\lambda(h)/p}$ produces \emph{unbounded}
powers of $\log x$ (Proposition~\ref{prop:saddle}). This shows that the constant-exponent ansatz is the wrong one,
whereas the exponent $\lambda(h)$ captures the solutions without any logarithm.

\subsection*{Applications to lunar theory and mathematical physics: the equations of Hill, Mathieu and Lam\'e}

\begin{figure}[ht]
\centering
\begin{tikzpicture}[scale=1.15,>=stealth]
  \draw[->,thick] (-2.9,0)--(3.1,0) node[right]{$X$};
  \draw[->,thick] (0,-2.9)--(0,3.1) node[above]{$Y$};
  \foreach \h in {0.15,0.45}{
    \draw[gray,domain=\h/2.8:2.8,samples=80] plot (\x,{\h/\x});
    \draw[gray,domain=-2.8:-\h/2.8,samples=80] plot (\x,{\h/\x});}
  \draw[very thick,domain=1/2.8:2.8,samples=80] plot (\x,{1/\x});
  \draw[very thick,domain=-2.8:-1/2.8,samples=80] plot (\x,{1/\x});
  \foreach \h in {2.2}{
    \draw[gray,dashed,domain=\h/2.8:2.8,samples=80] plot (\x,{\h/\x});
    \draw[gray,dashed,domain=-2.8:-\h/2.8,samples=80] plot (\x,{\h/\x});}
  \fill (0,0) circle (1.6pt);
  \node[right] at (2.05,0.62) {\small $XY=1$: Hill's equation};
  \draw[gray,->] (0.95,2.55)--(0.24,2.0);
  \node[gray,right] at (0.9,2.6) {\small $XY=H_0$};
  \node[left] at (-0.05,2.75) {\small separatrix $X=0$};
  \node[above] at (-2.1,0.05) {\small separatrix $Y=0$};
  \draw[->] (1.2,1.35) arc (60:350:0.25);
  \node at (1.65,1.65) {\small $e^{2\pi i\lambda(H_0)}$};
\end{tikzpicture}
\caption{Hill's equation as a leaf (real trace of the picture). The saddle $M=X\partial_X-Y\partial_Y$ has the
first integral $H=XY$; its leaves $\{XY=H_0\}$ are drawn in grey, and the separatrices $X=0$, $Y=0$ form the singular
leaf $H=0$. By Proposition~\ref{prop:hill-leaf}, the operator $M^2+B_\Theta$ induces on the leaf $\{XY=H_0\}$ Hill's
equation with the $k$-th Fourier mode scaled by $H_0^{k/2}$, and on the leaf $\{XY=1\}$ (thick) exactly Hill's
equation \eqref{eq:hill}. Each complex leaf is an annulus; going once around it multiplies the Floquet solution by
$e^{2\pi i\lambda(H_0)}$, where $\lambda(H)$ is holomorphic across the singular leaf $H=0$ (Corollary~\ref{cor:hill-exponent}).}
\label{fig:leaf}
\end{figure}
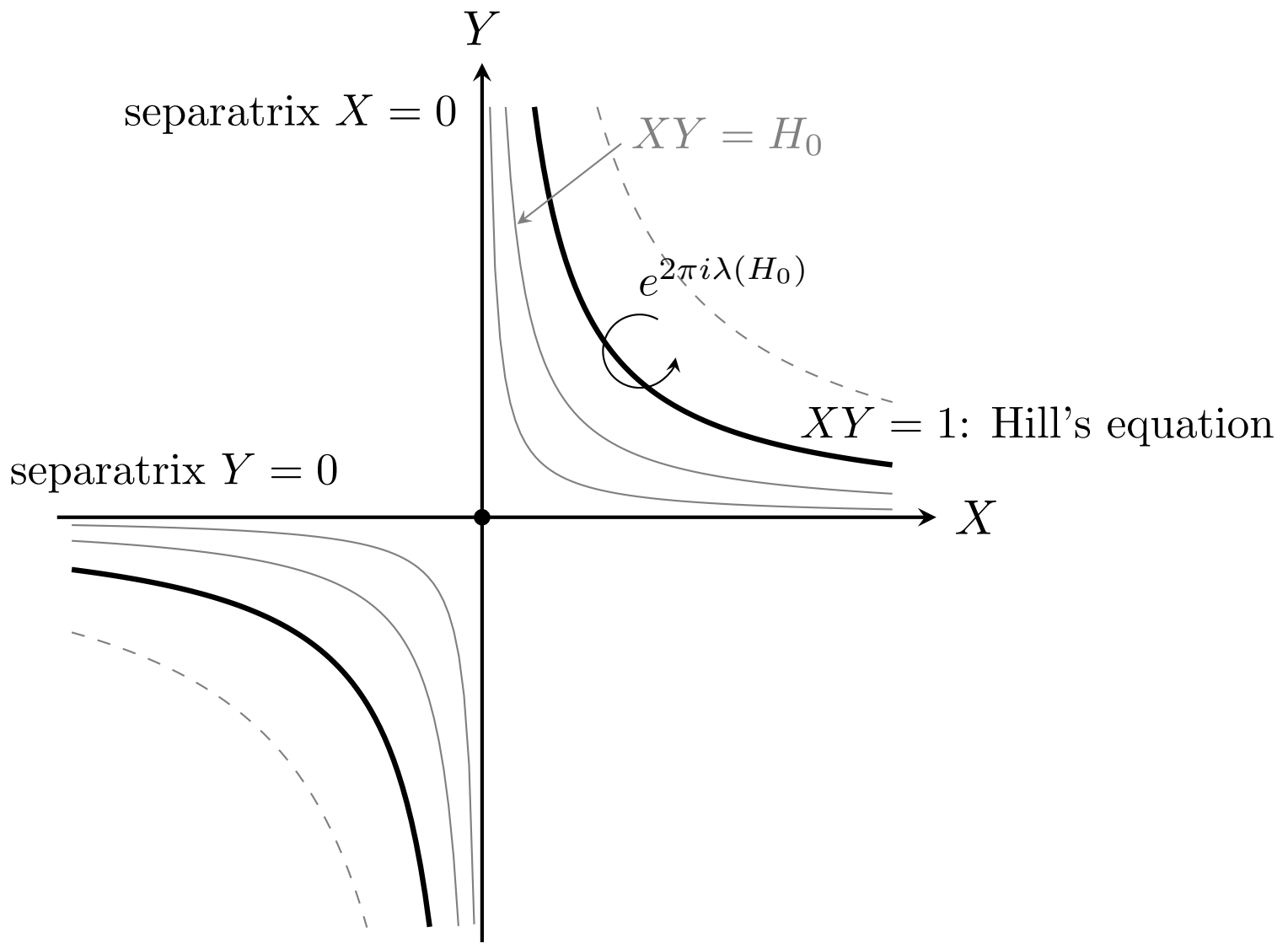

The structure described in Theorem~C has a concrete meaning for the classical equations mentioned at
the beginning of this introduction. Hill's equation $w''+\Theta(\tau)w=0$, with
$\Theta=\Theta_0+\sum_k\theta_k\cos2k\tau$, is not an isolated object from our point of view: it is the
equation induced on the leaf $\{XY=1\}$ by the operator $M^2+B_\Theta$, where $M=X\partial_X-Y\partial_Y$ and
$B_\Theta$ is built from the Fourier coefficients of $\Theta$. The neighbouring leaves $\{XY=H_0\}$ carry the
same equation with the $k$-th Fourier mode scaled by $H_0^{k/2}$ (Figure~\ref{fig:leaf}), so that Hill's equation becomes one
member of a holomorphic family, and its characteristic exponent becomes the value at $H=1$ of the
holomorphic exponent $\lambda(H)$ of Theorem~C.

\begin{theorem*}[D --- the lunar constants of Hill and Adams]
For the equations of Hill \cite{Hill1886} and Adams \cite{Adams1877}, with the coefficients printed in
their papers, the series $\lambda(H)$ converges up to $H=1$, and $2\lambda(1)$ is the characteristic exponent.
For Hill's equation this gives
\[
1.07158\,32774\,16010\,2695\le\mathbf c\le1.07158\,32774\,16010\,2723,
\]
and for Adams' equation $\mathbf g=1.08517\,13927\,46789\,1742\ldots$ with an error below $10^{-19}$.
\end{theorem*}

These enclosures are established by a rigorous computer-assisted proof, in which every inequality is
verified in exact rational arithmetic (Theorems~\ref{thm:hill-certified} and~\ref{thm:adams-certified}). The
conclusions hold not only for the coefficients as printed by Hill and Adams but for every set of data
consistent with their rounding, and the enclosures for the data as printed were confirmed independently by a
validated interval integration of both equations that does not use our recursion (Section~\ref{sec:hill-repro}). For Hill's data they determine
the fifteenth decimal of $\mathbf c$, which differs from Hill's printed value in the last digit. They also show
that his value of the motion of the perigee, $0.00857\,25730\,04864$ of the mean motion of the Moon (a
revolution of the perigee in about $8.73$ years), is correct in all its digits (Corollary~\ref{cor:hill-digits}).
For Adams' data they show that his value of $\mathbf g$ is correct to twelve decimals and that the last two of
the fourteen decimals he gives are not determined by his data, and they confirm his remark that Delaunay's
series in powers of $m$ is wrong in the eighth decimal. The recursion is also efficient: its first thirteen
terms already reproduce Hill's constant to fifteen decimals, and each new term gains about one decimal.
Table~\ref{tab:lunar} in Section~\ref{sec:moon} summarizes the comparison with the classical values.

The same mechanism applies to Mathieu's equation, which underlies the stability diagram of the Paul
trap. There, the curves along which the characteristic exponent stays constant are exactly the
classical Frobenius loci of Theorem~\ref{thm:isomonodromic}, and their coefficients are produced by the
recursion that defines these loci. For the iso-$\beta$ curves of the Paul trap, the values given by this
recursion agree to about $10^{-17}$ with those of the classical continued-fraction method, within a radius of
convergence that we determine and that shrinks as $\beta$ increases (Section~\ref{sec:paul}). For Lam\'e's equation the geometry is particularly transparent: for
all small values of the nome, the Lam\'e equations with a fixed value of $\Theta_0$ (that is, with the eigenvalue
parameter adjusted to the modulus) form the family of leafwise equations of a single explicit operator, whose potential is essentially a Weierstrass $\wp$-function written in the variables $X,Y$. The
first logarithmic obstruction of this operator at the $j$-th resonance turns out to be
\[
\kappa_j=\frac{(-1)^{j+1}}{((j-1)!)^2}\prod_{i=0}^{j-1}(\nu-i)(\nu+i+1)
\]
(Theorem~\ref{thm:lame}), and, comparing with a theorem of Volkmer, these numbers coincide, up to the factor
$8(-1)^{j+1}$, with the leading coefficients of the eigenvalue splittings at the instability intervals, whose absolute values are
the widths of these intervals (Corollary~\ref{cor:volkmer}). In
particular, for $\nu=n$ a non-negative integer the obstructions vanish exactly at the resonances $j>n$,
which is the leading-order form of the finite-gap property of Lam\'e's equation.

Our contribution to this classical subject is of four kinds. The first is geometric: Hill's, Mathieu's
and Lam\'e's equations become single leaves of holomorphic families attached to one operator at a saddle
point, and this explains why their characteristic exponents depend holomorphically on the size of the
perturbation, and why the Floquet solutions on nearby leaves glue together into one holomorphic germ.
The second is constructive: in place of Hill's infinite determinant we obtain a convergent recursion
that produces, in one pass and with a canonical choice of branch, both the characteristic exponent and
the Floquet solution. The third is rigorous: the same recursion, combined with an explicit contraction
argument, gives the certified enclosures described above. The fourth is structural: the characteristic
value curves of Mathieu's equation appear as classical Frobenius loci, and the leading coefficients of the
instability intervals of Lam\'e's equation appear as logarithmic obstructions, so that the finite-gap
property acquires a geometric meaning. At the same time, our statements concern the model equations as written
by Hill and Adams, not the full lunar problem, and the numerical corrections affect only the last
decimals and have no astronomical significance; the leading instability-interval formulas themselves are
due to Levy and Keller and to Volkmer \cite{LevyKeller1963,Volkmer2005}.

\subsection*{Structure of the paper}
The paper has two parts. Part~\ref{part:theory} (Sections~\ref{sec:riccati}--\ref{sec:siegel}) develops the theory of
square-type foliated equations: the Riccati reduction, the Euler model and its resonances, the Poincar\'e domain
(Theorem~B), the characterization of regular singularities (Theorem~A) and rational saddles (Theorem~C).
Part~\ref{part:applications} (Sections~\ref{sec:examples}--\ref{sec:lame}) contains explicit examples and the
applications: Hill's equation as a leaf, the lunar constants of Hill and Adams (Theorem~D), ion trapping and
Mathieu's equation, and Lam\'e's equation. Readers mainly interested in the applications can read
Section~\ref{sec:riccati} and then go directly to Part~\ref{part:applications}, which opens with a summary of the
results it uses; readers mainly interested in the theory can read Part~\ref{part:theory} and treat
Part~\ref{part:applications} as a collection of examples. We close with some remarks and open questions in
Section~\ref{sec:conclusion}.

\medskip
\noindent\textbf{Notation.} $\K\in\{\R,\C\}$, $\N=\{0,1,2,\dots\}$. For $Q=(q_1,q_2)\in\N^2$ we write
$X^Q=x^{q_1}y^{q_2}$, $|Q|=q_1+q_2$, and $P<Q$ means $P\le Q$ componentwise with $P\neq Q$.
$\cO_0$ denotes the ring of germs of holomorphic functions at $0\in\C^2$.

\part{Theory of square-type foliated equations}\label{part:theory}

\section{The Riccati reduction}\label{sec:riccati}

For a linear ordinary differential equation $u''+A(t)u'+B(t)u=0$, the substitution $v=u'/u$ leads to the
Riccati equation $v'=-v^2-Av-B$. This reduction goes back to the eighteenth century and is a standard tool of
the classical theory \cite{Ince1926}. It has a direct meaning for equations with periodic coefficients such as
Hill's: if $u=e^{\lambda t}p(t)$ is a Floquet solution, then $v=\lambda+p'/p$, so the characteristic exponent appears
as the constant part of a solution of the Riccati equation. In \cite{LSERA} we used this reduction, together with
holonomy arguments, to study linear second order ODEs in the complex domain. Square-type foliated equations are
the natural setting in which the same substitution works for PDEs. When the operator is a polynomial in a single
vector field $L$, the logarithmic derivative $v=Lu/u$ along $L$ again satisfies a Riccati equation, now a first
order PDE (Proposition~\ref{prop:riccati}). As we shall see, the leafwise exponents $\lambda_j(h)$ of
Section~\ref{sec:siegel} are exactly the constant parts, along each leaf, of the holomorphic solutions of this
Riccati PDE (Corollary~\ref{cor:rational}).

Let $U\subset\K^2$ be open and let $L=a(x,y)\,\partial_x+b(x,y)\,\partial_y$ be an analytic vector field on $U$.
For analytic $A,B\colon U\to\K$ consider the operator and the equation
\begin{equation}\label{eq:main}
\cP_{L,A,B}[u]:=L^2u+A\,Lu+B\,u=0 .
\end{equation}
Expanding, \[L^2u=a^2u_{xx}+2ab\,u_{xy}+b^2u_{yy}+L(a)\,u_x+L(b)\,u_y,\]so \eqref{eq:main} reads
\begin{equation}\label{eq:main-expanded}
a^2u_{xx}+2ab\,u_{xy}+b^2u_{yy}+\bigl(L(a)+aA\bigr)u_x+\bigl(L(b)+bA\bigr)u_y+B\,u=0 .
\end{equation}
In particular \eqref{eq:main} is parabolic wherever $L\neq0$, and its characteristic direction is $L$.

\begin{proposition}[Riccati reduction]\label{prop:riccati}
Let $u$ be analytic on $U$ and nowhere zero, and put $v=Lu/u$. Then
\[
\cP_{L,A,B}[u]=u\,\bigl(Lv+v^2+Av+B\bigr).
\]
Consequently $u$ solves \eqref{eq:main} if and only if $v$ solves the Riccati PDE
\begin{equation}\label{eq:riccati}
Lv=-v^2-A\,v-B .
\end{equation}
Conversely, if $v$ solves \eqref{eq:riccati} and $u$ is a nowhere-zero solution of the first order linear
equation $Lu=v\,u$, then $u$ solves \eqref{eq:main}.
\end{proposition}

\begin{proof}
$L$ is a derivation, so \[L^2u=L(vu)=(Lv)\,u+v\,Lu=(Lv+v^2)\,u.\] Hence
\[\cP_{L,A,B}[u]=(Lv+v^2)u+A\,vu+Bu.\]
The converse is the same computation read backwards.
\end{proof}

At a regular point, that is, if $L(p)\neq0$, the flow-box theorem gives local coordinates $(t,\sigma)$ with $L=\partial_t$. In these
coordinates \eqref{eq:main} becomes the second order linear ODE \[u_{tt}+A\,u_t+B\,u=0\] with parameter $\sigma$,
and \eqref{eq:riccati} becomes the classical Riccati ODE \[v_t=-v^2-Av-B.\] Thus, near regular points of $L$,
Proposition~\ref{prop:riccati} is the parametrized version of the classical reduction. For instance, for
$L=\partial_x+\partial_y$ one may take $t=y$ and $\sigma=x-y$. From now on we concentrate on
\emph{singular points} of $L$, where no such reduction is available.

Since our results are stated in coordinates where $L$ takes a simple form, we first check that square-type
foliated equations are preserved by changes of coordinates and by multiplying $L$ by a non-vanishing function.
The class therefore depends only on the foliation defined by $L$.

\begin{lemma}[Invariance under changes of coordinates and of the vector field]\label{lem:invariance}$\;$
\begin{enumerate}
\item If $\Phi$ is an analytic diffeomorphism, then
\[
\Phi_*\cP_{L,A,B}=\cP_{\Phi_*L,\;A\circ\Phi^{-1},\;B\circ\Phi^{-1}}.
\]
\item If $g$ is an analytic unit, then
\[
\cP_{gL,A,B}=g^2\,\cP_{L,\tilde A,\tilde B},\qquad \tilde A=\frac{A+L(g)}{g},\quad \tilde B=\frac{B}{g^2}.
\]
\end{enumerate}
In particular, whether an equation is a square-type foliated equation depends only on the \emph{foliation}
defined by $L$, not on $L$ itself.
\end{lemma}

\begin{proof} (1) Holds because $(Lu)\circ\Phi^{-1}=(\Phi_*L)(u\circ\Phi^{-1})$. (2) $(gL)^2u=gL(g\,Lu)=g^2L^2u+g\,L(g)\,Lu$, hence
\[\cP_{gL,A,B}[u]=g^2L^2u+\bigl(g\,L(g)+gA\bigr)Lu+Bu=g^2\,\cP_{L,\tilde A,\tilde B}[u].\]
\end{proof}

\section{The Euler model and its indicial conic}\label{sec:euler}

Fix $a,b\in\K^*$ and let $L_0=ax\,\partial_x+by\,\partial_y$. Then
\[L_0^2u=a^2x^2u_{xx}+2abxy\,u_{xy}+b^2y^2u_{yy}+a^2x\,u_x+b^2y\,u_y,\] and for analytic $A,B$ the equation
$\cP_{L_0,A,B}[u]=0$ is
\begin{equation}\label{eq:delta}
a^2x^2u_{xx}+2ab\,xy\,u_{xy}+b^2y^2u_{yy}+x\,a\bigl(a+A\bigr)u_x+y\,b\bigl(b+A\bigr)u_y+B\,u=0 .
\end{equation}
This is exactly the form \eqref{eq:regular-form} of Definition~\ref{def:regular}, with principal constants
$a^2,2ab,b^2$; they are real when $a,b\in\R$, as in \cite[Def.~6.1]{LS23}. So $0$ is a \emph{regular singularity},
parabolic in the classical sense ($(2ab)^2=4a^2b^2$). When $A,B$ are constants, \eqref{eq:delta} is an
\emph{Euler PDE} in the sense of \cite[\S3]{LS23}.

For $Q=(q_1,q_2)\in\N^2$ set $\br{Q}:=aq_1+bq_2$ and $\Gamma:=\{\br{Q}:Q\in\N^2\setminus\{0\}\}$.
Write $A_0=A(0)$ and $B_0=B(0)$, and let
\[
\chi(\lambda)=\lambda^2+A_0\lambda+B_0=(\lambda-\lambda_1)(\lambda-\lambda_2),\qquad \mu:=\lambda_2-\lambda_1 .
\]

We recall from \cite[Definition~6.2]{LS23} how the indicial polynomial of an equation in the regular form
\eqref{eq:regular-form} is defined: it is the polynomial
\[
P(r,s)=\alpha r^2+\beta rs+\gamma s^2+\bigl(a_1(0)-\alpha\bigr)r+\bigl(b_1(0)-\gamma\bigr)s+c_1(0),
\]
so that $x^{-r}y^{-s}$ times the Euler part of the operator (its coefficients frozen at $0$), applied to $x^ry^s$, equals $P(r,s)$, and its zero set
$\mathcal C\subset\C^2$ is the \emph{indicial conic}. In \cite[\S3]{LS23} these conics are classified, up to affine
changes of the variables $(r,s)$, into elliptic, hyperbolic and parabolic types and two degenerate types; the
second degenerate type, type~(v), is a pair of distinct parallel lines.

\begin{proposition}[Indicial conic]\label{prop:indicial}
The indicial polynomial of \eqref{eq:delta} is
\[
P(r,s)=\chi(ar+bs).
\]
Hence the indicial conic is the union of two parallel lines (a double line if $\mu=0$),
\[
\mathcal C=\ell_1\cup\ell_2,\qquad \ell_j=\{(r,s): ar+bs=\lambda_j\},
\]
that is, when $\mu\neq0$, the degenerate type~(v) of \cite[\S3]{LS23} rather than the parabolic type $s=r^2$.
Moreover $L_0(x^ry^s)=(ar+bs)\,x^ry^s$, so $v=L_0(x^ry^s)/(x^ry^s)\equiv ar+bs$. Therefore $\mathcal C$ is the
preimage, under $(r,s)\mapsto ar+bs$, of the constant (equilibrium) solutions $\lambda_1,\lambda_2$ of the
Euler--Riccati equation $L_0v=-\chi(v)$.
\end{proposition}

\begin{proof}
With the notation of \cite[Def.~6.2]{LS23}: $A=a^2$, $B=2ab$, $C=b^2$, $a(0,0)=a(a+A_0)$, $b(0,0)=b(b+A_0)$ and
$c(0,0)=B_0$. Therefore
\[P(r,s)=a^2r^2+2abrs+b^2s^2+aA_0\,r+bA_0\,s+B_0=(ar+bs)^2+A_0(ar+bs)+B_0.\]
\end{proof}

\begin{proposition}[General solution of the Euler PDE]\label{prop:euler-general}
Let $A\equiv A_0$ and $B\equiv B_0$ be constants. Let $D\subset(\K^*)^2$ be simply connected (when $\K=\R$, $D$
lies in an open quadrant), fix branches of $\log x$ and $\log y$ on $D$, and let $t=a^{-1}\log x$ and $h=x^{b}y^{-a}$; $h$
is a first integral of $L_0$, and $(t,h)$ are local holomorphic coordinates on $D$. Let $U\subset D$ be an open set on
which $(t,h)$ is a chart onto a product $T\times V$ of connected domains. Then every analytic solution of
\eqref{eq:delta} on $U$ has the form
\[
u=f(h)\,x^{\lambda_1/a}+g(h)\,x^{\lambda_2/a}\quad(\mu\neq0),\qquad
u=\bigl(f(h)+g(h)\log x\bigr)x^{\lambda_1/a}\quad(\mu=0),
\]
with $f,g$ analytic on $V$. In particular, every analytic solution on $D$ has this form locally. The monomial
solutions $x^ry^s$, $(r,s)\in\ell_j$, are the cases $f(h)=h^{\nu}$, with $r=\lambda_j/a+b\nu$ and $s=-a\nu$.
\end{proposition}

\begin{proof}
The Jacobian of $(t,h)$ with respect to $(x,y)$ is $-h/(xy)\neq0$, so $(t,h)$ are local coordinates. In them
$L_0t=1$ and $L_0h=0$, so $L_0=\partial_t$. On $U\cong T\times V$ the equation becomes \[u_{tt}+A_0u_t+B_0u=0\] with
parameter $h\in V$, whose solutions are as stated, with $e^{\lambda t}=x^{\lambda/a}$; since $T$ is connected, the
coefficients are functions of $h\in V$ alone.
\end{proof}

\begin{remark}\label{rem:euler-global}
The formula of Proposition~\ref{prop:euler-general} does not hold globally on every simply connected $D$: the
coefficients are local first integrals, which need not be single-valued functions of $h$ on $h(D)$. For example, let
$a=b=1$, $A_0=-1$, $B_0=0$ (so the roots are $0$ and $1$), $D=(\C\setminus(-\infty,0])^2$ with principal logarithms, and
$u=\operatorname{Log}x-\operatorname{Log}y$. Then $L_0u=0$, so $u$ solves \eqref{eq:delta}; a global formula would force
$u=f(h)$ with $h=x/y$. But $(x,y)=(i,-i)$ and $(-i,i)$ both have $h=-1$, while $u$ takes the values $i\pi$ and $-i\pi$.
All later results in the paper use the proposition only locally.
\end{remark}

As in the classical Frobenius method, a point $(r_0,s_0)$ of the indicial conic is called \emph{resonant} if some
shifted point $(r_0+q_1,s_0+q_2)$, with $(q_1,q_2)\in\N^2\setminus\{0\}$, also lies on the conic
\cite[Definition~6.3]{LS23}. Resonances are the obstruction to solving the recursion for a Frobenius series, and
for square-type foliated equations they can be described completely.

\begin{proposition}[Resonances]\label{prop:resonance}
Let $(r_0,s_0)\in\ell_1$. For $Q\in\N^2$ we have \[P(r_0+q_1,s_0+q_2)=\chi(\lambda_1+\br Q)=\br Q\,(\br Q-\mu).\]
Hence $(r_0,s_0)$ is resonant if and only if
\[
0\in\Gamma\quad\text{or}\quad \mu\in\Gamma .
\]
Resonance therefore depends only on the line, not on the point. More precisely:
\begin{enumerate}
\item $0\in\Gamma$ if and only if $a/b\in\Q_{<0}$. In that case \emph{every} point of $\mathcal C$ is resonant.
\item If $a/b\notin\Q_{<0}$, the line $\ell_1$ (respectively $\ell_2$) is entirely resonant if $\mu\in\Gamma$
(respectively $-\mu\in\Gamma$), and entirely non-resonant otherwise.
\item If $a/b\notin\R_{\le0}$ (Poincar\'e domain), $\Gamma$ lies in an open salient cone, so $\mu$ and $-\mu$ are
not both in $\Gamma$. Hence at least one of the lines $\ell_1,\ell_2$ is non-resonant.
\end{enumerate}
\end{proposition}

\begin{proof}
The identity follows from Proposition~\ref{prop:indicial}. For (1), $aq_1+bq_2=0$ with $Q\in\N^2\setminus\{0\}$
forces $q_1q_2\neq0$ and $a/b=-q_2/q_1$. Item (2) is immediate. For (3), \[\Gamma\subset\{ta+sb:t,s\ge0\}\setminus\{0\},\]
and this cone contains no line through $0$ when $a/b\notin\R_{\le 0}$.
\end{proof}

\begin{remark}[Comparison with \cite{LS23}]\label{rem:LS23}
Two points of comparison with our earlier paper are worth recording. In case (1) the resonant set is the whole of
$\mathcal C$, so the density statement of \cite[Lemma~6.2]{LS23} does not hold for this degenerate class; in its
proof some of the lines $R_Q$ coincide with the whole plane (when $\br Q=0$) or with a component of $\mathcal C$.
Moreover, for \eqref{eq:delta} the ``parabolic of real type'' condition of \cite[Definition~7.1]{LS23} reads
$\mathrm{Re}(ab)>0$, and it therefore fails for real $a,b$ with $ab<0$, which is the saddle case. We take the
opportunity to correct the text following \cite[Def.~7.1]{LS23}, which asserts that this condition is automatic in
the real case: the assertion should be restricted to $ab>0$.
\end{remark}

\section{The Poincar\'e domain: convergence}\label{sec:poincare}

We now show that non-resonant lines give convergent solutions when the eigenvalues $a,b$ lie in the
\emph{Poincar\'e domain}, that is, when $a/b\notin\R_{\le0}$, so that $0$ does not lie in the convex hull of
$\{a,b\}$. In this case the divisors $aq_1+bq_2$ grow linearly in $|Q|$, as in Poincar\'e's linearization theorem,
and the classical majorant argument shows that the Frobenius series converges.

\begin{theorem}[Frobenius solutions in the Poincar\'e domain]\label{thm:poincare}
Let $a,b\in\C^*$ with $a/b\notin\R_{\le0}$. Let $A,B$ be holomorphic in a polydisc $\Delta_R$ around $0$, and let
$\lambda_1,\lambda_2,\mu,\Gamma$ be as above. Assume $\mu\notin\Gamma$ (that is, $\ell_1$ is non-resonant). Then
there is a unique formal power series $F_1=\sum_{Q}d_QX^Q$ with $d_0=1$ such that, for every $(r_0,s_0)\in\ell_1$,
\[
u=x^{r_0}y^{s_0}\,F_1(x,y)
\]
is a formal solution of \eqref{eq:delta}. Moreover $F_1$ converges in a neighbourhood of $0$. The series $F_1$
does not depend on the point $(r_0,s_0)\in\ell_1$. Hence, for every analytic $f$,
$u=f(h)\,x^{\lambda_1/a}F_1$ is a solution, where $h=x^by^{-a}$.
\end{theorem}

\begin{proof}
Write \[A-A_0=\sum_{K\neq0}\alpha_KX^K\text{ and }B-B_0=\sum_{K\neq0}\beta_KX^K.\]Since
\[L_0(x^{r_0}y^{s_0}X^P)=(\lambda_1+\br P)\,x^{r_0}y^{s_0}X^P,\] substituting \[u=x^{r_0}y^{s_0}\sum d_PX^P\]into
\eqref{eq:delta} and comparing coefficients gives the recursion
\begin{equation}\label{eq:recursion}
\br Q\bigl(\br Q-\mu\bigr)\,d_Q=-\sum_{P<Q}\bigl(\alpha_{Q-P}(\lambda_1+\br P)+\beta_{Q-P}\bigr)\,d_P,\qquad Q\neq0 .
\end{equation}
It involves only $\lambda_1$, not $(r_0,s_0)$. Since $a/b\notin\R_{\le0}$, compactness of $[0,1]$ gives
$\delta>0$ with \[|ta+(1-t)b|\ge\delta,\;\text{for all }t\in[0,1],\] so $|\br Q|\ge\delta|Q|$. Hence
\[|\br Q-\mu|\ge\delta|Q|/2\;\text{ for }|Q|\ge 2|\mu|/\delta.\]The finitely many remaining factors are non-zero because
$\mu\notin\Gamma$. Thus there is $c>0$ with \[|\br Q(\br Q-\mu)|\ge c|Q|^2\;\text{ for all }Q\neq0.\]Consequently
\eqref{eq:recursion} determines the $d_Q$ uniquely.

Fix $0<\rho<R$. By Cauchy's estimates \[|\alpha_K|,|\beta_K|\le M\rho^{-|K|}.\] For $P<Q$ we have
\[|\lambda_1+\br P|+1\le C_1|Q|.\]Therefore
\[
|d_Q|\le \frac{C_1M}{c\,|Q|}\sum_{P<Q}\rho^{-|Q-P|}|d_P|\le C_2\sum_{P<Q}\rho^{-|Q-P|}|d_P| .
\]
Define \[e_0=1\text{ and }e_Q=C_2\sum_{P<Q}\rho^{-|Q-P|}e_P.\] By induction $|d_Q|\le e_Q$. The generating function
\[E=\sum e_QX^Q\text{ satisfies }E=1+C_2\Psi E\]with
\[\Psi(x,y)=\sum_{K\neq0}\rho^{-|K|}X^K=\bigl((1-x/\rho)(1-y/\rho)\bigr)^{-1}-1.\]
So $E=1/(1-C_2\Psi)$, which is holomorphic near $0$ because $\Psi(0)=0$. Hence $F_1$ converges. The last
statement follows from linearity and $L_0h=0$.
\end{proof}

\begin{corollary}\label{cor:always}
In the Poincar\'e domain, \eqref{eq:delta} always has a non-trivial family of convergent Frobenius type
solutions. This is the PDE analogue of choosing ``the root with larger real part'' in the classical Frobenius
theorem.
\end{corollary}

\begin{proof}
By Proposition~\ref{prop:resonance}(3), one of the lines is non-resonant. Apply Theorem~\ref{thm:poincare} to it.
\end{proof}

For real $a,b$ the hypothesis $a/b\notin\R_{\le0}$ is $ab>0$, which is the case covered by
\cite[Theorem~7.1]{LS23}. For complex $a,b$, Theorem~\ref{thm:poincare} replaces the condition $\mathrm{Re}(ab)>0$
by the (different) Poincar\'e condition, and it uses only the non-resonance of one line.

\begin{corollary}[Holomorphic solutions of the Riccati PDE]\label{cor:riccati-unique}
Under the hypotheses of Theorem~\ref{thm:poincare}:
\begin{enumerate}
\item Every germ $v\in\cO_0$ solving $L_0v=-v^2-Av-B$ satisfies $v(0)\in\{\lambda_1,\lambda_2\}$.
\item If $\ell_j$ is non-resonant, there is exactly one germ of holomorphic solution with $v(0)=\lambda_j$, namely
$v_j=\lambda_j+L_0F_j/F_j$.
\end{enumerate}
This is a Briot--Bouquet type statement for the Riccati PDE at the equilibria $\lambda_j$.
\end{corollary}

\begin{remark}
The Riccati PDE $L_0v=-v^2-Av-B$ is a singular first order nonlinear PDE of totally characteristic type, and in
the Poincar\'e domain Corollary~\ref{cor:riccati-unique}(2) is close to the Briot--Bouquet type theorems of
G\'erard and Tahara \cite{GT96}. The proof given here is independent of \cite{GT96}: it obtains uniqueness from
the linear Frobenius theory through the substitution $u=x^ry^se^{\varphi}$. The same mechanism is what makes
Theorem~\ref{thm:riccati-saddle} work at rational saddles, where every point of the indicial conic is resonant.
\end{remark}

\begin{proof}
(1) $L_0$ vanishes at $0$, so evaluating the equation at $0$ gives $\chi(v(0))=0$.
(2) Existence: $u=x^{r}y^{s}F_j$ with $(r,s)\in\ell_j$ gives $L_0u/u=\lambda_j+L_0F_j/F_j$, which is holomorphic,
and Proposition~\ref{prop:riccati} applies. Uniqueness: let $v$ be a solution with $v(0)=\lambda_j$. Write
\[v-\lambda_j=\sum_{Q\neq0}\psi_QX^Q\text{ and put }\varphi=\sum_{Q\neq0}(\psi_Q/\br Q)X^Q.\]This series converges
because $|\br Q|\ge\delta|Q|$, and it satisfies $L_0\varphi=v-\lambda_j$. Then $u=x^ry^se^{\varphi}$ satisfies
$L_0u=vu$, so $u$ solves \eqref{eq:delta} by Proposition~\ref{prop:riccati}. By uniqueness in
Theorem~\ref{thm:poincare}, $e^{\varphi}=F_j$, and therefore $v=v_j$.
\end{proof}

\section{Singular points of $L$: regular singularity versus linearization}\label{sec:characterization}

Let $L$ be a germ of holomorphic vector field at $0\in\C^2$ with an isolated singularity at $0$, and let
$A,B\in\cO_0$. Denote by $\cF_L$ the foliation defined by $L$. We say that $\cF_L$ is \emph{analytically
linearizable with non-degenerate diagonal linear part} if there are a germ of biholomorphism $\Phi$ and a unit
$g\in\cO_0$ with $g(0)=1$ such that
\[
\Phi_*(gL)=L_0=ax\,\partial_x+by\,\partial_y,\qquad ab\neq0 .
\]
Since $g(0)=1$ and $L(0)=0$, the numbers $a,b$ are then the eigenvalues of $DL(0)$.

In \cite[Def.~6.1]{LS23} a singular point is called regular if, after a change of coordinates and a division by
coefficients, the equation takes an Euler-type form whose principal coefficients are \emph{real} constants. Since
we work with holomorphic equations, we use the following complex version. It makes explicit the class of
coordinate changes and of divisions that are allowed.

\begin{definition}[Regular singularity, complex constants]\label{def:regular}
Let $\cP$ be a second order linear differential operator with holomorphic coefficients near $0\in\C^2$. We say that
$\cP$ has a \emph{regular singularity at $0$} if there are a germ of biholomorphism $\Psi$ fixing $0$ and a
holomorphic germ $f\not\equiv0$ such that $f^{-1}\,\Psi_*\cP$ equals
\begin{equation}\label{eq:regular-form}
\alpha x^2\partial_x^2+\beta xy\,\partial_x\partial_y+\gamma y^2\partial_y^2+x\,a_1(x,y)\,\partial_x+y\,b_1(x,y)\,\partial_y+c_1(x,y),
\end{equation}
with constants $\alpha,\beta,\gamma\in\C$ and holomorphic $a_1,b_1,c_1$. Here $\Psi_*\cP$ denotes the operator
$w\mapsto\cP[w\circ\Psi]\circ\Psi^{-1}$. If, in addition, $\Psi$ and $f$ can be chosen with
$\alpha,\beta,\gamma\in\R$, we say that the regular singularity has \emph{real principal constants}; this is the
normalization of \cite[Def.~6.1]{LS23}.
\end{definition}

Our definition is the complex-analytic counterpart of \cite[Def.~6.1]{LS23}. The latter explicitly
requires real principal constants, whereas \cite[Theorem~7.1 and \S8]{LS23} also studies
complex analytic Euler equations; we do not identify the two definitions verbatim, nor do we assume
that $f$ is a unit; for square-type foliated equations, the proof below shows that this is automatic. The
distinction between real and complex principal constants is not merely formal; see
Corollary~\ref{cor:real-constants}.

\begin{theorem}[Characterization]\label{thm:characterization}
Let $L$, $A$, $B$ be as above. The following are equivalent:
\begin{enumerate}
\item[(a)] \eqref{eq:main} has a regular singularity at $0$ in the sense of Definition~\ref{def:regular};
\item[(b)] $\cF_L$ is analytically linearizable with non-degenerate diagonal linear part.
\end{enumerate}
In this case, if $a,b$ are the eigenvalues of $DL(0)$ and $\lambda_1,\lambda_2$ are the roots of
$\lambda^2+A(0)\lambda+B(0)$, the indicial conic is
\[
\mathcal C=\{(r,s):ar+bs=\lambda_1\}\cup\{(r,s):ar+bs=\lambda_2\}.
\]
It depends on the linearizing coordinates only through the order of the eigenvalues: exchanging the two
eigendirections replaces $(a,b,r,s)$ by $(b,a,s,r)$.
\end{theorem}

\begin{proof}
(b)$\Rightarrow$(a). Put $M=gL$. Since $L=g^{-1}M$, Lemma~\ref{lem:invariance}(2), applied to the unit $g^{-1}$
and the field $M$, gives
\[
\cP_{L,A,B}=g^{-2}\,\cP_{M,A',B'},\qquad A'=gA+g\,M(g^{-1})=gA-L(g),\qquad B'=g^2B .
\]
By Lemma~\ref{lem:invariance}(1),
\[\Phi_*\cP_{L,A,B}=(g^{-2}\circ\Phi^{-1})\,\cP_{L_0,\tilde A,\tilde B},\]with $\tilde A=A'\circ\Phi^{-1}$ and
$\tilde B=B'\circ\Phi^{-1}$. Dividing by the unit $g^{-2}\circ\Phi^{-1}$ we obtain \eqref{eq:delta}, which is of
the form \eqref{eq:regular-form} with $\alpha=a^2$, $\beta=2ab$, $\gamma=b^2$. Moreover $g(0)=1$, and $L(g)(0)=0$
because $L(0)=0$. Hence $\tilde A(0)=A(0)$ and $\tilde B(0)=B(0)$, and Proposition~\ref{prop:indicial} gives
$\mathcal C$. Finally \[D(gL)(0)=g(0)\,DL(0)=DL(0),\]and $\Phi$ preserves eigenvalues, so $a,b$ are the eigenvalues
of $DL(0)$.

(a)$\Rightarrow$(b). Let $\Psi$, $f$ be as above and write \[M=\Psi_*L=m_1\partial_x+m_2\partial_y.\] Then $M$ has an
isolated zero at $0$, and \[\Psi_*\cP_{L,A,B}=\cP_{M,A\circ\Psi^{-1},B\circ\Psi^{-1}}\]by
Lemma~\ref{lem:invariance}(1). By \eqref{eq:main-expanded}, its second order part is
\[m_1^2\partial_x^2+2m_1m_2\,\partial_x\partial_y+m_2^2\partial_y^2.\] Comparing with \eqref{eq:regular-form} gives
\begin{equation}\label{eq:symbol}
m_1^2=\alpha\,x^2f,\qquad 2m_1m_2=\beta\,xy\,f,\qquad m_2^2=\gamma\,y^2f
\end{equation}
in $\cO_0$.

\emph{Step 1: $\alpha\gamma\neq0$.} If $\alpha=0$ then $m_1\equiv0$, so the zero set of $M$ is $\{m_2=0\}$. This is
a curve through $0$ (note $m_2(0)=0$ because $M(0)=0$), contradicting isolatedness. Likewise $\gamma\neq0$.

\emph{Step 2: $M$ is a holomorphic multiple of a diagonal linear field.} $\cO_0$ is a unique factorization domain
and $x$ is prime. Since $x^2\mid m_1^2$, we get $x\mid m_1$; write $m_1=\sqrt\alpha\,x\,\varphi$ for a fixed square
root $\sqrt\alpha$. Then \eqref{eq:symbol} gives $f=\varphi^2$. Similarly $m_2=\sqrt\gamma\,y\,\psi$ with
$\psi^2=f=\varphi^2$. Since $\cO_0$ is an integral domain, $\psi=\pm\varphi$; changing the sign of $\sqrt\gamma$
if necessary, $\psi=\varphi$. Therefore
\[
M=\varphi\,\bigl(\sqrt\alpha\,x\partial_x+\sqrt\gamma\,y\partial_y\bigr),
\]
and the middle equation of \eqref{eq:symbol} forces $\beta=2\sqrt\alpha\sqrt\gamma$.

\emph{Step 3: $\varphi$ is a unit.} If $\varphi(0)=0$, the zero set of $M$ contains the curve $\{\varphi=0\}$,
contradicting isolatedness. Hence $\varphi(0)\neq0$; in particular $f=\varphi^2$ is a unit.

\emph{Conclusion.} Put $g=\varphi(0)\,(\varphi^{-1}\circ\Psi)$, a unit with $g(0)=1$, and
$a=\varphi(0)\sqrt\alpha$, $b=\varphi(0)\sqrt\gamma$. Then
\[\Psi_*(gL)=(g\circ\Psi^{-1})\,M=a\,x\partial_x+b\,y\partial_y\;\text{ with }ab\neq0.\]This is (b) with $\Phi=\Psi$.

\emph{Independence of the regular form.} Finally, let \eqref{eq:regular-form} be \emph{any} regular form of
\eqref{eq:main}, obtained with some $\Psi$ and $f$. By Steps~1--3, $f=\varphi^2$ with $\varphi$ a unit and
$\Psi_*L=\varphi M_0$, where $M_0=\sqrt\alpha\,x\partial_x+\sqrt\gamma\,y\partial_y$. By Lemma~\ref{lem:invariance},
\[
f^{-1}\Psi_*\cP_{L,A,B}=\varphi^{-2}\,\cP_{\varphi M_0,\,A\circ\Psi^{-1},\,B\circ\Psi^{-1}}=\cP_{M_0,\hat A,\hat B},\qquad
\hat A=\frac{A\circ\Psi^{-1}+M_0(\varphi)}{\varphi},\quad \hat B=\frac{B\circ\Psi^{-1}}{\varphi^2}.
\]
Since $M_0(\varphi)(0)=0$ and $\Psi(0)=0$, we get $\hat A(0)=A(0)/\varphi(0)$ and $\hat B(0)=B(0)/\varphi(0)^2$. By
Proposition~\ref{prop:indicial} applied to $M_0$, the indicial polynomial of this regular form is
$(\sqrt\alpha\,r+\sqrt\gamma\,s)^2+\hat A(0)(\sqrt\alpha\,r+\sqrt\gamma\,s)+\hat B(0)$, whose zero set is
\[
\{\varphi(0)\sqrt\alpha\,r+\varphi(0)\sqrt\gamma\,s=\lambda_j\},\qquad j=1,2.
\]
Since $D(\Psi_*L)(0)=\varphi(0)\operatorname{diag}(\sqrt\alpha,\sqrt\gamma)$ is conjugate to $DL(0)$, we have
$(\varphi(0)\sqrt\alpha,\varphi(0)\sqrt\gamma)\in\{(a,b),(b,a)\}$. Hence every regular form has the indicial conic
stated in the theorem, up to the exchange $(a,b,r,s)\mapsto(b,a,s,r)$.
\end{proof}

The proof of (a)$\Rightarrow$(b) uses only the second order part of \eqref{eq:regular-form}. The conditions on
the first order part (the factors $x$ and $y$ in front of $\partial_x$ and $\partial_y$) are then automatically
satisfied, by the computation in (b)$\Rightarrow$(a). Thus, for square-type foliated equations, \emph{regularity is
decided by the principal symbol alone}. The hypothesis that the singularity is isolated is used in Steps~1 and~3.
For \[L=xy(a\partial_x+b\partial_y),\]$\{L=0\}$ is the union of the axes, so
Theorem~\ref{thm:characterization} does not apply. That case needs a separate (non-isolated) treatment.

The proof of Theorem~\ref{thm:characterization} also determines which regular forms are possible. As a result,
the real normalization of \cite[Def.~6.1]{LS23} is a genuinely stronger condition.

\begin{corollary}[Real principal constants]\label{cor:real-constants}
Let $L$, $A$, $B$ be as in Theorem~\ref{thm:characterization}. Then \eqref{eq:main} has a regular singularity at $0$
with real principal constants if and only if $\cF_L$ is analytically linearizable with non-degenerate diagonal
linear part \emph{and} the ratio $a/b$ of the eigenvalues of $DL(0)$ is real. Moreover, in any regular form
\eqref{eq:regular-form} of \eqref{eq:main} we have $\alpha\gamma\neq0$, $\beta=2\sqrt\alpha\sqrt\gamma$ for suitable
square roots, and
\[
\frac{\sqrt\alpha}{\sqrt\gamma}\in\Bigl\{\frac ab,\ \frac ba\Bigr\}.
\]
\end{corollary}

\begin{proof}
Let \eqref{eq:regular-form} be any regular form, obtained with $\Psi$ and $f$. Steps~1--3 of the proof of
(a)$\Rightarrow$(b) give $\alpha\gamma\neq0$, $\beta=2\sqrt\alpha\sqrt\gamma$ and
\[\Psi_*L=\varphi\,(\sqrt\alpha\,x\partial_x+\sqrt\gamma\,y\partial_y)\]with $\varphi(0)\neq0$. Hence
\[D(\Psi_*L)(0)=\varphi(0)\operatorname{diag}(\sqrt\alpha,\sqrt\gamma).\]This matrix is conjugate to $DL(0)$, so
\[\{a,b\}=\varphi(0)\{\sqrt\alpha,\sqrt\gamma\},\]which proves the last assertion. If $\alpha,\beta,\gamma$ are real,
then $\sqrt\alpha/\sqrt\gamma=\beta/(2\gamma)\in\R$, so $a/b\in\R$. Conversely, assume (b) of
Theorem~\ref{thm:characterization} and $a/b\in\R$. The proof of (b)$\Rightarrow$(a) gives a regular form with
\[(\alpha,\beta,\gamma)=(a^2,2ab,b^2).\] Dividing further by the constant $b^2$ gives the real constants
$((a/b)^2,\,2a/b,\,1)$.
\end{proof}

\begin{example}\label{ex:complex-constants}
Let $L=(1+i)\,x\partial_x+y\partial_y$, with arbitrary $A,B\in\cO_0$. Then \eqref{eq:main} is already of the form
\eqref{eq:regular-form} with $(\alpha,\beta,\gamma)=(2i,\,2+2i,\,1)$, so it has a regular singularity in the sense
of Definition~\ref{def:regular}. However, $a/b=1+i\notin\R$. By Corollary~\ref{cor:real-constants}, \emph{no}
biholomorphic change of coordinates and \emph{no} division by a holomorphic germ brings it to a regular form with
real principal constants. When $L$ is real analytic with real eigenvalues, as in the real setting of
\cite{LS23}, the two notions coincide.
\end{example}

By Theorem~\ref{thm:characterization}, deciding whether a square-type foliated equation has a regular singularity
is the same as deciding whether the foliation $\cF_L$ is analytically linearizable, a question studied extensively
since Poincar\'e. We recall the standard terminology; see \cite{IY} for a complete treatment and
\cite[Chapter~1]{Brunella2015} for a short introduction. Let $a,b$ be the eigenvalues of $DL(0)$. If $ab\neq0$, the
singularity is in the \emph{Poincar\'e domain} when $a/b\notin\R_{\le0}$; within that domain it is a \emph{node} when
$a/b>0$, and it is \emph{resonant} when $a/b$ or $b/a$ is an integer $n\ge2$. It is a \emph{saddle} when $a/b<0$. A
saddle has two smooth invariant curves through $0$, its \emph{separatrices}. The \emph{holonomy} of a separatrix is the
germ of diffeomorphism of a transversal disc obtained by following the leaves once around the separatrix; its
multiplier is $e^{2\pi i a/b}$. If exactly one eigenvalue vanishes, the singularity is a \emph{saddle-node}; if both
vanish but $DL(0)\neq0$, the linear part is \emph{nilpotent}. Finally, a real number $\alpha$ satisfies the
\emph{Brjuno condition} if $\sum_k(\log q_{k+1})/q_k<\infty$, where $p_k/q_k$ are the convergents of its continued
fraction. The following corollary translates the main results of linearization theory, case by case, into
statements about regular singularities.

\begin{corollary}[Dictionary]\label{cor:dictionary}
Let $L$ be a holomorphic vector field with an isolated singularity at $0$, let $a,b$ be the eigenvalues of $DL(0)$,
and let $A,B$ be holomorphic near $0$.
\begin{enumerate}
\item If $DL(0)$ is diagonalizable and the singularity is a non-resonant point of the Poincar\'e domain, then
equation~\eqref{eq:main} always has a regular singularity at $0$. In particular Theorem~\ref{thm:poincare} and
Corollaries~\ref{cor:always} and~\ref{cor:riccati-unique} apply.
\item If the singularity is a resonant node, with $a/b$ or $b/a$ equal to an integer $n\ge2$, then $\cF_L$ is
analytically equivalent to the foliation of $x\partial_x+(ny+\varepsilon x^n)\partial_y$ with $\varepsilon\in\{0,1\}$,
and equation~\eqref{eq:main} has a regular singularity at $0$ if and only if $\varepsilon=0$.
\item If the singularity is a saddle with irrational ratio $a/b$, and $-a/b$ satisfies the Brjuno condition, then
equation~\eqref{eq:main} has a regular singularity at $0$. If $-a/b$ does not satisfy it, there are vector fields $L$ with
these eigenvalues for which equation~\eqref{eq:main} has no regular singularity at $0$.
\item If the singularity is a saddle with rational ratio $a/b=-p/q$, then equation~\eqref{eq:main} has a regular singularity
at $0$ if and only if the holonomy of one (equivalently, of either) separatrix is linearizable, that is, periodic.
\item If the singularity is a saddle-node, or has nilpotent linear part, or $DL(0)=0$, or $DL(0)$ is not
diagonalizable, then equation~\eqref{eq:main} never has a regular singularity at $0$.
\end{enumerate}
\end{corollary}

\begin{proof}[Proof of Corollary~\ref{cor:dictionary}]
Each item follows from Theorem~\ref{thm:characterization} and a classical linearization result: Poincar\'e's
theorem \cite{Poincare} for (1), the Poincar\'e--Dulac normal form for (2), Brjuno's theorem \cite{Brjuno} and its
converse by Yoccoz \cite{Yoccoz}, P\'erez-Marco and Yoccoz \cite{PMY} and Mattei and Moussu \cite{MM} for (3), and
Mattei and Moussu \cite{MM} for (4); see \cite{IY} for all of these. Item~(5) holds because, by
Theorem~\ref{thm:characterization}, regularity forces a non-degenerate diagonal linear part.
\end{proof}

In item~(3) of Corollary~\ref{cor:dictionary}, the indicial conic is defined, but convergence of the Frobenius series is
a small-divisor problem, which we do not address here.

\section{Rational saddles}\label{sec:siegel}

Assume that $\cF_L$ is linearizable, so that we work with \eqref{eq:delta}, and that $a/b\in\R_{<0}$. The
irrational case $a/b\notin\Q$ is a small-divisor problem and is not treated here. We treat the
rational case completely when $\mu\notin\Z$, and the integral case under a one-sided hypothesis.

\subsection{Rational saddles}
Let $L_0=p\,x\partial_x-q\,y\partial_y$ with $p,q\in\N^*$ coprime, and let $h=x^qy^p$, a holomorphic first integral
($L_0h=0$). Here $\Gamma=\{pq_1-qq_2\}\ni0$, so \emph{every} point of $\mathcal C$ is resonant
(Proposition~\ref{prop:resonance}(1)). The resonant monomials are exactly the powers $h^n$. Every formal series
splits as \[F=\sum_{k\in\Z}F_k\]into $L_0$-weight components ($L_0F_k=kF_k$), and $F_0\in\C[[h]]$.

\subsubsection*{Leaves, annuli and leafwise monodromy}
All our solutions are defined near $0$, so the relevant leaves are \emph{local} ones. For $\rho>0$ let
$\Delta_\rho=\{|x|<\rho,\ |y|<\rho\}$, and for $h_0\neq0$ let $C_{h_0}=\{x^qy^p=h_0\}$.

\begin{lemma}\label{lem:annulus}
Let $h_0\neq0$ and fix $\zeta$ with $\zeta^p=h_0$.
\begin{enumerate}
\item The map $\sigma(s)=(s^p,\zeta s^{-q})$ is a biholomorphism from $\C^*$ onto the smooth curve $C_{h_0}$.
\item If $0<|h_0|<\rho^{p+q}$, then $C_{h_0}\cap\Delta_\rho=\sigma(\mathcal A)$, where
$\mathcal A=\{r_-<|s|<r_+\}$ is the non-empty annulus with $r_+=\rho^{1/p}$ and $r_-=(|h_0|^{1/p}/\rho)^{1/q}$. The
loop $\gamma(\theta)=\sigma(re^{2\pi i\theta})$, $\theta\in[0,1]$, with $r_-<r<r_+$, generates
\[\pi_1(C_{h_0}\cap\Delta_\rho)\cong\Z.\]Along $\gamma$, $x$ turns $p$ times and $y$ turns $-q$ times around $0$; in
particular any branch of $\log x$ increases by $2\pi i p$.
\item $\sigma_*(s\partial_s)=L_0$. Hence, if $A,B$ are holomorphic on $\Delta_\rho$, then for every $u$ defined on an
open subset of $\Delta_\rho\setminus\{xy=0\}$ we have \[\cP_{L_0,A,B}[u]\circ\sigma=\mathcal L_{h_0}[u\circ\sigma],\]where
\[
\mathcal L_{h_0}[w]=(s\partial_s)^2w+A(\sigma(s))\,(s\partial_s)w+B(\sigma(s))\,w
\]
is a linear ODE with holomorphic coefficients and \emph{no singular points} on $\mathcal A$.
\end{enumerate}
\end{lemma}

\begin{proof}
(1) $\sigma(s)\in C_{h_0}$ because $s^{pq}\zeta^ps^{-pq}=h_0$. If $\sigma(s_1)=\sigma(s_2)$, then $\omega=s_1/s_2$
satisfies \[\omega^p=\omega^q=1,\]so $\omega=1$ since $\gcd(p,q)=1$. Given $(x,y)\in C_{h_0}$, choose $s_0$ with
$s_0^p=x$. Then \[y^p=h_0x^{-q}=(\zeta s_0^{-q})^p,\]so $y=\varepsilon\zeta s_0^{-q}$ with $\varepsilon^p=1$. Since $q$
is invertible modulo $p$, there is a $p$-th root of unity $\omega$ with $\omega^{-q}=\varepsilon$, and then
$\sigma(\omega s_0)=(x,y)$. The curve is smooth because \[d(x^qy^p)=x^{q-1}y^{p-1}(qy\,dx+px\,dy)\neq0\]on it. A
holomorphic bijection between Riemann surfaces is biholomorphic.
(2) $|s^p|<\rho$ if and only if $|s|<r_+$, and $|\zeta s^{-q}|<\rho$ if and only if $|s|>r_-$. Moreover $r_-<r_+$
if and only if $|h_0|<\rho^{p+q}$. The winding numbers of $s^p$ and $s^{-q}$ along $|s|=r$ are $p$ and $-q$.
(3) $s\,\tfrac{d}{ds}(s^p)=p\,s^p$ and $s\,\tfrac{d}{ds}(\zeta s^{-q})=-q\,\zeta s^{-q}$, so $\sigma_*(s\partial_s)=L_0$.
Since $L_0$ is tangent to $C_{h_0}$, $L_0^ku$ restricted to the leaf depends only on $u$ restricted to the leaf.
\end{proof}

\begin{definition}\label{def:monodromy}
Let $A,B$ be holomorphic on $\Delta_\rho$ and $0<|h_0|<\rho^{p+q}$. The local solutions of
$\mathcal L_{h_0}[w]=0$ near a base point of $\mathcal A$ form a two-dimensional space $\mathcal S_{h_0}$. Since
$\mathcal L_{h_0}$ has no singular points on $\mathcal A$, every element of $\mathcal S_{h_0}$ continues analytically
along every path in $\mathcal A$. Continuation along the generator $\gamma$ defines the \emph{leafwise monodromy}
$M_{h_0}\in GL(\mathcal S_{h_0})$, well defined up to conjugacy. Since $\gamma$ generates
$\pi_1(C_{h_0}\cap\Delta_\rho)$, a solution germ extends to a \emph{single-valued} solution on the local leaf
$C_{h_0}\cap\Delta_\rho$ if and only if it is fixed by $M_{h_0}$. We say that the family is \emph{isomonodromic for an
eigenvalue branch} $e^{2\pi i\lambda(h_0)}$ of $M_{h_0}$ if this eigenvalue does not depend on $h_0$.
\end{definition}

The definition concerns the regular leaf $h=h_0\neq0$ and its annular domain, not the singular fiber $h=0$. The two notions should not be confused. Analytic continuation along $\gamma$ is always possible, by
Lemma~\ref{lem:annulus}(3). A non-trivial single-valued solution on the local leaf exists only when $M_{h_0}$ has the
eigenvalue $1$. The solutions constructed below are, in general, \emph{not} single-valued on the leaves: they are
multivalued, and $M_{h_0}$ multiplies them by explicit factors. To write them down we use the universal cover of the
annulus, \[\tilde{\mathcal A}=\{\tau\in\C:\ \log r_-<\operatorname{Re}\tau<\log r_+\},\]with covering map $\tau\mapsto
s=e^\tau$. On $\tilde{\mathcal A}$ the function $\tau$ is a global branch of $\log s$. The loop $\gamma$ lifts to a path
from $\tau_0$ to $\tau_0+2\pi i$, and the deck transformation $\tau\mapsto\tau+2\pi i$ corresponds to $\gamma$.

The simplest situation is when the coefficients are constant along the leaves, that is, $A,B\in\C\{h\}$. For each
$h$ put \[\chi_h(\lambda)=\lambda^2+A(h)\lambda+B(h).\]If $\mu=\lambda_2-\lambda_1\neq0$, the roots
$\lambda_1(h),\lambda_2(h)$ of $\chi_h$ are holomorphic near $h=0$, with $\lambda_j(0)=\lambda_j$.

\begin{proposition}[Leafwise exponents]\label{prop:saddle}
Let $A=A(h)$, $B=B(h)$ be holomorphic functions of $h=x^qy^p$, and assume $\mu\neq0$.
\begin{enumerate}
\item For every holomorphic $f$ of one variable, $u=f(h)\,x^{\lambda_j(h)/p}$ ($j=1,2$) solves \eqref{eq:delta}.
Locally on $(\C^*)^2$ near $0$, every solution is a sum of two such functions.
\item There is a (formal or convergent) Frobenius type solution with initial point on $\ell_j$ if and only if
$\lambda_j(h)\equiv\lambda_j$ is constant, equivalently $\chi_h(\lambda_j)\equiv0$. In that case
$u=x^ry^s$, $(r,s)\in\ell_j$, are exact solutions.
\item If $\lambda_j(h)$ is not constant, then
\[
x^{\lambda_j(h)/p}=x^{\lambda_j/p}\sum_{n\ge0}\frac1{n!}\Bigl(\frac{\lambda_j(h)-\lambda_j}{p}\Bigr)^n(\log x)^n,
\]
so, when expanded around the constant exponent $\lambda_j$, the solutions contain \emph{all} powers $(\log x)^n$,
each paired with $h^{n}$ or a higher power. In the classical Frobenius theory the power of the logarithm is
bounded by the multiplicity of the resonance. The unbounded powers here are produced by the expansion and do
not appear in the closed form $x^{\lambda_j(h)/p}$. The same phenomenon occurs for a linear ODE whose exponent
$\lambda(\varepsilon)$ depends on a parameter, when $x^{\lambda(\varepsilon)}$ is expanded in powers of $\varepsilon$.
What is specific to the PDE setting is that the role of the parameter is played by the first integral $h$.
\end{enumerate}
\end{proposition}

\begin{proof}
(1) $L_0(p^{-1}\log x)=1$ and $L_0h=0$. Hence \[L_0\bigl(f(h)x^{\lambda(h)/p}\bigr)=\lambda(h)\,f(h)x^{\lambda(h)/p},\]
and therefore \[\cP_{L_0,A,B}[f(h)x^{\lambda(h)/p}]=\chi_h(\lambda(h))\,f(h)x^{\lambda(h)/p}.\] This vanishes for
$\lambda=\lambda_j$. Locally, \[(t,h)=(p^{-1}\log x,h)\]are coordinates with $L_0=\partial_t$, and the equation becomes
\[u_{tt}+A(h)u_t+B(h)u=0,\] whose solutions are \[f(h)e^{\lambda_1(h)t}+g(h)e^{\lambda_2(h)t}.\]

(2) Let $u=x^ry^sF$, where $(r,s)\in\ell_j$, and let $F=\sum_kF_k$ be the weight decomposition of $F$. Then
\[L_0(x^ry^sF_k)=(\lambda_j+k)\,x^ry^sF_k.\]Since $A$ and $B$ have weight $0$, we get
\[\cP[x^ry^sF_k]=\chi_h(\lambda_j+k)\,x^ry^sF_k,\]and these terms have distinct weights. Hence
$\cP[u]=0$ if and only if $\chi_h(\lambda_j+k)F_k=0$ for every $k$. Since $F_0(0)=d_0=1$, this forces $\chi_h(\lambda_j)\equiv0$ in
$\C[[h]]$. By continuity and $\mu\neq0$, this means $\lambda_j(h)\equiv\lambda_j$. Conversely, if
$\chi_h(\lambda_j)\equiv0$ then $F\equiv1$ works.

(3) Expand $\exp\bigl(p^{-1}(\lambda_j(h)-\lambda_j)\log x\bigr)$. Since $\lambda_j(h)-\lambda_j=O(h)$, its $n$-th
power is $O(h^n)$.
\end{proof}

\begin{corollary}[Monodromy and isomonodromy]\label{cor:monodromy}
In the setting of Proposition~\ref{prop:saddle}, a Frobenius type solution on $\ell_j$ exists if and only if the family
of leafwise equations is isomonodromic for the eigenvalue branch $e^{2\pi i\lambda_j(h_0)}$, that is, if and only if
this eigenvalue of the monodromy of $\mathcal L_{h_0}$ does not depend on $h_0$.
\end{corollary}

\begin{proof}
In the coordinate $s$ of Lemma~\ref{lem:annulus}, $A$ and $B$ are constant on each leaf and
\[\mathcal L_{h_0}[w]=(s\partial_s)^2w+A(h_0)\,s\partial_sw+B(h_0)\,w,\]
which is defined on the whole curve $C_{h_0}\cong\C^*$. On the universal cover, where $s\partial_s=\partial_\tau$, a basis
of solutions is $e^{\lambda_1(h_0)\tau}$, $e^{\lambda_2(h_0)\tau}$. Continuation along a generator of $\pi_1(C_{h_0})$
corresponds to $\tau\mapsto\tau+2\pi i$, so it multiplies $e^{\lambda_j(h_0)\tau}$ by $e^{2\pi i\lambda_j(h_0)}$. Since
$\lambda_j$ is continuous on a punctured disc, $e^{2\pi i\lambda_j}$ is constant exactly when $\lambda_j$ is constant,
and by Proposition~\ref{prop:saddle}(2) this is exactly when a Frobenius type solution on $\ell_j$ exists.
\end{proof}

This is the PDE counterpart of the resonance obstruction, and it ties the Frobenius theory to the monodromy point
of view of foliation theory. When $\mu\notin\Z$ and $h_0$ is small, the two eigenvalues are distinct, since
$\lambda_1(h_0)-\lambda_2(h_0)$ is close to $-\mu$.

For general $A,B$, not constant along the leaves, Proposition~\ref{prop:saddle} suggests an exponent that depends
on the first integral, and this is indeed the right ansatz. Theorem~\ref{thm:rational} below contains
Proposition~\ref{prop:saddle} as the special case $A,B\in\C\{h\}$, in which the recursion has no source of non-zero
weight, so that $G=1$ and $\lambda(h)=\lambda_1(h)$.

\begin{theorem}[Frobenius method over $\C\{h\}$]\label{thm:rational}
Let $L_0=p\,x\partial_x-q\,y\partial_y$ with $p,q\in\N^*$ coprime, and $h=x^qy^p$. Let $A,B$ be holomorphic in a
polydisc $\Delta_R$ around $0$, let $\lambda_1,\lambda_2$ be the roots of $\chi(\nu)=\nu^2+A(0)\nu+B(0)$, and
assume
\[
\mu=\lambda_2-\lambda_1\notin\Z .
\]
Then there exist a unique formal series $\lambda(h)\in\C[[h]]$ with $\lambda(0)=\lambda_1$ and a unique
$G\in\C[[x,y]]$ with $G(0)=1$ and no monomials $h^n$, $n\ge1$ (i.e.\ the weight-$0$ part of $G$ is $1$), such that
\[
u=x^{\lambda(h)/p}\,G(x,y)
\]
is a formal solution of \eqref{eq:delta}, meaning that $\cP_{L_0,A,B}[u]=x^{\lambda(h)/p}\,\mathcal N_\lambda[G]$ with
$\mathcal N_\lambda[G]=0$ in $\C[[x,y]]$ (see Step~1 of the proof). Moreover $\lambda\in\C\{h\}$ and $G\in\C\{x,y\}$ converge, so $u$ is a
genuine solution on every simply connected domain $D\subset\Delta_r\setminus\{x=0\}$, with a branch of $\log x$
fixed on $D$. Here $r>0$ is such that $G$ converges on $\Delta_r$ and $\lambda$ converges for $|h|<r^{p+q}$ (both hold
for $r=\rho_0$ in the proof below); since $|h|<r^{p+q}$ on $\Delta_r$, $u$ is defined on $D$.
\end{theorem}

\begin{proof}[Proof of Theorem~\ref{thm:rational}]
In the notation of the weight decomposition we write \[G=\sum_{k\in\Z}e_k(h)\,\mathfrak m_k\text{ with }e_0\equiv1,\]where
$\mathfrak m_k$ is the monomial of weight $k$ not divisible by $h$.

Write $\br Q=pq_1-qq_2\in\Z$. Since $\gcd(p,q)=1$, $\br Q=0$ if and only if $X^Q=h^n$ for some $n$.
Recall that $\chi(\lambda_1)=0$ and \[\chi'(\lambda_1)=2\lambda_1+A(0)=\lambda_1-\lambda_2=-\mu.\]

\emph{Step 1: reduction.} Since $L_0h=0$ and $L_0(\log x)=p$, we have \[L_0\bigl(x^{\lambda(h)/p}\bigr)=\lambda(h)\,x^{\lambda(h)/p}\]
for every $\lambda\in\C[[h]]$. Hence
\[
\cP_{L_0,A,B}\bigl[x^{\lambda(h)/p}G\bigr]=x^{\lambda(h)/p}\,\mathcal N_\lambda[G],\qquad
\mathcal N_\lambda[G]=L_0^2G+(2\lambda+A)\,L_0G+(\lambda^2+A\lambda+B)\,G .
\]

\emph{Step 2: a single unknown.} Write $A=A(0)+\mathfrak a$, $B=B(0)+\mathfrak b$ with $\mathfrak a(0)=\mathfrak b(0)=0$, $\lambda=\lambda_1+\Lambda$ with
$\Lambda\in h\,\C[[h]]$, and $G=1+\hat G$, where $\hat G$ has no monomials $h^n$ ($n\ge0$). Using
\[\chi(\lambda_1+\Lambda)=-\mu\Lambda+\Lambda^2\;\text{ and }\;2\lambda+A=-\mu+2\Lambda+\mathfrak a,\]we get
\[
\mathcal N_\lambda[G]=(L_0^2-\mu L_0)\hat G-\mu\Lambda+\Bigl\{2\Lambda L_0\hat G+\mathfrak a\,L_0\hat G-\mu\Lambda\hat G
+\bigl(\Lambda^2+\lambda_1\mathfrak a+\mathfrak a\Lambda+\mathfrak b\bigr)(1+\hat G)\Bigr\}.
\]
The monomials of $\Lambda$ are exactly those with $\br Q=0$ (and $Q\neq0$), and the monomials of $\hat G$ are
exactly those with $\br Q\neq0$. So we may combine both unknowns into one series
\[
W:=\Lambda+\hat G=\sum_{Q\neq0}w_QX^Q .
\]
Let $\mathcal D$ be the diagonal operator $\mathcal D X^Q=\delta_QX^Q$, with
\[
\delta_Q=\br Q\bigl(\br Q-\mu\bigr)\ \text{ if }\br Q\neq0,\qquad \delta_Q=-\mu\ \text{ if }\br Q=0 .
\]
Then \[(L_0^2-\mu L_0)\hat G-\mu\Lambda=\mathcal DW,\] and the equation $\mathcal N_\lambda[G]=0$ becomes
\begin{equation}\label{eq:fixed-point}
\mathcal DW=\mathcal R(W):=-\Bigl\{2\Lambda L_0\hat G+\mathfrak a\,L_0\hat G-\mu\Lambda\hat G+\bigl(\Lambda^2+\lambda_1\mathfrak a+\mathfrak a\Lambda+\mathfrak b\bigr)(1+\hat G)\Bigr\}.
\end{equation}

\emph{Step 3: formal existence and uniqueness.} Since $\mu\notin\Z$ and $\br Q\in\Z$, we have $\delta_Q\neq0$ for all
$Q\neq0$. Every term of $\mathcal R(W)$ is either independent of $W$ ($\lambda_1\mathfrak a+\mathfrak b$) or a product of a $W$-dependent
factor with another factor vanishing at $0$ ($\Lambda$, $\hat G$, $\mathfrak a$ or $\mathfrak b$). Hence the coefficient
$\mathcal R(W)_Q$ depends only on the $w_P$ with $|P|<|Q|$. The relation $w_Q=\delta_Q^{-1}\mathcal R(W)_Q$
therefore determines $W$ uniquely by induction on $|Q|$. This gives the unique $\lambda=\lambda_1+\Lambda$ and
$G=1+\hat G$.

\emph{Step 4: lower bound for the divisors.} Let $d=\operatorname{dist}(\mu,\Z)>0$ and $c=d/2$. For $k\in\Z\setminus\{0\}$,
\[|k(k-\mu)|\ge|k|\,d\ge c(1+|k|),\text{ and }|\mu|\ge d\ge c.\] Thus
\begin{equation}\label{eq:divisor-bound}
|\delta_Q|\ge c\,\bigl(1+|\br Q|\bigr)\qquad(Q\neq0).
\end{equation}
This lower bound absorbs the derivative $L_0$ in $\mathcal R$. In the term $\Lambda L_0\hat G$ the coefficient at
$Q$ is $\sum_{K+P=Q}\lambda_K\br P\hat g_P$ with $\br K=0$, so that $\br P=\br Q$ and $|\br P|/|\delta_Q|\le1/c$. In
the term $\mathfrak a\,L_0\hat G$ the coefficient at $Q$ is $\sum_{K+P=Q}\mathfrak a_K\br P\hat g_P$; here
$|\br P|\le|\br Q|+m_0|K|$ with $m_0=\max(p,q)$, and hence $|\br P|/|\delta_Q|\le(1+m_0|K|)/c$.

\emph{Step 5: majorants.} For a series $F=\sum f_QX^Q$ let \[\|F\|(\tau)=\sum_{n}\bigl(\sum_{|Q|=n}|f_Q|\bigr)\tau^n.\]
Then \[\|FG\|\ll\|F\|\,\|G\|,\text{ and }\|\Lambda\|,\|\hat G\|\ll\omega:=\|W\|.\] Put \[\alpha_1=\|\mathfrak a\|,\;\alpha_2=\|\mathfrak b\|\text{ and
}\tilde\alpha_1=\Bigl\|\sum_K|\mathfrak a_K|(1+m_0|K|)X^K\Bigr\|.\]These are convergent for $|\tau|<R$, have non-negative
coefficients, and vanish at $\tau=0$. By Step~4 and $|\delta_Q|\ge c$, for every $n\ge1$,
\[
\omega_n\le\frac1c\Bigl[\Theta\bigl(\tau,\omega(\tau)\bigr)\Bigr]_n,\qquad
\Theta(\tau,\omega)=(2+|\mu|)\,\omega^2+\tilde\alpha_1\,\omega+\bigl(\omega^2+|\lambda_1|\alpha_1+\alpha_1\,\omega+\alpha_2\bigr)(1+\omega),
\]
where $[\cdot]_n$ is the coefficient of $\tau^n$. As in Step~3, the right-hand side involves only $\omega_k$ with
$k<n$. Now $\Theta$ is a polynomial in $\omega$ whose coefficients are convergent series with non-negative
coefficients, and $\Theta(0,0)=\partial_\omega\Theta(0,0)=0$. By the implicit function theorem,
$\psi=c^{-1}\Theta(\tau,\psi)$ has a unique holomorphic solution with $\psi(0)=0$. Its coefficients satisfy
$\psi_n=c^{-1}[\Theta(\tau,\psi)]_n$, so they are non-negative, and by induction $\omega_n\le\psi_n$ for all $n$.
Hence $\omega$ converges for $|\tau|<\rho_0$, the radius of $\psi$. Therefore \[\sum_Q|w_Q|\,\rho^{|Q|}<\infty\text{ for
}\rho<\rho_0,\]so $W$ converges on $\{|x|,|y|<\rho_0\}$. Then $\hat G$ (a sub-series of $W$) converges, and so does
\[\Lambda(h)=\sum_n\lambda_nh^n,\text{ for }|h|<\rho_0^{\,p+q}.\]
\end{proof}

\begin{corollary}\label{cor:rational}
Under the hypotheses of Theorem~\ref{thm:rational}:
\begin{enumerate}
\item Equation~\eqref{eq:delta} has two solutions $u_j=x^{\lambda_j(h)/p}G_j$, $j=1,2$, as in
Theorem~\ref{thm:rational}, with $\lambda_j(0)=\lambda_j$. Locally on $\{xy\neq0\}$ near $0$, every solution of
\eqref{eq:delta} is of the form $f(h)\,u_1+g(h)\,u_2$, with $f$ and $g$ holomorphic.
\item There is a Frobenius type solution on $\ell_j$ if and only if $\lambda_j(h)\equiv\lambda_j$ is constant.
\item Let $\rho>0$ be smaller than the radii of convergence of $G_1,G_2$ and $|h_0|$ small. On the local leaf
$C_{h_0}\cap\Delta_\rho$, the restrictions of $u_1,u_2$ form a basis of $\mathcal S_{h_0}$ of eigenvectors of $M_{h_0}$,
with eigenvalues $e^{2\pi i\lambda_1(h_0)}\neq e^{2\pi i\lambda_2(h_0)}$. In particular $u_j$ is not single-valued on the
local leaf unless $\lambda_j(h_0)\in\Z$, and the conclusion of Corollary~\ref{cor:monodromy} holds in this generality:
a Frobenius type solution on $\ell_j$ exists if and only if the family is isomonodromic for the eigenvalue branch
$e^{2\pi i\lambda_j(h_0)}$.
\item $v_j=\lambda_j(h)+L_0G_j/G_j$ is a holomorphic solution of the Riccati PDE near $0$, with $v_j(0)=\lambda_j$.
\end{enumerate}
\end{corollary}

\begin{proof}
(1) Since $-\mu\notin\Z$, Theorem~\ref{thm:rational} applies to both roots $\lambda_1$ and $\lambda_2$, which gives
$u_1$ and $u_2$. Let $\mathcal W=u_1L_0u_2-u_2L_0u_1$. Then
\[\mathcal W=x^{(\lambda_1+\lambda_2)(h)/p}\bigl[G_1G_2(\lambda_2(h)-\lambda_1(h))+G_1L_0G_2-G_2L_0G_1\bigr],\] and the bracket
equals $\mu\neq0$ at $0$. Near a point of $\{xy\neq0\}$, $h$ is a submersion and $L_0$ spans $\ker dh$, so
$(t,h)=(p^{-1}\log x,h)$ are local coordinates with $L_0=\partial_t$. In these coordinates \eqref{eq:delta} is a
second order linear ODE in $t$ depending on the parameter $h$, with fundamental system $u_1,u_2$. So
$u=fu_1+gu_2$ with $L_0f=L_0g=0$, and functions killed by $L_0$ are locally functions of $h$.

(2) If $\lambda_j(h)\equiv\lambda_j$, then $u_j=x^{\lambda_j/p}G_j$ is of Frobenius type with initial point
$(\lambda_j/p,0)\in\ell_j$. Conversely, let $x^ry^sF$ be a Frobenius solution with $(r,s)\in\ell_j$. Since
$pr-qs=\lambda_j$, we have $x^ry^s=x^{\lambda_j/p}h^{s/p}$. Let $F_0(h)$ be the weight-$0$ part of $F$; then
$F_0(0)=1$. Dividing by the function $h^{s/p}F_0(h)$, which is annihilated by $L_0$, yields a normalized solution
$x^{\lambda_j/p}\,(F/F_0)$ with constant exponent. By uniqueness in Theorem~\ref{thm:rational}, $\lambda_j(h)\equiv\lambda_j$.

(3) Take $|h_0|<\rho^{p+q}$, also small enough for $\lambda_j(h_0)$ to be defined, so that the local leaf is an annulus
by Lemma~\ref{lem:annulus}, and use that lemma
and the universal cover $\tilde{\mathcal A}$ described after Definition~\ref{def:monodromy}. Let
$\tilde\sigma(\tau)=\sigma(e^\tau)$. Since $x\circ\tilde\sigma=e^{p\tau}$, the function $p\tau$ is a global branch of
$\log x$ on $\tilde{\mathcal A}$. Lifting $u_j$ with this branch, we get
\[
u_j\circ\tilde\sigma(\tau)=e^{\lambda_j(h_0)\tau}\,G_j\bigl(\sigma(e^\tau)\bigr),\qquad \tau\in\tilde{\mathcal A}.
\]
Here $G_j\circ\sigma$ is a single-valued holomorphic function on $\mathcal A$, because $\sigma(\mathcal A)\subset\Delta_\rho$
lies in the domain of convergence of $G_j$. So the second factor is invariant under $\tau\mapsto\tau+2\pi i$, while the
first one is multiplied by $e^{2\pi i\lambda_j(h_0)}$. Hence continuation along $\gamma$ multiplies the germ of $u_j$ by
$e^{2\pi i\lambda_j(h_0)}$. Both restrictions solve $\mathcal L_{h_0}$ by Lemma~\ref{lem:annulus}(3). They are
independent because the Wronskian in (1) does not vanish near $0$. The eigenvalues are distinct, since
$\lambda_1(h_0)-\lambda_2(h_0)\to-\mu\notin\Z$. Finally, $e^{2\pi i\lambda_j(h_0)}$ is independent of $h_0$ if and only
if $\lambda_j$ is constant, because $\lambda_j$ is continuous on a connected punctured disc.

(4) $L_0u_j/u_j=\lambda_j(h)+L_0G_j/G_j$, which is holomorphic since $G_j(0)=1$; apply Proposition~\ref{prop:riccati}.
\end{proof}

In particular, the Riccati PDE has holomorphic solutions through both equilibria even when no Frobenius solution
exists, so the uniqueness of Corollary~\ref{cor:riccati-unique} fails at a saddle; the obstruction is the weight-$0$
part of $v_j-\lambda_j$, which is absorbed into $\lambda_j(h)$.

\subsection{Relation with Floquet theory}\label{sec:floquet}
On a fixed regular leaf, the leafwise equation $\mathcal L_{h_0}$ of Lemma~\ref{lem:annulus} is a linear ODE on an
annulus, whose coefficients $A(\sigma(s))$, $B(\sigma(s))$ are Laurent series in $s$. In the variable $\tau=\log s$ its
coefficients are $2\pi i$-periodic. When the monodromy has distinct eigenvalues, classical Floquet theory
\cite{Floquet1883,MagnusWinkler1966}, \cite[Chapter~3]{CL} gives solutions of the form $s^{\lambda}\,P(s)$, with $P$ a
Laurent series on the annulus and $e^{2\pi i\lambda}$ an eigenvalue of the monodromy $M_{h_0}$. By
Corollary~\ref{cor:rational}(3), the restriction of $u_j$ to the leaf is exactly of this form, with
$\lambda=\lambda_j(h_0)$ and $P=G_j\circ\sigma$. So, leaf by leaf, Theorem~\ref{thm:rational} recovers the Floquet
decomposition.

The theorem adds two things to Floquet theory. First, Floquet theory determines $\lambda(h_0)$ only modulo $\Z$, and
only for $h_0\neq0$. Theorem~\ref{thm:rational} singles out a choice $\lambda_j(h)$ that remains holomorphic
\emph{across} $h=0$, precisely where the annuli degenerate. Second, the Floquet factors $G_j\circ\sigma$, which a priori
live on separate annuli, turn out to be restrictions of a \emph{single} germ $G_j\in\C\{x,y\}$ at the singular point.
It is this holomorphic dependence on the leaf, through the singular leaf, that we exploit for Hill's equation in
Part~\ref{part:applications}.

\subsection{The Riccati PDE at a rational saddle}
In the Poincar\'e domain, holomorphic solutions of the Riccati PDE through a non-resonant equilibrium exist and
are unique (Corollary~\ref{cor:riccati-unique}). The same holds at a rational saddle, although \emph{every} point of
the indicial conic is resonant there.

\begin{theorem}[Briot--Bouquet for the Riccati PDE]\label{thm:riccati-saddle}
In the setting of Theorem~\ref{thm:rational} ($\mu\notin\Z$), every germ $v\in\cO_0$ solving
$L_0v=-v^2-Av-B$ satisfies $v(0)\in\{\lambda_1,\lambda_2\}$. For each $j$ there is exactly one such germ with
$v(0)=\lambda_j$, namely
\[
v_j=\lambda_j(h)+\frac{L_0G_j}{G_j},
\]
where $u_j=x^{\lambda_j(h)/p}G_j$ is the solution of Theorem~\ref{thm:rational}.
\end{theorem}

\begin{proof}
$L_0$ vanishes at $0$, so $\chi(v(0))=0$. Existence is Corollary~\ref{cor:rational}(4). For uniqueness, let
$v(0)=\lambda_1$ and decompose $v=v_0(h)+v'$, where $v'$ has no monomials of weight $0$. Put
\[\varphi=\sum_{\br Q\neq0}(v'_Q/\br Q)\,X^Q.\]This series converges because $|\br Q|\ge1$, and it satisfies
$L_0\varphi=v'$. Then $u=x^{v_0(h)/p}e^{\varphi}$ satisfies $L_0u=v\,u$, so $u$ solves \eqref{eq:delta} by
Proposition~\ref{prop:riccati}. Let $G=e^{\varphi}$, and let $G_0(h)$ be its weight-$0$ part; then $G_0(0)=1$.
Division by a function of $h$ preserves the weight decomposition, so $\tilde G=G/G_0$ has weight-$0$ part $1$.
Moreover \[\tilde u=u/G_0(h)=x^{v_0(h)/p}\tilde G\] is still a solution, since $L_0(G_0(h))=0$. By the uniqueness in
Theorem~\ref{thm:rational}, $v_0=\lambda_1(h)$ and $\tilde G=G_1$. Hence \[v=L_0u/u=L_0\tilde u/\tilde u=v_1.\]
\end{proof}

\subsection{The first exponent coefficient}
The exponent $\lambda(h)$ is computable. Its first coefficient already shows that $\lambda(h)$ is, in general, not
constant.

\begin{proposition}\label{prop:first-coefficient}
Let $p=q=1$ (so $h=xy$), let $\mathfrak a_{ij},\mathfrak b_{ij}$ be the Taylor coefficients of $A-A(0)$ and $B-B(0)$, and put
$\gamma_{ij}=\lambda_1\mathfrak a_{ij}+\mathfrak b_{ij}$. Under the hypotheses of Theorem~\ref{thm:rational},
\[
\lambda'(0)=\frac1\mu\left[\gamma_{11}-\frac{\gamma_{10}\,(\gamma_{01}+\mathfrak a_{01})}{1-\mu}
-\frac{\gamma_{01}\,(\gamma_{10}-\mathfrak a_{10})}{1+\mu}\right].
\]
\end{proposition}

\begin{proof}
In the recursion of Step~3 of the proof of Theorem~\ref{thm:rational}, the degree one terms are
$\hat g_{10}=-\gamma_{10}/\delta_{(1,0)}$ and $\hat g_{01}=-\gamma_{01}/\delta_{(0,1)}$, with $\delta_{(1,0)}=1-\mu$ and
$\delta_{(0,1)}=1+\mu$. At $Q=(1,1)$ we have $\delta_Q=-\mu$ and $w_{11}=\lambda'(0)$. The only contributions to
$\mathcal R(W)_{(1,1)}$ come from $\lambda_1\mathfrak a+\mathfrak b$, from $\mathfrak a\,L_0\hat G$ and from $(\lambda_1\mathfrak a+\mathfrak b)\hat G$:
\[
\mathcal R(W)_{(1,1)}=-\bigl[\gamma_{11}+\mathfrak a_{01}\hat g_{10}-\mathfrak a_{10}\hat g_{01}+\gamma_{01}\hat g_{10}+\gamma_{10}\hat g_{01}\bigr].
\]
Substituting and dividing by $-\mu$ gives the formula.
\end{proof}

For instance, for $A=y$ and $B=-(\tfrac19+x+xy)$, so that $\lambda_1=\tfrac13$ and $\mu=-\tfrac23$, the formula gives
$\lambda'(0)=-\tfrac65$, and continuing the recursion of Step~3 gives
\[\lambda(h)=\tfrac13-\tfrac65h+\tfrac{1593}{500}h^2-\cdots.\]

\subsection{The classical Frobenius locus}
Fix one indicial line, say $\ell_1$, i.e.\ the root $\lambda_1$. By Corollary~\ref{cor:rational}(2), classical
Frobenius solutions on $\ell_1$ exist if and only if $\lambda_1(h)\equiv\lambda_1$. Let $\pi_0$ denote the projection
of a series onto its weight-$0$ part. We decompose every germ $B$ uniquely as
\[
B=B'+\beta,\qquad \beta:=\pi_0\bigl(B-B(0)\bigr)\in h\,\C\{h\},\qquad B':=B-\beta .
\]
Thus $B'$ has no monomials $h^n$ with $n\ge1$, and $B'(0)=B(0)$. We show that the classical Frobenius condition on
$\ell_1$ reads $\beta=\Phi(A,B')$ for a locally holomorphic map $\Phi$.

\subsubsection*{Function spaces}
For $r>0$, let $H_r$ denote the space of power series \[F=\sum f_QX^Q\] such that
\[\|F\|_r:=\sum_Q|f_Q|\,r^{|Q|}<\infty.\]We use four elementary facts.
\begin{enumerate}
\item[(H1)] $H_r$ is a Banach algebra: $\|FG\|_r\le\|F\|_r\|G\|_r$. Its elements are holomorphic on $\Delta_r$, and
every germ belongs to some $H_r$.
\item[(H2)] $\|\pi_0F\|_r\le\|F\|_r$ and $\|(1-\pi_0)F\|_r\le\|F\|_r$.
\item[(H3)] If $\beta$ has weight $0$ and $G$ has no weight-$0$ part, then $\beta G$ has no weight-$0$ part, i.e.
$\pi_0(\beta G)=0$. Weights add under multiplication.
\item[(H4)] Let $m_0=\max(p,q)$, $F\in H_R$ with $F(0)=0$, and $r\le R/2$. Put \[\tilde F=\sum_K|f_K|(1+m_0|K|)X^K.\] Then
\[\|F\|_r\le\frac rR\|F\|_R\text{ and }\|\tilde F\|_r\le(1+m_0)\frac rR\|F\|_R.\]Indeed, for $k\ge1$ we have
\[(r/R)^k\le(r/R)\text{ and }(1+m_0k)(r/R)^{k}\le(r/R)(1+m_0k)2^{1-k}\le(1+m_0)(r/R).\]
\end{enumerate}
Let $H'_R=\{B'\in H_R:\ \pi_0(B'-B'(0))=0\}$, a closed subspace.

The next theorem says that, once $A$ and the part $B'$ of $B$ are fixed, there is exactly one choice of the
remaining part $\beta$, a power series in the first integral $h$, for which a classical Frobenius solution exists on
$\ell_1$. This choice depends holomorphically on $(A,B')$, so equations with classical Frobenius solutions form the
graph of a holomorphic map, and in this sense they are exceptional.

\begin{theorem}[The classical Frobenius locus]\label{thm:isomonodromic}
Let $L_0=p\,x\partial_x-q\,y\partial_y$ with $p,q$ coprime, and $h=x^qy^p$. Let $A,B'$ be germs at $0$, with $B'$
having no monomials $h^n$ for $n\ge1$. Let $\lambda_1$ be a root of $\nu^2+A(0)\nu+B'(0)$, and assume
$\mu=-A(0)-2\lambda_1\notin\Z$.
\begin{enumerate}
\item There is a unique $\beta\in h\,\C\{h\}$ such that, for $B=B'+\beta$, equation \eqref{eq:delta} has a
Frobenius type solution on $\ell_1$ (equivalently, $\lambda_1(h)\equiv\lambda_1$). We write $\beta=\Phi(A,B')$.
\item \emph{(Local holomorphy.)} Let $R>0$ and $(A^0,B'^0)\in H_R\times H'_R$ satisfy the hypotheses, and choose
$\lambda_1$ holomorphically in $(A(0),B'(0))$ near $(A^0(0),B'^0(0))$. Then there are a neighbourhood $U$ of
$(A^0,B'^0)$ in $H_R\times H'_R$ and a radius $r\in(0,R/2]$ such that $\Phi\colon U\to\pi_0H_r$ is holomorphic.
\end{enumerate}
Consequently, the equations admitting a classical Frobenius solution on $\ell_1$ are exactly those satisfying
$B=B'+\Phi(A,B')$. Here the equality $B=B'+\Phi(A,B')$ is understood as an equality of germs: $\Phi(A,B')$
converges for $|h|<r^{p+q}$, but need not belong to $\pi_0H_R$.
\end{theorem}

\begin{proof}
Write $\mathfrak a=A-A(0)$ and $\mathfrak b'=B'-B'(0)$. Then $\mathfrak a(0)=\mathfrak b'(0)=0$ and $\pi_0\mathfrak b'=0$.

\emph{Step 1: the functional equation.} We look for $u=x^{\lambda_1/p}(1+\hat G)$ with $\hat G(0)=0$ and
$\pi_0\hat G=0$. Repeating the computation of Step~2 in the proof of Theorem~\ref{thm:rational} with $\Lambda=0$ and
$B-B(0)=\mathfrak b'+\beta$, we get $\cP[u]=x^{\lambda_1/p}\,\mathcal N$, where
\[
\mathcal N=(L_0^2-\mu L_0)\hat G+\mathfrak a\,L_0\hat G+(\lambda_1\mathfrak a+\mathfrak b'+\beta)(1+\hat G).
\]
The operator $L_0^2-\mu L_0$ preserves weights, so $\pi_0(L_0^2-\mu L_0)\hat G=0$; and $\pi_0(\beta\hat G)=0$ by (H3).
Hence $\mathcal N=0$ if and only if
\[
\beta=-\pi_0\bigl[\mathfrak a\,L_0\hat G+(\lambda_1\mathfrak a+\mathfrak b')(1+\hat G)\bigr],\qquad
(L_0^2-\mu L_0)\hat G=-(1-\pi_0)\bigl[\mathfrak a\,L_0\hat G+(\lambda_1\mathfrak a+\mathfrak b')(1+\hat G)+\beta\hat G\bigr].
\]
Put $W=\beta+\hat G$, so that $\beta=\pi_0W$ and $\hat G=(1-\pi_0)W$. Let $\mathcal D'$ be the diagonal operator with
$\delta'_Q=1$ if $\br Q=0$ and $\delta'_Q=\br Q(\br Q-\mu)$ if $\br Q\neq0$. The system is equivalent to
\begin{equation}\label{eq:locus}
\mathcal D'W=-\Psi(W),\qquad \Psi(W):=\mathfrak a\,L_0\hat G+(\lambda_1\mathfrak a+\mathfrak b')(1+\hat G)+\beta\hat G,
\end{equation}
with $W(0)=0$. Conversely, every solution of \eqref{eq:locus} gives a Frobenius type solution
$x^{\lambda_1/p}(1+\hat G)$ on $\ell_1$, with initial point $(\lambda_1/p,0)$, for $B=B'+\pi_0W$.

\emph{Step 2: formal uniqueness.} Since $\mu\notin\Z$, $\delta'_Q\neq0$ for all $Q\neq0$. Each term of $\Psi(W)$ is
either independent of $W$ (namely $\lambda_1\mathfrak a+\mathfrak b'$), or a product of a $W$-dependent factor with a factor vanishing
at $0$ ($\mathfrak a$, $\lambda_1\mathfrak a+\mathfrak b'$, or a component of $W$). Hence $\Psi(W)_Q$ depends only on the $w_P$ with $|P|<|Q|$, and
\eqref{eq:locus} has exactly one formal solution with $W(0)=0$. Now suppose $B=B'+\tilde\beta$ admits a Frobenius
solution on $\ell_1$. Then $\lambda_1(h)\equiv\lambda_1$ by Corollary~\ref{cor:rational}(2), and the normalized
solution of Theorem~\ref{thm:rational} has the form $x^{\lambda_1/p}(1+\hat G_1)$ with $\pi_0\hat G_1=0$. So
$\tilde W=\tilde\beta+\hat G_1$ is a formal solution of \eqref{eq:locus}, hence $\tilde W=W$ and $\tilde\beta=\pi_0W$.
This proves the uniqueness in (1). It also shows that any convergent solution must coincide with the one
constructed below.

\emph{Step 3: estimates for the new operator.} Let $c=\min\{1,\tfrac12\operatorname{dist}(\mu,\Z)\}$. If $\br Q=0$ then
$|\delta'_Q|=1\ge c$. If $k=\br Q\neq0$ then \[|k(k-\mu)|\ge|k|\operatorname{dist}(\mu,\Z)\ge c(1+|k|).\]Thus
\[|\delta'_Q|\ge c(1+|\br Q|)\text{ for all }Q\neq0.\]For $F,G\in H_r$, and for $\mathfrak a$ such that $\tilde{\mathfrak a}\in H_r$, this gives
\begin{equation}\label{eq:E1E2}
\text{(E1)}\ \ \|\mathcal D'^{-1}F\|_r\le c^{-1}\|F\|_r,\qquad
\text{(E2)}\ \ \|\mathcal D'^{-1}(\mathfrak a\,L_0G)\|_r\le c^{-1}\|\tilde{\mathfrak a}\|_r\,\|G\|_r .
\end{equation}
For (E2): the coefficient of $X^Q$ in $\mathfrak a\,L_0G$ is \[\sum_{K+P=Q}\mathfrak a_K\br Pg_P,\text{ and }|\br P|\le|\br Q|+m_0|K|.\]So
$|\br P|/|\delta'_Q|\le(1+m_0|K|)/c$, and summing gives the claim. The gain of $\mathcal D'^{-1}$ thus compensates
the derivative $L_0$. The estimates \eqref{eq:E1E2} will be referred to as (E1) and (E2), respectively.
Compared with \eqref{eq:fixed-point}, the weight-$0$ divisor is now $1$ instead of $-\mu$,
which is why $c\le1$.

\emph{Step 4: contraction.} Let \[X_r=\{W\in H_r:\ W(0)=0\}\text{ and }T(W)=-\mathcal D'^{-1}\Psi(W).\]Put
$\alpha_r=\|\tilde{\mathfrak a}\|_r$ and $\varepsilon_r=\|\lambda_1\mathfrak a+\mathfrak b'\|_r$. Using (E1) and (E2) in \eqref{eq:E1E2}, together with (H1) and (H2), term by term,
\[
\|T(W)\|_r\le c^{-1}\bigl(\alpha_r\|W\|_r+\varepsilon_r(1+\|W\|_r)+\|W\|_r^2\bigr).
\]
Also, since \[\Psi(W_1)-\Psi(W_2)=\mathfrak a\,L_0(\hat G_1-\hat G_2)+(\lambda_1\mathfrak a+\mathfrak b'+\beta_1)(\hat G_1-\hat G_2)+(\beta_1-\beta_2)\hat G_2,\]
\[
\|T(W_1)-T(W_2)\|_r\le c^{-1}\bigl(\alpha_r+\varepsilon_r+\|W_1\|_r+\|W_2\|_r\bigr)\|W_1-W_2\|_r .
\]
By (H4), \[\alpha_r\le(1+m_0)\tfrac rR\|\mathfrak a\|_R\text{ and }\varepsilon_r\le\tfrac rR(|\lambda_1|\|\mathfrak a\|_R+\|\mathfrak b'\|_R).\]Choose
$r\le R/2$ with $\alpha_r\le c/8$ and $\varepsilon_r\le c^2/64$, and let $K_\eta=\{W\in X_r:\|W\|_r\le\eta\}$ with
$\eta=c/4$. For $W,W_1,W_2\in K_\eta$, using $c\le1$,
\begin{gather*}
\|T(W)\|_r\le\tfrac{c}{32}+\tfrac{c}{64}\bigl(1+\tfrac c4\bigr)+\tfrac{c}{16}\le\tfrac{29}{256}\,c<\eta,\\
\|T(W_1)-T(W_2)\|_r\le\bigl(\tfrac18+\tfrac{c}{64}+\tfrac12\bigr)\|W_1-W_2\|_r\le\tfrac34\|W_1-W_2\|_r .
\end{gather*}
The map $T$ preserves $X_r$: since $\mathfrak a(0)=\mathfrak b'(0)=0$ and $W(0)=0$, every term in $\Psi(W)$ has zero constant coefficient; the diagonal inverse preserves this property. Thus $T(K_\eta)\subseteq X_r$, and the preceding estimate gives $T(K_\eta)\subseteq K_\eta$. So $T$ is a contraction of the complete metric space $K_\eta$, and has a unique fixed point $W$ there. Its
coefficients satisfy the recursion of Step~2, so $W$ is the formal solution. Hence $\hat G=(1-\pi_0)W$ converges on
$\Delta_r$, and \[\beta=\pi_0W=\sum_n\beta_nh^n\]converges for $|h|<r^{p+q}$, because
\[\|\beta\|_r=\sum_n|\beta_n|r^{n(p+q)}.\] This proves (1).

\emph{Step 5: holomorphic dependence.} Let $d_0=\operatorname{dist}(\mu^0,\Z)>0$ at the base point. Let $V\subset\C^2$
be a neighbourhood of $(A^0(0),B'^0(0))$ on which $\lambda_1$ is holomorphic and bounded and
$|\mu-\mu^0|<d_0/4$; then \[\operatorname{dist}(\mu,\Z)>d_0/2,\text{ so }c\ge c_0:=\min\{1,d_0/4\}\text{ on }V.\]Let $U$ be the set of
$(A,B')\in H_R\times H'_R$ with $(A(0),B'(0))\in V$, $\|a-a^0\|_R<1$ and $\|b'-b'^0\|_R<1$. On $U$ the quantities
$|\lambda_1|$, $\|\mathfrak a\|_R$, $\|\mathfrak b'\|_R$ are bounded. By (H4), one radius $r$ satisfies the smallness conditions of
Step~4, with $c$ replaced by $c_0$, for all parameters in $U$. Hence $T_\theta$ ($\theta\in U$) is a contraction of
$K_{c_0/4}$ with the uniform constant $3/4$.

The map $(\theta,W)\mapsto T_\theta(W)$ is holomorphic on $U\times X_r$. Here one must not regard
$L_0\colon H_r\to H_r$ as a bounded operator: in general it is unbounded. Instead, for each admissible $\mu$ define
\[\mathcal K_\mu(\mathfrak a,G):=\mathcal D'^{-1}_\mu(\mathfrak a\,L_0G)\]for $(\mathfrak a,G)\in H_R\times H_r$ by the coefficient formula: if
$\mathfrak a=\sum_K \mathfrak a_KX^K$ and $G=\sum_Pg_PX^P$, its coefficient at $X^Q$ is
\[
[\mathcal K_\mu(\mathfrak a,G)]_Q=\sum_{K+P=Q}\frac{\br P}{\delta'_Q(\mu)}\,\mathfrak a_Kg_P .
\]
For arbitrary $\mathfrak a\in H_R$ and $0<r\le R/2$, the weighted series
\[\tilde{\mathfrak a}=\sum_K|\mathfrak a_K|(1+m_0|K|)X^K\]belongs to $H_r$; indeed,
\[\|\tilde{\mathfrak a}\|_r\le C_{r/R,m_0}\|\mathfrak a\|_R,\]where
\[C_{t,m_0}:=\sup_{n\ge0}(1+m_0n)t^n<\infty\text{ for }0<t<1.\]
Consequently, (E2) in \eqref{eq:E1E2} yields
\[\|\mathcal K_\mu(\mathfrak a,G)\|_r\le c^{-1}C_{r/R,m_0}\|\mathfrak a\|_R\|G\|_r.\]
This proves absolute convergence of the coefficient series and boundedness of the bilinear map, uniformly for $\mu$ in the parameter neighbourhood used below.
When $\mathfrak a(0)=0$, as in the fixed-point argument, the stronger estimate (H4) gives
\[\|\tilde{\mathfrak a}\|_r\le(1+m_0)(r/R)\|\mathfrak a\|_R;\]it is this additional smallness that is used in Step~4. Thus the derivative always occurs inside the compensated operator
$\mathcal K_\mu$, not as an independent bounded factor. The other terms in $T_\theta$ are compositions
of bounded projections, multiplication in $H_r$, and the diagonal inverse $\mathcal D'^{-1}_\mu$. The diagonal entries $1/\delta'_Q(\mu)$ and $\br P/\delta'_Q(\mu)$ are holomorphic in
$\mu$ on $\{|\mu-\mu^0|<d_0/2\}$, and they obey the bounds of Step~3 uniformly. For every $\mu_1$ with $|\mu_1-\mu^0|<d_0/4$ (that is, for every value of $\mu$ on $V$), the closed disk of radius
$d_0/4$ centred at $\mu_1$ lies in $\{|\mu-\mu^0|<d_0/2\}$, where these bounds hold. The Cauchy estimates
therefore bound the Taylor coefficients at $\mu_1$ by $C\,(d_0/4)^{-j}$, uniformly in $Q,P$ and in $\mu_1$. Consequently, $\mu\mapsto\mathcal D'^{-1}_\mu$ is holomorphic in operator norm and
$\mu\mapsto\mathcal K_\mu$ is holomorphic in the norm of bounded bilinear maps
$H_R\times H_r\to H_r$: apply the same Cauchy bounds to the coefficient multipliers of
$\mathcal K_\mu$, using the summable majorant $\tilde{\mathfrak a}$ from (E2). Both maps therefore admit locally
norm-convergent power series. Compositions of holomorphic maps between Banach spaces are holomorphic \cite{Mujica}.

The iterates $W_0=0$, $W_{n+1}(\theta)=T_\theta(W_n(\theta))$ are therefore holomorphic maps $U\to K_{c_0/4}$. They
converge uniformly on $U$, since \[\|W_{n+1}-W_n\|_r\le(3/4)^n\|W_1-W_0\|_r\text{ and }\|W_1\|_r\le c_0/4.\] A uniform limit
of holomorphic maps with values in a Banach space is holomorphic \cite{Mujica}. Hence $\theta\mapsto W(\theta)$ and
$\Phi=\pi_0\circ W\colon U\to\pi_0H_r$ are holomorphic. This proves (2). The final statement follows from (1) and the
uniqueness of the decomposition $B=B'+\beta$.
\end{proof}

\begin{corollary}[Codimension of the classical Frobenius locus]\label{cor:generic}
Let $U$, $r$ and $\Phi$ be as in Theorem~\ref{thm:isomonodromic}(2), and write \[\Phi(A,B')=\sum_{n\ge1}\Phi_n(A,B')\,h^n.\]
For $B\in H_R$ let $b_n$ be its coefficient of $h^n$, so that \[B'=B-\sum_{n\ge1}b_nh^n,\text{ and let
}\tilde U=\{(A,B)\in H_R\times H_R\colon\ (A,B')\in U\}.\]
\begin{enumerate}
\item For fixed $(A,B')$, among all $B=B'+\gamma(h)$ with $\gamma\in h\,\C\{h\}$, exactly one admits a classical
Frobenius solution on $\ell_1$, namely $\gamma=\Phi(A,B')$.
\item Each $\Phi_n\colon U\to\C$ is holomorphic. For every $N\ge1$, the set
\[
Z_N=\{(A,B)\in\tilde U:\ b_n=\Phi_n(A,B')\ \text{for }1\le n\le N\}
\]
is a closed complex submanifold of $\tilde U$ of codimension exactly $N$.
\item The set of $(A,B)\in\tilde U$ for which \eqref{eq:delta} has a classical Frobenius solution on $\ell_1$ is
\[\bigcap_{N\ge1}Z_N.\]In particular, it is contained in complex submanifolds of $\tilde U$ of arbitrarily large
finite codimension.
\end{enumerate}
\end{corollary}

\begin{proof}
(1) is Theorem~\ref{thm:isomonodromic}(1). For (2), the coefficient functional $\beta\mapsto\beta_n$ is continuous and
linear on $\pi_0H_r$, since \[|\beta_n|\le r^{-n(p+q)}\|\beta\|_r.\]Hence $\Phi_n$ is holomorphic. The maps $B\mapsto b_n$
(since $|b_n|\le R^{-n(p+q)}\|B\|_R$) and $B\mapsto B'$ (by (H2)) are continuous and linear on $H_R$. For each $N$, let $E_N$ be the closed Banach subspace of $H_R\times H_R$ consisting of pairs $(A,C)$ for which the coefficients of $h,h^2,\ldots,h^N$ in $C$ vanish. The bounded coefficient projections give a topological splitting
$H_R\times H_R\cong E_N\times\C^N$, $(A,B)\mapsto\bigl((A,C),(b_1,\dots,b_N)\bigr)$ with
\[C=B-\sum_{n\le N}b_nh^n=B'+\sum_{n>N}b_nh^n,\text{ and }\tilde U\]is open. Since
\[B'=C-\sum_{n>N}c_nh^n,\]where $c_n$ is the coefficient of $h^n$ in $C$, the map $(A,C)\mapsto(A,B')$ is continuous
and linear, so $(A,C)\mapsto\Phi_n(A,B')$ is holomorphic. In these coordinates $Z_N$ is the graph of the holomorphic
map $(A,C)\mapsto(\Phi_1(A,B'),\dots,\Phi_N(A,B'))$, defined on the open subset \[V_N=\{(A,C)\in E_N:(A,B')\in U\}\]of $E_N$ with values in $\C^N$. Indeed, $V_N$ is open because $(A,C)\mapsto(A,B')$ is continuous and $U$ is open. The graph of a holomorphic map
from an open subset of a Banach space to $\C^N$ is a relatively closed complex submanifold of codimension $N$ in $V_N\times\C^N\cong\tilde U$. Indeed, adding any finite linear combination of $h,\ldots,h^N$ preserves membership in $H_R$, and the projection $B\mapsto B'$ is unchanged by these additions. (3) follows from
Theorem~\ref{thm:isomonodromic}(1): a classical Frobenius solution on $\ell_1$ exists if and only if
\[\pi_0(B-B(0))=\Phi(A,B'),\]i.e.\ if and only if $b_n=\Phi_n(A,B')$ for all $n\ge1$.
\end{proof}

In the sense of (2), the conditions $b_n=\Phi_n(A,B')$ are independent: the first $N$ of them cut out a submanifold
of codimension exactly $N$. For $p=q=1$ the first condition is $\lambda'(0)=0$
(Proposition~\ref{prop:first-coefficient}). Indeed, $\lambda'(0)$ is affine in $b_1=\mathfrak b_{11}$ with leading coefficient
$1/\mu\neq0$, so it vanishes for exactly one value of $b_1$. Since a classical Frobenius solution forces
$\lambda'(0)=0$, that value is $\Phi_1(A,B')$.

There is a loss of radius in Theorem~\ref{thm:isomonodromic}(2): $\Phi$ takes values in $\pi_0H_r$ for some $r\le R/2$, and not necessarily
in $\pi_0H_R$. For this reason we describe the classical Frobenius locus through the finite truncations $Z_N$
of Corollary~\ref{cor:generic}, rather than as the graph $\{\pi_0(B-B(0))=\Phi(A,B')\}$ in $H_R\times H_R$. If
$\Phi$ mapped a neighbourhood in $H_R\times H'_R$ holomorphically into $\pi_0H_R$, this graph would be a closed
complex Banach submanifold of $H_R\times H_R$ (of infinite codimension), parametrized by $(A,B')$. We do not know
whether the loss of radius is genuine or only an artifact of the contraction argument.

Theorem~\ref{thm:isomonodromic} produces many equations with classical Frobenius solutions on $\ell_1$
that are not of the trivial form $A,B\in\C\{h\}$. For example, for $A=y$ and $B'=-(\tfrac19+x-y^2-3x^2y)$, the recursion of Step~2 gives
\[\beta(h)=-\tfrac95h-\tfrac{8883}{2000}h^2-\cdots.\]By Theorem~\ref{thm:isomonodromic}, for $B=B'+\beta$ the exponent
$\lambda_1(h)$ is constant, so the equation has a classical Frobenius solution on $\ell_1$.

\subsection{The integral case}\label{sec:integral}
When $\mu\in\Z$, the divisor $\delta_Q$ vanishes for $\br Q=\mu$ (if $\mu\neq0$), or for $\br Q=0$ (if $\mu=0$). In the
classical Frobenius theorem the root with larger real part always works, because the weights $q\ge0$ all have one
sign. At a saddle the weights $\br Q$ take both signs, so neither root is protected in advance. Protection is
restored when the coefficients are \emph{one-sided}, that is, when all their monomials have weights of one sign.
For $m\in\Z$, let $\mathcal S_m\subset\C[[x,y]]$ be the set of series whose monomials $X^Q$ satisfy
$m\,\br Q\le0$. It is a subring, stable under $L_0$, and it contains $\C[[h]]$.

\begin{theorem}[Integral case, one-sided coefficients]\label{thm:integral}
In the setting of Theorem~\ref{thm:rational}, assume instead that $\mu=m\in\Z\setminus\{0\}$, and that
$A-A(0)$ and $B-B(0)$ belong to $\mathcal S_m$, i.e.\ all their monomials have weight of the sign of $-m$ (or
weight $0$). Then the conclusions of Theorem~\ref{thm:rational} hold for the root $\lambda_1$: there are
$\lambda\in\C\{h\}$, $\lambda(0)=\lambda_1$, and $G\in\C\{x,y\}\cap\mathcal S_m$, normalized as before, such that
$x^{\lambda(h)/p}G$ solves \eqref{eq:delta}. The pair $(\lambda,G)$ is unique among pairs with $G\in\mathcal S_m$.
\end{theorem}

\begin{proof}[Proof of Theorem~\ref{thm:integral}]
Repeat the proof of Theorem~\ref{thm:rational} inside $\mathcal S_m$. If $W\in\mathcal S_m$, then
$\mathcal R(W)\in\mathcal S_m$, so the recursion only involves monomials with $k=\br Q$ satisfying $mk\le0$. For
those, $\delta_Q\neq0$: if $k=0$ then $\delta_Q=-m$; if $k\neq0$ then $k$ and $k-m$ are non-zero integers of the
same sign, so \[|\delta_Q|=|k|\,(|k|+|m|)\ge1+|k|.\] Thus \eqref{eq:divisor-bound} holds with $c=1$ for every $Q$ in
the support. The functional equation \eqref{eq:fixed-point} is the same as in Theorem~\ref{thm:rational}; only its
unknown is restricted to $\mathcal S_m$, which is a subring stable under $L_0$ and contains $a,b$. Therefore
\[\mathcal R(\mathcal S_m)\subset\mathcal S_m.\]Step~3 (lower-degree dependence) and Step~5 (the majorant $\Theta$ and the
comparison with $\psi=c^{-1}\Theta(\tau,\psi)$) use only \eqref{eq:divisor-bound}, the algebra property of
$\|\cdot\|$, and the vanishing of $a,b,\Lambda,\hat G$ at $0$. So they hold without change in $\mathcal S_m$. Uniqueness
within $\mathcal S_m$ follows from Step~3.
\end{proof}

Without the one-sided hypothesis, the recursion for $\lambda_1$ meets $\delta_Q=0$ at the monomials
of weight $\mu$. For each degree $n$ there is at most one such $Q$ with $|Q|=n$, since $\br Q=(p+q)q_1-qn$ determines
$q_1$, so the first of them is well defined. All coefficients of lower degree are determined uniquely, and
solvability at this $Q$ therefore requires a computable number $\kappa=\mathcal R(W)_Q$, the analogue of the classical
coefficient of the logarithm, to vanish. If $\kappa\neq0$, there is no formal solution $x^{\lambda(h)/p}G$ with
$\lambda(0)=\lambda_1$ (after normalizing the weight-$0$ part of $G$ to $1$, which is no loss of generality). As in the
classical theory one then expects logarithmic terms; Example~\ref{ex:bessel} shows them explicitly, although we do
not prove their appearance in general. These numbers $\kappa$ reappear in Section~\ref{sec:lame} as the leading
coefficients of the eigenvalue splittings of Lam\'e's equation.

The double root $\mu=0$ behaves differently. There $\delta_Q$ vanishes for every monomial $h^n$, so the equation for
$\lambda(h)$ can no longer be solved by the implicit function theorem, and the exponent may be ramified in $h$. For
example, the equation of Example~\ref{ex:saddle} with $c=0$ has the solutions $x^{\pm\sqrt{xy}}$, whose exponents
$\pm\sqrt h$ are Puiseux series. A complete treatment, via a Newton polygon in $(h,\lambda)$, is beyond the scope of
this paper.

\part{Hill's Lunar equation and other applications}\label{part:applications}

Part~\ref{part:applications} applies the theory to explicit equations. It can be read after the introduction and
Section~\ref{sec:riccati}, using only the following facts from Part~\ref{part:theory}. A square-type foliated equation
restricts, on each leaf of $\cF_L$, to a linear second order ODE, and $v=Lu/u$ turns it into a Riccati PDE
(Section~\ref{sec:riccati}). It has a regular singularity exactly when $\cF_L$ is analytically linearizable
(Theorem~\ref{thm:characterization}). At a rational saddle $L_0=p\,x\partial_x-q\,y\partial_y$, with first integral
$h=x^qy^p$ and $\mu\notin\Z$, it has solutions $x^{\lambda_j(h)/p}G_j(x,y)$ in which the exponent $\lambda_j(h)$ is a
convergent power series in $h$, computed together with $G_j$ by an explicit triangular recursion
(Theorem~\ref{thm:rational} and Corollary~\ref{cor:rational}); on each leaf, $\lambda_j(h_0)$ is a Floquet exponent
(Section~\ref{sec:floquet}). Solutions with a constant exponent exist exactly on the classical Frobenius locus of
Theorem~\ref{thm:isomonodromic}, and when $\mu\in\Z$ the recursion may meet a computable obstruction $\kappa$
(Section~\ref{sec:integral}).

Section~\ref{sec:examples} illustrates these results on equations whose solutions can be written in closed form.
Section~\ref{sec:hill} places Hill's equation on a leaf of a single operator and turns the recursion into an algorithm
for its characteristic exponent. The remaining sections apply this to three physical problems: the motion of the
Moon (Section~\ref{sec:moon}), ion trapping (Section~\ref{sec:mathieu}) and Lam\'e's equation
(Section~\ref{sec:lame}). Page numbers in our citations of Hill's memoir refer to its reprint in \emph{Acta
Mathematica} \cite{Hill1886}.

\section{Examples at a singular point of $L$}\label{sec:examples}

The following examples illustrate the main results in explicit cases, where the solutions can be written in
closed form, often in terms of Bessel-type functions. They cover a resonant node, first linear and then linearizable but
not linear, followed by a node that cannot be linearized and so has no regular singularity, a rational saddle
whose exponents vary from leaf to leaf, and an integral case in which logarithms appear.

The first three examples below share the eigenvalues $a=1$, $b=2$ (resonant node, item (2) of
Corollary~\ref{cor:dictionary}) and the Riccati data $A\equiv-1$, so $\chi(\lambda)=\lambda^2-\lambda$,
$\lambda_1=0$, $\lambda_2=1$, $\mu=1$. Here $\Gamma=\{q_1+2q_2\}=\{1,2,3,\dots\}$, so $\ell_1=\{r+2s=0\}$ is
resonant ($\mu=\br{(1,0)}$) and $\ell_2=\{r+2s=1\}$ is non-resonant.

A useful observation: if $\varphi$ is weighted homogeneous of degree $k$ for $L_0$ (that is, $L_0\varphi=k\varphi$)
and \[B-B_0=-\varphi,\;A\equiv A_0,\]then \[F=\sum_n c_n\varphi^n\] and the recursion \eqref{eq:recursion} becomes
$nk(nk-\mu)\,c_n=c_{n-1}$. Hence \[F=\,_0F_1\!\left(;1-\tfrac{\mu}{k};\tfrac{\varphi}{k^2}\right),\]a Bessel type
function of $\varphi$.

\begin{example}[Linear resonant node; Bessel type solutions]\label{ex:linear}
Let \[L_0=x\partial_x+2y\partial_y,\;A\equiv-1,\;B=-(x^2+y).\]By \eqref{eq:delta} the equation is
\begin{equation*}
x^2u_{xx}+4xy\,u_{xy}+4y^2u_{yy}+2y\,u_y-(x^2+y)\,u=0 .
\end{equation*}
Put \[z=x^2+y\;(\text{so }L_0z=2z)\text{ and }h=x^2/y\;(\text{so }L_0h=0).\]With $k=2$ the observation above gives
\[F_1=\,_0F_1(;\tfrac12;\tfrac z4)=\cosh\sqrt z\text{ and }F_2=\,_0F_1(;\tfrac32;\tfrac z4)=\sinh\sqrt z/\sqrt z,\]
both entire in $(x,y)$. Therefore:
\begin{enumerate}
\item for $(r,s)\in\ell_1$: $u=x^ry^s\cosh\sqrt{x^2+y}$; in particular $u=\cosh\sqrt{x^2+y}$ is entire;
\item for $(r,s)\in\ell_2$: $u=x^ry^s\,\dfrac{\sinh\sqrt{x^2+y}}{\sqrt{x^2+y}}$; in particular
$u=x\,\dfrac{\sinh\sqrt{x^2+y}}{\sqrt{x^2+y}}$ is entire;
\item general solution locally on simply connected coordinate neighborhoods where $h=x^2/y$ is a
valid transverse coordinate:
\[
u=f\!\left(\tfrac{x^2}{y}\right)\cosh\sqrt{x^2+y}+g\!\left(\tfrac{x^2}{y}\right)x\,\frac{\sinh\sqrt{x^2+y}}{\sqrt{x^2+y}} .
\]
\end{enumerate}
\emph{Resonance without obstruction.} On $\ell_1$, the recursion at $Q=(1,0)$ reads $0\cdot d_{(1,0)}=-\beta_{(1,0)}=0$,
so $d_{(1,0)}$ is free. The Frobenius solution on $\ell_1$ is therefore not unique: $F_1+c\,xF_2$ works for every $c$.
If instead $B=-(x^2+y)-\varepsilon x$ with $\varepsilon\neq0$, the same equation reads $0=\varepsilon$, so there is
no Frobenius solution on $\ell_1$ and logarithmic terms must appear. This mirrors exactly the classical Frobenius
picture for roots differing by an integer.

\emph{Riccati PDE.} By Corollary~\ref{cor:riccati-unique}, $Lv=-v^2+v+(x^2+y)$ has the holomorphic solutions
\begin{gather*}
v_1=\sqrt{x^2+y}\,\tanh\sqrt{x^2+y}\quad(v_1(0)=\lambda_1=0),\\
v_2=\sqrt{x^2+y}\,\coth\sqrt{x^2+y}\quad(v_2(0)=\lambda_2=1).
\end{gather*}
Both are even functions of $\sqrt{x^2+y}$, hence holomorphic. $v_2$ is the unique solution with $v(0)=1$
($\ell_2$ is non-resonant). Through $v(0)=0$ there is a one-parameter family $v=L_0u/u$ with
$u=F_1+c\,xF_2$, which reflects the resonance of $\ell_1$.
\end{example}

\begin{example}[Non-linear but linearizable $L$]\label{ex:linearizable}
Let \[\Phi(x,y)=(x,\,y+x^3)\text{ and }L=\Phi^*L_0=x\partial_x+(2y-x^3)\partial_y.\]Take \[A\equiv-1\text{ and
}B=B_{\rm ex\ref{ex:linear}}\circ\Phi=-(x^2+y+x^3).\]By \eqref{eq:main-expanded} the equation is
\[
x^2u_{xx}+2x(2y-x^3)\,u_{xy}+(2y-x^3)^2u_{yy}+(2y-4x^3)\,u_y-(x^2+y+x^3)\,u=0 .
\]
In these coordinates it is \emph{not} in the form \eqref{eq:regular-form}: the $u_{xy}$ coefficient is not a
constant multiple of $xy$. Nevertheless, by Theorem~\ref{thm:characterization} it has a regular singularity,
and its solutions are those of Example~\ref{ex:linear} composed with $\Phi$. With $Z=x^2+y+x^3$:
\[
u=x^{r}(y+x^3)^{s}\cosh\sqrt Z\ \ (r+2s=0),\qquad
u=x^{r}(y+x^3)^{s}\,\frac{\sinh\sqrt Z}{\sqrt Z}\ \ (r+2s=1).
\]
In particular $\cosh\sqrt{x^2+y+x^3}$ and $x\sinh\sqrt{x^2+y+x^3}\big/\sqrt{x^2+y+x^3}$ are entire solutions.
This example shows why the Frobenius ansatz must be made in linearizing coordinates.
\end{example}

\begin{example}[Non-linearizable resonant node]\label{ex:nonlinearizable}
Let \[L=x\partial_x+(2y+x^2)\partial_y,\]the Poincar\'e--Dulac normal form with $\varepsilon=1$, which is not
linearizable. Take $A\equiv-1$, $B\equiv0$. The equation is
\[
x^2u_{xx}+2x(2y+x^2)\,u_{xy}+(2y+x^2)^2u_{yy}+(3x^2+2y)\,u_y=0 .
\]
By Theorem~\ref{thm:characterization} it has no regular singularity at $0$, in any coordinates. Concretely,
$H=y/x^2-\log x$ is a (multivalued) first integral of $L$, and $L(\log x)=1$. In the coordinates $(t,H)=(\log x,H)$
we have $L=\partial_t$, and the general solution locally on simply connected coordinate neighborhoods in $\{x\neq0\}$
where $(t,H)$ are holomorphic coordinates is
\[
u=f\!\left(\frac{y}{x^2}-\log x\right)+g\!\left(\frac{y}{x^2}-\log x\right)x .
\]
Thus solutions such as \[u=y/x^2-\log x\text{ or }u=x\,(y/x^2-\log x)^2\]carry logarithms, although $\chi$ has simple roots
and there is no ``Frobenius'' obstruction coming from the indicial data. Compare with Example~\ref{ex:linear}:
the linear node has the meromorphic first integral $x^2/y$ and a whole line of monomial solutions $x^ry^s$, while
here the first integral is transcendental.
\end{example}

\begin{example}[Rational saddle: exponents varying along the leaves]\label{ex:saddle}
Let \[L_0=x\partial_x-y\partial_y\;(p=q=1,\text{ first integral }h=xy),\] $A\equiv0$ and $B=-(c^2+xy)$, with
$c\in\C^*$. By \eqref{eq:delta} the equation is
\begin{equation*}
x^2u_{xx}-2xy\,u_{xy}+y^2u_{yy}+x\,u_x+y\,u_y-(c^2+xy)\,u=0 .
\end{equation*}
Here $\chi(\lambda)=\lambda^2-c^2$, $\lambda_{1,2}=\pm c$, $\mu=-2c$, and the indicial conic is
$\{(r-s)^2=c^2\}$, i.e.\ $\ell_\pm=\{r-s=\pm c\}$. The leafwise polynomial is $\chi_h(\lambda)=\lambda^2-(c^2+h)$,
so the leafwise exponents are $\lambda_\pm(h)=\pm\sqrt{c^2+xy}$, where the square-root branch is chosen holomorphically near $h=0$ with value $c$ at $0$. All powers $x^{\lambda_\pm(h)}$ below are understood on a simply connected domain in $\{x\neq0\}$, with a fixed
branch of $\log x$. By Proposition~\ref{prop:saddle}:
\begin{enumerate}
\item the general solution, locally on $(\C^*)^2$ near $0$, is
\[
u=f(xy)\,x^{\sqrt{c^2+xy}}+g(xy)\,x^{-\sqrt{c^2+xy}} ;
\]
\item there is \emph{no} Frobenius type solution on either line, because $\lambda_\pm(h)$ is not constant.
Directly: the recursion \eqref{eq:recursion} at $Q=(1,1)$, where $\br Q=0$, reads $0\cdot d_{(1,1)}=-\beta_{(1,1)}=1$;
\item when expanded around the constant exponent $c$, the solutions contain unbounded powers of $\log x$:
\[
x^{\sqrt{c^2+xy}}=x^{c}\Bigl(1+\frac{xy}{2c}\log x+\frac{(xy)^2}{8c^2}\bigl((\log x)^2-\tfrac{1}{c}\log x\bigr)+O\bigl((xy)^3\bigr)\Bigr);
\]
\item here $A,B$ are constant on the leaves, so the leafwise equation lives on the whole curve
$\{xy=h_0\}\cong\C^*$ (Corollary~\ref{cor:monodromy}). Its monodromy has eigenvalues $e^{\pm2\pi i\sqrt{c^2+h_0}}$, which
vary with $h_0$, so the family is not isomonodromic. The solutions $x^{\pm\sqrt{c^2+xy}}$ continue analytically
around the leaf, but they are single-valued on it only when $\sqrt{c^2+h_0}\in\Z$;
\item the Riccati PDE $L_0v=-v^2+c^2+xy$ has the holomorphic solutions $v_\pm=\pm\sqrt{c^2+xy}$, with
$v_\pm(0)=\pm c$, although no Frobenius solution exists (Corollary~\ref{cor:rational}(4)).
\end{enumerate}
If instead $B=-c^2$ (Euler case), then $\lambda_\pm(h)\equiv\pm c$ and every $x^ry^s$ with $(r,s)\in\ell_\pm$ is a
solution. The single term $-xy$ in $B$, which is resonant, destroys the whole family of Frobenius solutions.
\end{example}

\begin{example}[Integral case: a Bessel logarithm]\label{ex:bessel}
Let $L_0=x\partial_x-y\partial_y$, $A\equiv0$ and $B=-(\tfrac14+x)$. The equation is
\[
x^2u_{xx}-2xy\,u_{xy}+y^2u_{yy}+x\,u_x+y\,u_y-\bigl(\tfrac14+x\bigr)u=0 ,
\]
with $\chi(\lambda)=\lambda^2-\tfrac14$, so the roots are $\pm\tfrac12$ and $\mu\in\{\pm1\}\subset\Z$.
\begin{enumerate}
\item For the root $\tfrac12$ we have $m=-1$, and $B-B(0)=-x$ has weight $1$, so it lies in $\mathcal S_{-1}$.
Theorem~\ref{thm:integral} applies; here $\lambda(h)\equiv\tfrac12$ and
\[
u_1=\sqrt x\sum_{n\ge0}\frac{x^n}{n!\,(n+1)!}=I_1(2\sqrt x).
\]
\item For the root $-\tfrac12$ we have $m=1$, and the source $-x$ has weight $1=m$. The recursion meets
$\delta_{(1,0)}=0$ with $\mathcal R_{(1,0)}=1$, so $\kappa=1\neq0$ (Section~\ref{sec:integral}). The second
solution is $u_2=K_1(2\sqrt x)$, which contains the term \[\tfrac12\sqrt x\,\log x\,(1+O(x)).\]
\end{enumerate}
Indeed, for functions of $x$ alone the equation reads \[(x\partial_x)^2u-(\tfrac14+x)u=0.\]With $z=2\sqrt x$ this is
the modified Bessel equation of order $1$. Locally, the general solution is
\[f(xy)\,I_1(2\sqrt x)+g(xy)\,K_1(2\sqrt x).\]This is the exact PDE analogue of the classical integral case of
Frobenius.
\end{example}

\section{Hill's equation as a leaf}\label{sec:hill}
\label{sec:hill-leaf}
\emph{Notation.} From now on $p=q=1$: we work with the saddle $M=X\partial_X-Y\partial_Y$ and its first integral $H=XY$,
in coordinates $(X,Y)$. In each of the following sections, some letters keep their classical meanings: in
Section~\ref{sec:moon}, $m$ is the ratio of mean motions and $q^2,q_k$ are Adams' coefficients; in
Section~\ref{sec:mathieu}, $a$ and $q$ are the parameters of Mathieu's equation and $\beta$ its characteristic exponent;
in Section~\ref{sec:lame}, $q$ is the nome. They should not be confused with the eigenvalues $a,b$, the integers $p,q$
and the series $\beta(h)$ of Part~\ref{part:theory}.

Let
\begin{equation}\label{eq:hill}
w''(\tau)+\Theta(\tau)\,w(\tau)=0,\qquad \Theta(\tau)=\Theta_0+\sum_{k\ge1}\theta_k\cos(2k\tau),
\end{equation}
with $\Theta_0,\theta_k\in\C$ and $|\theta_k|\le C\rho^{-k}$ for some $C>0$ and $\rho>1$, so that $\Theta$ is analytic in a
horizontal strip around the real axis. A solution of the form
$e^{i\mathbf c\tau}P(\tau)$, with $P$ $\pi$-periodic, is a \emph{Floquet solution} and $\mathbf c$ (defined modulo
$2\Z$ and up to sign) is the \emph{characteristic exponent}. In Hill's lunar theory, $\mathbf c$ is the ratio of the
synodic month to the anomalistic month \cite{Hill1886}.

We work in coordinates $(X,Y)$ with $M=X\partial_X-Y\partial_Y$ (so $p=q=1$) and first integral $H=XY$. Put
\[
B_\Theta(X,Y)=-\frac{\Theta_0}{4}-\sum_{k\ge1}\frac{\theta_k}{8}\,\bigl(X^k+Y^k\bigr),
\]
which is holomorphic on $\{|X|,|Y|<\rho\}$, and consider $\cP_{M,0,B_\Theta}[u]=M^2u+B_\Theta u=0$.

\begin{proposition}\label{prop:hill-leaf}
Let $0<|H_0|<\rho^2$, fix a square root $\sqrt{H_0}$, and parametrize the leaf $\{XY=H_0\}$ by
$s=\sqrt{H_0}\,e^{2i\tau}$, i.e.\ $(X,Y)=(s,H_0/s)$ as in Lemma~\ref{lem:annulus}. Then the leafwise equation
$\mathcal L_{H_0}$ becomes
\[
w''+\Theta_{H_0}(\tau)\,w=0,\qquad \Theta_{H_0}(\tau)=\Theta_0+\sum_{k\ge1}\theta_k\,H_0^{k/2}\cos(2k\tau),
\]
where $H_0^{k/2}=(\sqrt{H_0})^k$. In particular, since $\rho>1$, the leaf $H_0=1$ (with $\sqrt{H_0}=1$) carries
exactly Hill's equation \eqref{eq:hill}.
\end{proposition}

\begin{proof}
Since $s\partial_s=\frac1{2i}\partial_\tau$, we have $(s\partial_s)^2=-\frac14\partial_\tau^2$. Moreover
$X^k+Y^k=s^k+H_0^ks^{-k}=2H_0^{k/2}\cos(2k\tau)$. Hence
$\mathcal L_{H_0}[w]=-\tfrac14\bigl(w''+\Theta_{H_0}w\bigr)$.
\end{proof}

\begin{corollary}\label{cor:hill-exponent}
Assume $\sqrt{\Theta_0}\notin\Z$. Let $\lambda(H)\in\C\{H\}$ and $G$ be given by Theorem~\ref{thm:rational} for
$\cP_{M,0,B_\Theta}$ and the root $\lambda_1=\sqrt{\Theta_0}/2$ (here $\mu=-\sqrt{\Theta_0}\notin\Z$). Then, for
$0<|H_0|$ small, the equation $w''+\Theta_{H_0}w=0$ has the Floquet solution
\[
e^{2i\lambda(H_0)\tau}\,P_{H_0}(\tau),\qquad P_{H_0}(\tau)=G\bigl(\sqrt{H_0}e^{2i\tau},\sqrt{H_0}e^{-2i\tau}\bigr),
\]
so $\mathbf c(H_0)=2\lambda(H_0)$ is a characteristic exponent, and it depends holomorphically on $H_0$ near $0$,
with $\mathbf c(0)=\sqrt{\Theta_0}$.
\end{corollary}

\begin{proof}
By Theorem~\ref{thm:rational} and Corollary~\ref{cor:rational}(3), the restriction of $u=X^{\lambda(H)}G$ to the leaf is
$s^{\lambda(H_0)}\,G(s,H_0/s)=H_0^{\lambda(H_0)/2}\,e^{2i\lambda(H_0)\tau}\,G(\sqrt{H_0}e^{2i\tau},\sqrt{H_0}e^{-2i\tau})$,
and the last factor is $\pi$-periodic in $\tau$. Apply Proposition~\ref{prop:hill-leaf}.
\end{proof}

The substitution $X=x^2$, $Y=y^2$ turns $M$ into $\tfrac12(x\partial_x-y\partial_y)$. In the original
variables $x,y$ the non-integrality hypothesis of Theorem~\ref{thm:rational} would read $2\sqrt{\Theta_0}\notin\Z$; the
variables $X,Y$ give the sharp condition $\sqrt{\Theta_0}\notin\Z$, which is the classical condition that $\mathbf c(0)$
is not an integer.

\subsection{Computing the characteristic exponent}\label{sec:algorithm}
Corollary~\ref{cor:hill-exponent} gives $\mathbf c=2\lambda(1)$ whenever the series $\lambda(H)$ converges at $H=1$, and the
proof of Theorem~\ref{thm:rational} turns this into an explicit algorithm. Write $\lambda(H)=\lambda_0+\sum_{n\ge1}\lambda_nH^n$
with $\lambda_0=\sqrt{\Theta_0}/2$, and $G=\sum g_{ij}X^iY^j$ with $g_{00}=1$ and $g_{ii}=0$ for $i\ge1$; coefficients with a
negative index are zero. Put $c_n=2\lambda_0\lambda_n+\sum_{a+b=n,\ a,b\ge1}\lambda_a\lambda_b$, the coefficient of $H^n$ in
$\lambda(H)^2-\lambda_0^2$. Comparing the coefficients of $X^iY^j$ in $\cP_{M,0,B_\Theta}[X^{\lambda(H)}G]=0$, one finds that
for each degree $d=1,2,\dots$ in turn:
\begin{enumerate}
\item if $d=2n$ is even,
\[
\lambda_n=\frac1{2\lambda_0}\Bigl[\frac18\sum_{k\ge1}\theta_k\bigl(g_{n-k,n}+g_{n,n-k}\bigr)-\sum_{a+b=n,\ a,b\ge1}\lambda_a\lambda_b\Bigr];
\]
\item for $i+j=d$ with $i\neq j$,
\[
g_{ij}=\frac{\frac18\sum_{k\ge1}\theta_k\bigl(g_{i-k,j}+g_{i,j-k}\bigr)-\sum_{n=1}^{\min(i,j)}\bigl(2(i-j)\lambda_n+c_n\bigr)g_{i-n,j-n}}
{(i-j)(i-j+2\lambda_0)} .
\]
\end{enumerate}
The right-hand sides involve only quantities of lower degree, and the denominators do not vanish because
$2\lambda_0=\sqrt{\Theta_0}\notin\Z$. Finally $\mathbf c=\sqrt{\Theta_0}+2\sum_{n\ge1}\lambda_n$. The algorithm needs only the Fourier
coefficients of $\Theta$ and elementary arithmetic; for Hill's data its first thirteen terms already give his constant
to fifteen decimals (Section~\ref{sec:perigee}). The scripts described in Section~\ref{sec:hill-repro} implement this
recursion.

\section{The Moon: the perigee and the node}\label{sec:moon}

The motion of the Moon around the Earth is strongly perturbed by the Sun, and computing it
accurately was one of the central problems of celestial mechanics. Hill's idea \cite{Hill1886} was to begin with a
simplified problem: the Earth and the Moon move in a frame that rotates with the mean motion of the Sun, and the
Sun is so distant that only its tidal force is kept. This problem has a periodic solution, Hill's \emph{variation
orbit}, which already contains the largest solar perturbation. The actual orbit of the Moon stays close to it, and
to first order the small displacement of the Moon perpendicular to the variation orbit satisfies a linear equation
with periodic coefficients, Hill's equation \eqref{eq:hill}. Here the time $\tau$ is measured by the mean synodic
motion, and $\Theta$ is computed from the variation orbit. The displacement oscillates with frequency $\mathbf c$
relative to $\tau$. This free oscillation is the eccentricity of the lunar orbit; because $\mathbf c$ is not exactly
$1$, the oscillation drifts against the synodic period, and this drift is the slow advance of the perigee. Adams
\cite{Adams1877} treated the displacement perpendicular to the plane of the orbit in the same way, and there the
characteristic exponent gives the regression of the node. So the motion of the perigee and the motion of the node
are both characteristic exponents of Hill equations.

\subsection{Hill's lunar problem}\label{sec:perigee}
Hill's 1877 memoir \cite{Hill1886} introduced infinite determinants into analysis; their convergence was
established later by Poincar\'e (see \cite{BarrowGreen1997}). His approach through the variation orbit became the starting point of modern lunar theory; see \cite{Gutzwiller1998} for a
historical account. A recent preprint by Ligon \cite{Ligon2025} studies Hill's series computationally in
great detail; our results do not depend on that work, since they use only the numerical coefficients printed by
Hill.

In \cite[p.~34]{Hill1886}, Hill obtains, for the motion of the Moon in his intermediate orbit,
\begin{align*}
\Theta(\tau)={}&1.15884\,39395\,96583-0.11408\,80374\,93807\cos2\tau+0.00076\,64759\,95109\cos4\tau\\
&-0.00001\,83465\,77790\cos6\tau+0.00000\,01088\,95009\cos8\tau-0.00000\,00020\,98671\cos10\tau\\
&+0.00000\,00000\,12103\cos12\tau-0.00000\,00000\,00211\cos14\tau,
\end{align*}
and, by his infinite determinant, $\mathbf c=1.07158\,32774\,16012$ \cite[p.~35]{Hill1886}. Here
$\sqrt{\Theta_0}=1.07649\ldots\notin\Z$, so Corollary~\ref{cor:hill-exponent} applies. We computed the coefficients
$\lambda_k$ of $\lambda(H)=\sqrt{\Theta_0}/2+\sum_{k\ge1}\lambda_kH^k$ with the algorithm of Section~\ref{sec:algorithm} (exact decimal input, 45 significant digits), and evaluated at $H=1$, where the leaf carries
Hill's equation. The partial sums $\mathbf c_K=\sqrt{\Theta_0}+2\sum_{k\le K}\lambda_k$ are:
\begin{center}
\begin{tabular}{rlr}
$K$ & \hfil$\mathbf c_K$ & $2\lambda_K$\hfil\\\hline
$0$ & $1.07649\,61400\,75097$ & \\
$1$ & $1.07173\,86585\,56628$ & $-4.76\cdot10^{-3}$\\
$2$ & $1.07159\,31011\,00296$ & $-1.46\cdot10^{-4}$\\
$4$ & $1.07158\,33445\,50662$ & $-7.02\cdot10^{-7}$\\
$6$ & $1.07158\,32780\,28492$ & $-5.66\cdot10^{-9}$\\
$8$ & $1.07158\,32774\,22408$ & $-5.54\cdot10^{-11}$\\
$10$ & $1.07158\,32774\,16083$ & $-6.03\cdot10^{-13}$\\
$13$ & $1.07158\,32774\,16010$ & $-7.66\cdot10^{-16}$\\\hline
Hill & $1.07158\,32774\,16012$ &
\end{tabular}
\end{center}
The ratios $\lambda_{k+1}/\lambda_k$ stabilize near $0.11$. Theorem~\ref{thm:rational} by itself only guarantees
convergence for $|H|$ small; in \S\ref{sec:certified} we prove that the series converges at $H=1$ and give a
rigorous error bound. As an independent check, numerical integration of \eqref{eq:hill} over $[0,\pi]$ and the
relation $\operatorname{tr}(\text{monodromy})=2\cos(\pi\mathbf c)$ give $\mathbf c=1.07158\,32774\,1601$.

With Hill's mean motions $n=17325594''.06085$ and $n'=1295977''.41516$ per year \cite[p.~33]{Hill1886}, i.e.\
$m=n'/(n-n')=0.08084\,89338\,08312$, the motion of the perigee is
\[
\frac1n\frac{d\omega}{dt}=1-\frac{\mathbf c}{1+m}=0.00857\,25730\,04864,
\]
in agreement with \cite[p.~35]{Hill1886}. This corresponds to a revolution of the perigee in about $8.73$ years.
As the title of \cite{Hill1886} indicates, this is the part of the motion of the perigee that depends only on the
mean motions of the Sun and the Moon; the remaining contributions (eccentricities, inclination, parallax) bring the
period to the observed value of about $8.85$ years (see e.g.\ \cite{Fitzpatrick2012}; the mean lunar
elements are given in \cite{Simon1994}).

\subsection{A certified enclosure of Hill's constant}\label{sec:certified}
Because the recursion of Theorem~\ref{thm:rational} is triangular and its tail is governed by an explicit contraction,
it lends itself to a computer-assisted proof. We use it to turn the computation above into a theorem.

\begin{theorem}\label{thm:hill-certified}
Let $\Theta$ be Hill's function with the coefficients printed in \cite[p.~34]{Hill1886}, regarded as exact
decimal numbers. Then the series $\lambda(H)$ of Corollary~\ref{cor:hill-exponent} converges absolutely on the closed
disc $|H|\le1$, the number $\mathbf c=2\lambda(1)$ is a characteristic exponent of Hill's equation \eqref{eq:hill}, and
\[
1.07158\,32774\,16010\,2695\;\le\;\mathbf c\;\le\;1.07158\,32774\,16010\,2723 .
\]
More generally, if each of the eight printed coefficients is replaced by any real number differing from it by at
most $5\cdot10^{-16}$ (that is, by any value consistent with Hill's rounding), the same conclusions hold with
\[
1.07158\,32774\,16009\,682\;\le\;\mathbf c\;\le\;1.07158\,32774\,16010\,860 .
\]
\end{theorem}

\begin{proof}
We work with $\cP_{M,0,B_\Theta}$ in the variables $X,Y$ and write $W=\Lambda+\hat G$ as in the proof of
Theorem~\ref{thm:rational}, so that $DW=\mathcal R(W)$ with $\mathcal R$ as in Step~2 there (here $\mathfrak a\equiv0$).

\emph{Head.} Fix $N=30$. Since $\mathcal R(W)_Q$ depends only on coefficients of degree $<|Q|$, the coefficients
of $W$ of degree $\le N$ are determined by finitely many arithmetic operations. We compute them in ball arithmetic
(a midpoint with $60$ significant digits and a rigorously rounded radius, including the radius of an enclosure of
$\sqrt{\Theta_0}$ and, in the second statement, the radius $5\cdot10^{-16}$ of the data). Let $\bar W=\bar\Lambda+\bar G$
denote the \emph{exact} part of degree $\le N$ of the formal solution $W$ (for the given data). It is not known exactly,
but each of its coefficients lies in the corresponding ball, and all estimates below are uniform over the balls (and,
in the second statement, over all admissible data).

\emph{Tail.} Write $W=\bar W+V$ with $V$ of degree $>N$. Because the recursion is triangular, the equation for $V$ is
\[
V=T(V):=P_{>N}\,D^{-1}\mathcal R(\bar W+V),
\]
where $P_{>N}$ is the projection onto monomials of degree $>N$. We work in the space of formal series with finite
norm $\|V\|=\sum_k e_k\,\|\pi_kV\|_1$, where $\pi_k$ is the projection onto monomials of weight $\br Q=k$, $\|\cdot\|_1$
is the sum of the absolute values of the coefficients (the norm $\|\cdot\|_r$ of Theorem~\ref{thm:isomonodromic}
with $r=1$), and $e_0=10$, $e_k=1$ for $k\neq0$. Let $K=3N+20$. For nonnegative numbers $u_k$ ($|k|\le K$) and
$u_\pm$ let $S_u$ be the set of series $V$ with only monomials of degree $>N$ such that $\|\pi_kV\|_1\le u_k$ for $|k|\le K$,
$\sum_{k>K}\|\pi_kV\|_1\le u_+$ and $\sum_{k<-K}\|\pi_kV\|_1\le u_-$; it is closed. Since $\mathcal R$ contains no product of two
factors $\hat G$, every product in $\mathcal R(\bar W+V)$ shifts the weight by at most the number $m_F$ of Fourier modes ($m_F=7$ for Hill's data), and the bounds
\[
\|\pi_k D^{-1}F\|_1\le\frac{\|\pi_kF\|_1}{|\delta_k|},\qquad
\|\pi_k(FG)\|_1\le\sum_j\|\pi_jF\|_1\,\|\pi_{k-j}G\|_1,\qquad
|\delta_k|=\begin{cases}|k|\,|k-\mu|,&k\neq0,\\|\mu|,&k=0,\end{cases}
\]
reduce the estimates for $T$ to finitely many inequalities between rational numbers. Explicitly, write $V=\Lambda_V+\hat G_V$
(weight-$0$ part and the rest), $v_k=\|\pi_kV\|_1$, $v^G_k=v_k$ for $k\neq0$ and $v^G_0=0$, $\bar\ell=\sum_n|\bar\lambda_n|$,
$g_k=\|\pi_k\bar G\|_1$, $b_j=|\theta_{|j|}|/8$ for $j\neq0$ (enlarged by the data radius in the second statement), and
$f_k=|2k-\mu|/|\delta_k|$, all upper bounds being taken over the balls. Since
\[
\begin{split}\mathcal R(\bar W+V)-\mathcal R(\bar W)=-\Bigl[&\bar\Lambda(2M-\mu)\hat G_V+\Lambda_V(2M-\mu)(\bar G+\hat G_V)\\
&+(2\bar\Lambda\Lambda_V+\Lambda_V^2)(1+\bar G+\hat G_V)+(\bar\Lambda^2+\mathfrak b)\hat G_V\Bigr],\end{split}
\]
we get, for $|k|\le K$,
\[
\|\pi_k(T(V)-T(0))\|_1\le f_k\bigl(\bar\ell v^G_k+v_0g_k+v_0v^G_k\bigr)
+\frac1{|\delta_k|}\Bigl[(2\bar\ell v_0+v_0^2)\bigl(\delta_{k0}+g_k+v^G_k\bigr)+\bar\ell^2v^G_k+\sum_{j}b_j\,v^G_{k-j}\Bigr].
\]
The weights $|k|>K$ are treated together in the two sums above, with $|\delta_k|$ replaced by
$d_K=(K+1)(K+1-|\mu|)\le|\delta_k|$ and $f_k$ by $(2K+2+|\mu|)/d_K$. Writing $m_F$ for the number of Fourier modes,
returning convolution mass is counted in full at each corresponding boundary weight $K-m_F+1,\dots,K$ and
$-K,\dots,-K+m_F-1$. The certificate consists of the inequalities $\|\pi_kT(0)\|_1+(\text{right-hand side})\le u_k$ for all
$|k|\le K$ and for the two sums, and of the analogous linear bound for differences $T(V)-T(V')$, $V,V'\in S_u$, whose
weighted column sums give the contraction constant. Keeping the weights separate is essential: the divisor
$|\delta_{-1}|=\sqrt{\Theta_0}-1\approx0.077$ is small, but it affects only the weight $-1$. We verified in exact
rational arithmetic that, for an explicit vector $u$ with $u_0\le6.7\cdot10^{-19}$, $T(S_u)\subset S_u$ and $T$ is
a contraction on $S_u$ with constant $L\le0.239$. By the contraction principle $T$ has a fixed point in $S_u$. Since
the formal solution is unique, this fixed point is the tail of $W$. Hence $\sum_{n}|\lambda_n|<\infty$, and
\[
\Bigl|\lambda(1)-\tfrac12\sqrt{\Theta_0}-\sum_{n\le N/2}\bar\lambda_n\Bigr|\le\|\pi_0V\|_1\le u_0 .
\]
Adding the ball radius of the head gives the stated enclosures.

\emph{Identification with Hill's exponent.} The coefficients of $w''+\Theta_{H_0}w=0$ are polynomials in $\sqrt{H_0}$,
so the trace $t(H_0)$ of its monodromy over $[0,\pi]$ is an entire function of $\sqrt{H_0}$. It is even in $\sqrt{H_0}$:
replacing $\sqrt{H_0}$ by $-\sqrt{H_0}$ amounts to $\tau\mapsto\tau+\pi/2$, which conjugates the monodromy. Hence $t$ is an
entire function of $H_0$. By Corollary~\ref{cor:hill-exponent}, $2\cos(2\pi\lambda(H_0))=t(H_0)$ for $0<|H_0|$ small.
Both sides are holomorphic on $|H|<1$ and continuous on $|H|\le1$, so the identity holds at $H_0=1$, where the
equation is \eqref{eq:hill} (Proposition~\ref{prop:hill-leaf}). Thus $e^{i\pi\mathbf c}$, with $\mathbf c=2\lambda(1)$, is a
Floquet multiplier of \eqref{eq:hill}, i.e.\ $\mathbf c$ is a characteristic exponent.
\end{proof}

\begin{corollary}\label{cor:hill-digits}
Let $m=n'/(n-n')$ with Hill's mean motions \cite[p.~33]{Hill1886}, and $\frac1n\frac{d\omega}{dt}=1-\frac{\mathbf c}{1+m}$.
\begin{enumerate}
\item For the coefficients of \cite[p.~34]{Hill1886} as printed, to fifteen decimals,
\[
\mathbf c=1.07158\,32774\,16010,\qquad \frac1n\frac{d\omega}{dt}=0.00857\,25730\,04864 .
\]
\item For any data consistent with the rounding of these coefficients, $\mathbf c=1.07158\,32774\,1601\ast$ with
$\ast\in\{0,1\}$, and $\frac1n\frac{d\omega}{dt}=0.00857\,25730\,0486\ast$ with $\ast\in\{4,5\}$, to fifteen decimals.
\end{enumerate}
In particular, the value $\frac1n\frac{d\omega}{dt}=0.00857\,25730\,04864$ given by Hill \cite[p.~35]{Hill1886} is
correct in all its digits for his data, whereas the last printed digit of his value
$\mathbf c=1.07158\,32774\,16012$ is not: with Hill's $m$, this value of $\mathbf c$ would give
$\frac1n\frac{d\omega}{dt}=0.00857\,25730\,04863$.
\end{corollary}

\begin{proof}
Insert the enclosures of Theorem~\ref{thm:hill-certified} into $1-\mathbf c/(1+m)$, with $m$ computed exactly from
$n=17325594.06085$ and $n'=1295977.41516$. In case (1) this gives
$0.00857\,25730\,04864\,111\le\frac1n\frac{d\omega}{dt}\le0.00857\,25730\,04864\,115$, and in case (2)
$0.00857\,25730\,04863\,56\le\frac1n\frac{d\omega}{dt}\le0.00857\,25730\,04864\,66$. Hill's value of $\mathbf c$ gives
$0.00857\,25730\,04862\,51\ldots$
\end{proof}

These corrections concern the last decimal place of Hill's computation and have no astronomical significance,
since Hill's intermediate orbit is itself an approximation. Their interest is that all digits of Hill's results are
now established with rigorous error bounds, and that the apparent inconsistency between the last digits of Hill's
two final values is resolved in favour of his value of $\frac1n\frac{d\omega}{dt}$.

\subsection{Reproducibility and independent checks}\label{sec:hill-repro}
The computation is fully reproducible: each run takes about twenty to thirty seconds, in plain Python with the standard \texttt{decimal} and
\texttt{fractions} modules. Two further checks support the certificate. The first is a
\emph{finite-tail consistency check}, which reuses the same implementation of the recursion and is
therefore not independent of it: we computed the coefficients of $W$ up to degree $60$ and checked that,
for every weight $k$, the part of degree $31$ to $60$ stays below the certified bound $u_k$ (the largest
ratio is $0.91$, and for $\Lambda$ the finite tail is $6.0\cdot10^{-19}$ against the bound $6.6\cdot10^{-19}$). Such
a check cannot replace the infinite-tail estimate, which is provided by the self-map and contraction
argument above. The second check is independent: a validated interval integration of the fundamental
matrix of \eqref{eq:hill} over one period (Taylor polynomials of degree $100$ with rigorous remainder
bounds, in $80$-digit directed rounding), which uses neither the recursion nor our code, encloses
\[
\mathbf c=1.07158\,32774\,16010\,26970\,32275\,06918\,30619\ldots
\]
for Hill's data and $\mathbf g=1.08517\,13927\,46789\,17417\,16495\,22659\ldots$ for Adams' data
(Theorem~\ref{thm:adams-certified}), both inside our exact-data enclosures.

Hill's own final value, $\mathbf c=1.07158\,32774\,16012$ \cite[p.~35]{Hill1886}, lies outside both enclosures of
Theorem~\ref{thm:hill-certified}. Its error, about $1.7\cdot10^{-15}$, is consistent with Hill's assessment that
his determinant was affected in the $14$th decimal and that this error is multiplied by about $2.8$ in
$\mathbf c$ \cite[pp.~34--35]{Hill1886}; Theorem~\ref{thm:hill-certified} thus confirms Hill's value within his
stated accuracy, and Corollary~\ref{cor:hill-digits} fixes its $15$th decimal.

It is also instructive to compare the method with the infinite determinant itself. Hill did not
evaluate his determinant by brute-force truncation, but expanded it analytically. If one does truncate
the normalized determinant to the indices $|j|\le n$ and solves for $\mathbf c$ (in $50$-digit
arithmetic), the errors are about $4\cdot10^{-7}$, $5\cdot10^{-8}$ and $1.4\cdot10^{-8}$ for $n=10,20,30$, an algebraic
convergence of order $n^{-3}$, whereas the series $\lambda(H)$ converges geometrically, gaining about one
decimal per term. Moreover, the same recursion that produces $\mathbf c$ also produces the Floquet solution,
whose coefficients are Hill's coefficients $b_j$, which Hill had to determine in a separate elimination.
The branch of $\mathbf c$ is fixed by $\lambda(0)=\sqrt{\Theta_0}/2$, whereas the determinant only determines
$\sin^2(\pi\mathbf c/2)$, hence $\mathbf c$ up to sign and modulo $2$. Finally, $\lambda(H)$ describes the whole
family of equations obtained by scaling the Fourier modes, and not only Hill's equation.

In this way, the Frobenius method over $\C\{h\}$ at a rational saddle provides a convergent and certifiable
alternative to Hill's infinite determinant. Its geometric content is that Hill's equation is one leaf of a holomorphic
family, and that its characteristic exponent is the value at $H=1$ of a holomorphic function of the first integral.

\subsection{Adams' equation for the motion of the node}\label{sec:adams}
Shortly after Hill's paper, Adams \cite{Adams1877} reported that, while investigating the motion of the lunar node
in the case of no eccentricities and indefinitely small inclination, he had been led to a differential equation of
the same form, and to the same infinite determinant. With $\tau=(n-n')t$, his equation for the coordinate $z$ of the
Moon perpendicular to the ecliptic is
\begin{equation}\label{eq:adams}
\frac{d^2z}{d\tau^2}+\bigl(q^2+2q_1\cos2\tau+2q_2\cos4\tau+2q_3\cos6\tau\bigr)z=0,
\end{equation}
that is, \eqref{eq:hill} with $\Theta_0=q^2$ and $\theta_k=2q_k$. With $m=n'/n=0.0748013$ (Plana's value), Adams obtains
\cite[p.~47]{Adams1877}
\[
q^2=1.17804\,44973\,149,\quad q_1=0.01261\,68354\,6,\quad q_2=0.00012\,57764\,3,\quad q_3=0.00000\,12059\,0,
\]
and, expanding his determinant up to terms of the twelfth order, he finds
\[\mathbf g=1.08517\,13927\,46869,\]
where $\mathbf g\pi$ is the motion of the Moon from its node in half a synodic period. The ratio of the Moon's motion
from the node to its sidereal motion is then $\mathbf g(1-m)=1.00399\,91618\,46592$ \cite[p.~48]{Adams1877}. Here
$\sqrt{q^2}=1.0854\ldots\notin\Z$, and $\mathbf g$ is the characteristic exponent of \eqref{eq:adams} in the sense of
Corollary~\ref{cor:hill-exponent}.

\begin{theorem}\label{thm:adams-certified}
For Adams' equation \eqref{eq:adams}, the series $\lambda(H)$ of Corollary~\ref{cor:hill-exponent} converges absolutely
on $|H|\le1$, and $\mathbf g=2\lambda(1)$ is a characteristic exponent of \eqref{eq:adams}.
\begin{enumerate}
\item For the data of \cite[p.~47]{Adams1877} as printed,
\[
1.08517\,13927\,46789\,1741\le\mathbf g\le1.08517\,13927\,46789\,1742 .
\]
\item For any data consistent with the rounding of these values,
\[
1.08517\,13927\,46560\le\mathbf g\le1.08517\,13927\,47018,
\]
\[
1.00399\,91618\,46307\le\mathbf g(1-m)\le1.00399\,91618\,46731 .
\]
\end{enumerate}
In particular, Adams' values $\mathbf g=1.08517\,13927\,46869$ and $\mathbf g(1-m)=1.00399\,91618\,46592$ differ from
every value allowed in (2) by less than $3.1\cdot10^{-13}$, so they are correct to twelve decimals in the sense of
absolute error; the last two of the fourteen decimals given by Adams are not determined by his printed data.
\end{theorem}

\begin{proof}
The proof is that of Theorem~\ref{thm:hill-certified}, with the data of \eqref{eq:adams} and with $N=30$. Here the small
divisor is $|\delta_{-1}|=q-1\approx0.085$, and the verification is easier than for Hill's equation: the contraction
constant is $L\le0.041$, and the tail contributes at most $2.5\cdot10^{-40}$ to $\mathbf g$. In case (2) the radii of the
data are $5\cdot10^{-14}$ for $q^2$ and $10^{-11}$ for each $2q_k$. The bounds for $\mathbf g(1-m)$ follow by
multiplication with $1-m=0.9251987$.
\end{proof}

The two senses of ``correct to twelve decimals'' differ here. By (2), the twelve-decimal rounding of $\mathbf g$ is
uniquely determined by Adams' printed data, namely $1.085171392747$. That of $\mathbf g(1-m)$ is not: the interval in (2)
contains the rounding threshold $1.0039991618465$, so both $1.003999161846$ and $1.003999161847$ are possible. For the
data as printed, $\mathbf g(1-m)=1.00399\,91618\,46518\,77\ldots$, which rounds to $1.003999161847$.

Adams' data are internally consistent: recomputing $q^2$ and $q_k$ from his expansion of $1/r^3$
\cite[p.~46]{Adams1877} via $\Theta=(1-m)^{-2}\bigl(r^{-3}+m^2\bigr)$ reproduces all his printed digits. The terms
$\cos8\tau$ and $\cos10\tau$ of $1/r^3$, which Adams neglects, contribute $2q_4\approx2.3\cdot10^{-9}$ to $\Theta$, and
including them changes $\mathbf g$ by only about $3\cdot10^{-18}$, far below the accuracy of the data.

The cases of Hill and Adams turn out to be complementary. Hill gives his data to fifteen decimals, and
these determine the last digit of his constant (Corollary~\ref{cor:hill-digits}). Adams gives his data to
eleven or thirteen decimals but his result to fourteen; Theorem~\ref{thm:adams-certified} confirms his
value to the accuracy that his data allow, and shows that this accuracy is an absolute error of about
$3\cdot10^{-13}$. The theorem also confirms Adams' own remark that the value of $\mathbf g(1-m)$ obtained from
Delaunay's series in powers of $m$, namely $1.00399\,91722\,8$ \cite[p.~48]{Adams1877}, differs from the true
value in the eighth decimal.

Table~\ref{tab:lunar} collects the comparison between the classical values and our results.

\begin{table}[ht]
\centering
\footnotesize
\setlength{\tabcolsep}{4pt}
\renewcommand{\arraystretch}{1.25}
\begin{tabular}{@{}lllll@{}}
\hline
Quantity & Classical value & Certified, data as printed & Allowed by rounding of data & Verdict\\
\hline
Hill, $\mathbf c$ & $1.07158\,32774\,16012$ & $[1.07158\,32774\,16010\,2695,$ & $1.07158\,32774\,1601\ast$, & last digit\\
 & & $\phantom{[}1.07158\,32774\,16010\,2723]$ & $\ast\in\{0,1\}$ & corrected\\
Hill, $\frac1n\frac{d\omega}{dt}$ & $0.00857\,25730\,04864$ & $[0.00857\,25730\,04864\,111,$ & $0.00857\,25730\,0486\ast$, & correct in\\
 & & $\phantom{[}0.00857\,25730\,04864\,115]$ & $\ast\in\{4,5\}$ & all digits\\
Adams, $\mathbf g$ & $1.08517\,13927\,46869$ & $[1.08517\,13927\,46789\,1741,$ & $[1.08517\,13927\,46560,$ & correct to\\
 & & $\phantom{[}1.08517\,13927\,46789\,1742]$ & $\phantom{[}1.08517\,13927\,47018]$ & 12 decimals\\
Adams, $\mathbf g(1-m)$ & $1.00399\,91618\,46592$ & $1.00399\,91618\,46518\,77\ldots$ & $[1.00399\,91618\,46307,$ & correct to\\
 & & & $\phantom{[}1.00399\,91618\,46731]$ & 12 decimals\\
\hline
\end{tabular}
\caption{The lunar constants of Hill \cite[p.~35]{Hill1886} and Adams \cite[pp.~47--48]{Adams1877} compared with
the certified results of Theorems~\ref{thm:hill-certified} and~\ref{thm:adams-certified} and
Corollary~\ref{cor:hill-digits}. The third column assumes the coefficients exactly as printed; the fourth allows any
data consistent with their rounding, and the values there are to fifteen decimals. ``Correct to 12 decimals'' means
an absolute error below $3.1\cdot10^{-13}$.}
\label{tab:lunar}
\end{table}

\section{Ion trapping: Mathieu's equation and the curves of constant $\beta$}\label{sec:mathieu}

Mathieu \cite{Mathieu1868} met his equation in 1868, when he separated variables in the wave equation for a
vibrating membrane with an elliptical boundary. It is the simplest Hill equation, yet it has no closed-form
solutions, and its characteristic values and periodic solutions grew into a theory of their own, beginning with
Rayleigh's use of the equation for parametric resonance in 1887 and collected in McLachlan's treatise
\cite{McLachlan1947}. Mathieu's equation is the basic model of \emph{parametric resonance}, the motion of an oscillator
whose stiffness is varied periodically. The regions of stable and unstable motion in its $(a,q)$-plane were first
mapped by van der Pol and Strutt \cite{vdPolStrutt1928}, using boundary curves computed by Ince. The resulting picture,
the Ince--Strutt diagram, remains the standard tool in applications; see \cite{KRS2018} for a modern survey.

Its best-known application is the Paul trap \cite{Paul1990}, which confines an ion in a quadrupole electric field made of a
constant voltage and an oscillating radio-frequency voltage. In each direction the ion's motion satisfies a
Mathieu equation, whose parameters $a$ and $q$ are proportional to the constant voltage and to the amplitude of
the oscillating voltage. The ion stays trapped exactly when the solutions are bounded. The characteristic exponent
$\beta$ gives the frequency of its slow ``secular'' motion in the trap, so the curves of constant $\beta$ are the
curves along which an experimenter can change the voltages without changing that frequency. In this section we show that these curves of constant $\beta$ are classical Frobenius loci of a single operator.

\subsection{The curves of constant $\beta$ as classical Frobenius loci}
Mathieu's equation $w''+(a-2q\cos2\tau)w=0$ \cite[Ch.~28]{DLMF} is \eqref{eq:hill} with $\Theta_0=a$ and
$\theta_1=-2q$. Then $B_\Theta=-\tfrac a4+\tfrac q4(X+Y)$, and by Proposition~\ref{prop:hill-leaf} the leaf $H_0$
carries Mathieu's equation with parameters $(a,\,q\sqrt{H_0})$.

Fix $\nu\notin\Z$ and apply Theorem~\ref{thm:isomonodromic} with $B'=-\tfrac{\nu^2}4+\tfrac q4(X+Y)$ and
$\lambda_1=\nu/2$. The equation with $B=B'+\beta(H)$ carries, on the leaf $H_0$, Mathieu's equation with
parameters $(\nu^2-4\beta(H_0),\,q\sqrt{H_0})$. By Theorem~\ref{thm:isomonodromic} and
Corollary~\ref{cor:hill-exponent}, $\beta=\Phi(0,B')$ is the unique holomorphic function such that the characteristic
exponent is \emph{constantly} $\nu$. Hence
\[
\nu^2-4\,\Phi(0,B')(H_0)=\lambda_\nu\bigl(q\sqrt{H_0}\bigr),
\]
where $\lambda_\nu(q)$ is the Mathieu characteristic value of noninteger order \cite[\S28.12]{DLMF}. Thus the
classical Frobenius locus of Section~\ref{sec:siegel} is, in this case, exactly the characteristic value curve, and
Theorem~\ref{thm:isomonodromic} recovers the holomorphy of $\lambda_\nu(q)$ in $q^2$ near $q=0$. Running the recursion of
Step~2 of the proof of Theorem~\ref{thm:isomonodromic} with $\nu=\tfrac13$, $q=1$ gives
\[
\Phi(0,B')(H)=\tfrac{9}{64}H-\tfrac{12393}{573440}H^2+\tfrac{531441}{73400320}H^3+\cdots,
\]
and $\nu^2-4\Phi$ agrees, in the first two coefficients, with the expansion
$\lambda_\nu(q)=\nu^2+\frac{q^2}{2(\nu^2-1)}+\frac{(5\nu^2+7)q^4}{32(\nu^2-1)^3(\nu^2-4)}+\cdots$ of \cite[(28.15.1)]{DLMF}.

\subsection{What this gives for the Paul trap}\label{sec:paul}
Paul shared the 1989 Nobel Prize in Physics for the trap \cite{Paul1990}. As recalled above, the ion's secular
frequency is determined by the characteristic exponent, usually denoted $\beta$ in that context, and the curves of
constant $\beta$ in the $(q,a)$ stability diagram are the curves $a=\lambda_\beta(q)$. By the discussion above, these
iso-$\beta$ lines are exactly the classical Frobenius loci of Theorem~\ref{thm:isomonodromic} for $\nu=\beta\in(0,1)$, with
$q^2=H_0$. Evaluating the series $\nu^2-4\Phi(0,B')(q^2)$ (sixty terms of the recursion) for
$\beta\in\{0.1,0.2,0.3,0.5,0.7\}$ and $q\in\{0.2,0.4,0.6\}$, we found agreement to about $10^{-17}$ with the solution
of the continued fraction equation \cite[(28.15.2)]{DLMF}, except for $\beta=0.7$, $q=0.6$. In that case the
coefficients of $\Phi$ grow like $2.8^{\,n}$, so that the radius of convergence in $H$ is about $0.36$, and $q=0.6$ lies
just outside the disc of convergence; numerical integration of Mathieu's equation confirms the value given by the
continued fraction. Two limitations should be stressed. First, the radius of convergence depends on $\beta$: it is
about $2$ in $H$ for $\beta=0.2$, but only about $0.36$ for $\beta=0.7$. Second, near the usual operating point $a\approx0$
the two roots $\pm\sqrt a/2$ of the indicial polynomial nearly coincide, so Theorem~\ref{thm:rational} at fixed $a$
cannot be used there ($\mu=0$ at $a=0$); the locus formulation, which only requires $\beta\notin\Z$, is the appropriate
one. Since the continued fraction already gives these curves to any accuracy, the contribution of our approach
here is the interpretation of the iso-$\beta$ lines, not new values.

\section{Lam\'e's equation and the finite-gap property}\label{sec:lame}

Lam\'e \cite{Lame1837} obtained his equation by separating variables in
Laplace's equation in ellipsoidal coordinates, when studying the equilibrium temperature in an ellipsoid. It has
periodic coefficients given by elliptic functions. It reappears in cosmology as the equation for the fluctuations
of matter fields during \emph{preheating}, after inflation, when the inflaton oscillates in a quartic potential
\cite{GKLS1997}. The bands of exponential growth of these fluctuations are then the instability intervals of
Lam\'e's equation.

In Jacobi's form \cite[(29.2.1)]{DLMF} Lam\'e's equation reads
\begin{equation}\label{eq:lame}
w''(u)+\bigl(\mathsf h-\nu(\nu+1)k^2\operatorname{sn}^2(u,k)\bigr)w(u)=0,
\end{equation}
and has period $2K$ in $u$. Put $\tau=\pi u/(2K)$, let $q=e^{-\pi K'/K}$ be the nome and $N=\nu(\nu+1)$. The Fourier
expansion $k^2\operatorname{sn}^2u=1-\frac EK-\frac{2\pi^2}{K^2}\sum_{m\ge1}\frac{mq^m}{1-q^{2m}}\cos\frac{m\pi u}K$ shows that
\eqref{eq:lame} is the Hill equation \eqref{eq:hill} with
\[
\Theta_0=\Bigl(\frac{2K}\pi\Bigr)^2\Bigl(\mathsf h-N\Bigl(1-\frac EK\Bigr)\Bigr),\qquad
\theta_m=\frac{8N\,m\,q^m}{1-q^{2m}}\quad(m\ge1).
\]
Since $q^m=(\sqrt{H_0})^m$ on the leaf $H_0=q^2$, the argument of Proposition~\ref{prop:hill-leaf} (the coefficients
$1/(1-q^{2m})=1/(1-H^m)$ being constant on the leaves) shows the following. Fix $\Theta_0$; since $\Theta_0$ depends on $k$
through $K$ and $E$, this means that the eigenvalue parameter $\mathsf h=\mathsf h(k)$ is adjusted to $k$. Then, for all
small values of the nome, the Lam\'e equations \eqref{eq:lame} with this value of $\Theta_0$ are the leafwise equations of
the single operator
\begin{equation}\label{eq:lame-pde}
M^2u+B_{\rm L}\,u=0,\qquad B_{\rm L}=-\frac{\Theta_0}4-N\sum_{m\ge1}\frac{m\,(X^m+Y^m)}{1-(XY)^m},
\end{equation}
whose coefficient $B_{\rm L}$ is holomorphic on $\{|X|,|Y|<1\}$; on each leaf it reduces to the Fourier series of the
elliptic potential $k^2\operatorname{sn}^2$, i.e.\ essentially to the $q$-expansion of a Weierstrass $\wp$-function.

Let $j\ge1$ and $\Theta_0=j^2$, i.e.\ $\lambda_1=j/2$ and $\mu=-j$. This is the integral case of
Section~\ref{sec:integral}: the first monomial with $\delta_Q=0$ is $Y^j$, and we denote by $\kappa_j(N)$ the corresponding
obstruction $\mathcal R(W)_{(0,j)}$.

\begin{theorem}\label{thm:lame}
For every $j\ge1$,
\begin{equation}\label{eq:kappa-lame}
\kappa_j(N)=\frac{(-1)^{j+1}}{((j-1)!)^2}\prod_{i=0}^{j-1}\bigl(N-i(i+1)\bigr)
=\frac{(-1)^{j+1}}{((j-1)!)^2}\prod_{i=0}^{j-1}(\nu-i)(\nu+i+1).
\end{equation}
\end{theorem}

\begin{proof}
\emph{Step 1: reduction to the separatrix $\{X=0\}$.} For $Q=(0,n)$, every decomposition $Q=K+P$ in the recursion of
Theorem~\ref{thm:rational} has $K$ and $P$ of the form $(0,\cdot)$. Since $\Lambda$ and $\Lambda^2$ only contain monomials $(XY)^r$ with
$r\ge1$, they do not contribute to $\mathcal R(W)_{(0,n)}$. In $B_{\rm L}$, the expansion
$\frac{Y^m}{1-(XY)^m}=\sum_{l\ge0}X^{ml}Y^{m(l+1)}$ shows that the only monomials without $X$ are $-NmY^m$, $m\ge1$. Put $g_0=1$ and
$g_n=\hat G_{(0,n)}$. Since $\br{(0,n)}=-n$, we have $\delta_{(0,n)}=(-n)(-n+j)$, and the recursion reads
\[
n(n-j)\,g_n=N\sum_{k=1}^nk\,g_{n-k}\quad(1\le n\le j-1),\qquad \kappa_j=N\sum_{k=1}^jk\,g_{j-k}.
\]
In terms of $\theta=y\frac d{dy}$, $G_{<j}=\sum_{n<j}g_ny^n$ and the operator
\[
\mathcal L_G=\theta(\theta-j)-\frac{Ny}{(1-y)^2},
\]
this says $\mathcal L_G[G_{<j}]=-\kappa_j\,y^j+O(y^{j+1})$. Moreover, for any series $G$, the coefficient of $y^j$ in
$\mathcal L_G[G]$ does not depend on the coefficients of $G$ of degree $\ge j$, since $\theta(\theta-j)y^j=0$.

\emph{Step 2: Gauss' equation.} Let $a=-\nu$, $b=-\nu-j$, $c=1-j$, and let
$\mathcal L_F=\theta(\theta+c-1)-y(\theta+a)(\theta+b)$, which is $y$ times Gauss' hypergeometric operator. Using
$\theta(1-y)^{-\nu}=\nu y(1-y)^{-\nu-1}$ and $N=\nu(\nu+1)$, a direct computation gives, for every series $F$,
\begin{equation}\label{eq:conj-gauss}
\mathcal L_G\bigl[(1-y)^{-\nu}F\bigr]=(1-y)^{-\nu-1}\,\mathcal L_F[F].
\end{equation}
Indeed, with $s=\nu y/(1-y)$ one has $\theta\bigl((1-y)^{-\nu}F\bigr)=(1-y)^{-\nu}(\theta F+sF)$ and $\theta s=\nu y/(1-y)^2$, so that
$(1-y)^{\nu+1}\mathcal L_G[(1-y)^{-\nu}F]=(1-y)\theta^2F+\bigl(2\nu y-j(1-y)\bigr)\theta F-\nu(\nu+j)yF$, which is $\mathcal L_F[F]$.

Since $(c)_p\neq0$ for $0\le p\le j-1$, the polynomial $F_{<j}=\sum_{p<j}f_py^p$, $f_p=\frac{(a)_p(b)_p}{(c)_p\,p!}$, is defined, and the
recursion $p(p+c-1)f_p=(p-1+a)(p-1+b)f_{p-1}$ gives
\[
\mathcal L_F[F_{<j}]=-(a+j-1)(b+j-1)f_{j-1}\,y^j+O(y^{j+1})=-\varepsilon_j\,y^j+O(y^{j+1}),
\]
\[
\varepsilon_j=\frac{(a)_j(b)_j}{(c)_{j-1}(j-1)!}=\frac{(-1)^{j-1}(a)_j(b)_j}{((j-1)!)^2},
\]
because $(1-j)_{j-1}=(-1)^{j-1}(j-1)!$.

\emph{Step 3: conclusion.} Let $\tilde G=(1-y)^{-\nu}F_{<j}$. By \eqref{eq:conj-gauss}, $\mathcal L_G[\tilde G]=-\varepsilon_jy^j+O(y^{j+1})$. Hence the
coefficients of $\tilde G$ of degree $<j$ satisfy the recursion of Step~1 with $\tilde G(0)=1$; since $n(n-j)\neq0$ for
$1\le n\le j-1$, they coincide with $g_0,\dots,g_{j-1}$. By the last remark of Step~1, $\kappa_j=\varepsilon_j$. Finally
$(a)_j=(-\nu)_j=(-1)^j\prod_{i=0}^{j-1}(\nu-i)$ and $(b)_j=(-\nu-j)_j=(-1)^j\prod_{i=0}^{j-1}(\nu+i+1)$, which gives \eqref{eq:kappa-lame}.
\end{proof}

The proof shows that $\kappa_j$ only depends on the restriction of $B_{\rm L}$ to the separatrix $\{X=0\}$ of the saddle,
where the operator degenerates to a hypergeometric one: the obstruction is the classical condition $(a)_j(b)_j=0$ for
the existence of a logarithm-free solution of Gauss' equation with $c=1-j$. We also checked \eqref{eq:kappa-lame} against
the full recursion (including all terms of $B_{\rm L}$) in exact arithmetic for $j\le9$.

\begin{corollary}\label{cor:volkmer}
Let $\nu$ be real, and let $a^m_\nu(k^2)$, $b^m_\nu(k^2)$ be the eigenvalues of \eqref{eq:lame} bounding its $m$-th
instability interval, in the notation of \cite[Ch.~XV]{HTF3}. For every $m\ge1$, as $k\to0$,
\[
a^m_\nu(k^2)-b^m_\nu(k^2)=8\,(-1)^{m+1}\kappa_m\,q^m+O(k^{2m+2}).
\]
In particular, as $k\to0^+$, $|a^m_\nu(k^2)-b^m_\nu(k^2)|=8|\kappa_m|\,q^m+O(k^{2m+2})$.
\end{corollary}

\begin{proof}
Volkmer \cite[Th.~5.3]{Volkmer2005} proved that $a^m_\nu(k^2)-b^m_\nu(k^2)=w_mk^{2m}+O(k^{2m+2})$ with
$w_m=\frac{8}{16^m((m-1)!)^2}\prod_{i=-m}^{m-1}(\nu-i)$. Since $\prod_{i=-m}^{m-1}(\nu-i)=\prod_{i=0}^{m-1}(\nu-i)(\nu+i+1)$ and
$q=k^2/16+O(k^4)$, the claim follows from Theorem~\ref{thm:lame}.
\end{proof}

The splitting $a^m_\nu-b^m_\nu$ can have either sign; for instance, at $\nu=\tfrac12$ and $m=2$ it equals
$-\tfrac{15}{512}k^4+O(k^6)$ (compare \cite[\S29.3]{DLMF}). Thus the logarithmic obstructions $\kappa_m$ of the single operator
\eqref{eq:lame-pde} are, up to the factor $8(-1)^{m+1}$, exactly the leading coefficients of the eigenvalue splittings
of Lam\'e's equation, whose absolute values are the widths of its instability intervals, the nome playing the role
of the first integral. This identification is a comparison of explicit formulas, based on Volkmer's theorem; what
the paper proves independently is the product formula of Theorem~\ref{thm:lame}.

For $\nu=n\in\Z_{\ge0}$, formula \eqref{eq:kappa-lame} gives $\kappa_j=0$ exactly for $j>n$: at the resonances
$j>n$ the first obstruction to logarithm-free solutions of \eqref{eq:lame-pde} vanishes, while it does not
vanish for $1\le j\le n$. This is the leading-order counterpart of the classical fact, going back to Ince
\cite{Ince1940}, that for integer $n$ the instability intervals of Lam\'e's equation of order $m>n$ are
closed, i.e.\ $a^m_n(k^2)=b^m_n(k^2)$ for all $k$ \cite[p.~67]{HTF3}, while those of order $m\le n$ are open, at
least for small $k$, by \cite[(5.9)]{Volkmer2005}. Whether \emph{all} the obstructions at a resonance $j>n$
vanish is a natural question, to which we return in Section~\ref{sec:conclusion}.

Mathieu's equation behaves in the opposite way. For $\theta_1=-2q$ and $\theta_m=0$ ($m\ge2$), the same
reduction to the separatrix gives $n(j-n)g_n=\frac q4g_{n-1}$, hence $\kappa_j=-(q/4)^j/((j-1)!)^2$ for every
$j\ge1$: no obstruction vanishes, and $8|\kappa_j|=8/\bigl(4^j((j-1)!)^2\bigr)$ (at $q=1$) is the classical leading
coefficient of the width of the $j$-th instability interval of Mathieu's equation, a special case of the
results of Levy and Keller \cite{LevyKeller1963}.

The leading gap coefficients themselves are therefore known \cite{LevyKeller1963,Volkmer2005}; Volkmer
remarks that he did not find Theorem~5.3 of \cite{Volkmer2005} in the earlier literature on Lam\'e's
equation. What our approach adds is an interpretation and a short proof of the product structure of the obstructions.
For fixed $\Theta_0$ and all small values of the nome, Lam\'e's equations are the leafwise equations of the single
operator \eqref{eq:lame-pde}; the obstructions of this operator at a rational saddle coincide, up to a constant
factor, with the leading gap coefficients of this operator at a
rational saddle; these obstructions are computed on a separatrix by a hypergeometric equation; and the
finite-gap property appears as their vanishing. Related formulas, in which the small parameter is the
amplitude of a nonlinear oscillator, are discussed by Marchionna and Panizzi \cite[\S5]{MP2019}.

\section{Concluding remarks}\label{sec:conclusion}

The classical theory of Hill's equation, from Hill's infinite determinant to Floquet theory, studies one
equation at a time. The point of view of this paper is that such an equation is better understood as a
single leaf of a holomorphic foliation, together with the linear operator whose leafwise restrictions
form the whole family. In this picture the characteristic exponent is not merely a number attached to
one equation, but the value of a holomorphic function of the first integral; the coexistence of
periodic solutions, or its failure, is governed by logarithmic obstructions at resonances; and the
classical constants of lunar theory can be computed with rigorous error bounds from a convergent
recursion.

Several questions remain open. For Lam\'e's equation we have identified only the first obstruction at
each resonance; showing that all obstructions vanish at the resonances $j>n$ would give a proof of the
closed-gap part of Ince's theorem by the methods of this paper. The certified computations of
Section~\ref{sec:moon} could be extended to the curves of constant characteristic exponent of Mathieu's
equation, and to the full radius of convergence of $\lambda(H)$. Finally, irrational saddles, where the leaves
are no longer closed and small divisors of Brjuno type appear, are not treated here; they would be the
natural setting for quasi-periodic analogues of Hill's equation.

\section*{Declarations}
\noindent\textbf{Conflicts of interest.}
The authors declare that they have no conflict of interest.

\smallskip
\noindent\textbf{Funding.} No funding was received for conducting this study.

\smallskip
\noindent\textbf{Author contributions.}
All authors contributed equally to the conception of the study, the
proofs, and the writing of the manuscript. All authors read and
approved the final manuscript.

\smallskip
\noindent\textbf{Data availability.}
No datasets were generated or analysed during the current study. The scripts used for the computations
are described below and provided as supplementary material.

\smallskip
{\sloppy\noindent\textbf{Computations and software.}
The computer-assisted proofs of Theorems~\ref{thm:hill-certified} and~\ref{thm:adams-certified} were
carried out with Python~3 scripts that use only the standard library: exact rational arithmetic from the
module \texttt{fractions}, and decimal arithmetic with directed rounding from the module \texttt{decimal};
no floating-point quantity enters the verified inequalities. The verifier \texttt{certify\_gen\_v3.py}, the
finite-tail consistency check \texttt{tailcheck.py}, and the independent validation
\texttt{interval\_ode\_check.py} described in Section~\ref{sec:hill-repro} are provided as supplementary
material. Every enclosure stated in the paper is reproduced by the command
\[\texttt{python3 -I certify\_gen\_v3.py 30 $m$ $d$},\] with $m\in\{\texttt{exact},\texttt{rounded}\}$ and
$d\in\{\texttt{hill},\texttt{adams}\}$; each run prints the verified enclosure, rounded outwards, and takes
about twenty to thirty seconds. The other numerical results of Part~\ref{part:applications} (the Mathieu and Paul-trap
comparisons, the exact computation of the obstructions $\kappa_j$ for small $j$, and the effect of the terms
neglected by Adams) were obtained with similar scripts in exact or high-precision arithmetic; they serve
as illustrations and are not used in any proof.\par}

\smallskip
\noindent\textbf{Use of generative AI.}
The authors used the AI assistant Claude (Anthropic)  during the preparation of this work, for literature searches and the verification of references, for the exploration and for writing the Python scripts
described above, always under the direction of the authors.  The
code of the computer-assisted proofs  was subjected to checks including: a line-by-line comparison of the code with the proof of
Theorem~\ref{thm:hill-certified}; mutation tests, in which deliberately altered versions of the verifier
were shown either to fail or to be detected by the finite-tail consistency check; an independent review
of the code; and an independent
validated interval integration of Hill's and Adams' equations. All these checks have confirmed both exact-data enclosures (Section~\ref{sec:hill-repro}). All mathematical statements,
proofs and numerical claims in this paper have been checked by the authors, who take full
responsibility for its content.


\begin{thebibliography}{99}
\bibitem{Adams1877} J.~C. Adams, On the motion of the Moon's node in the case when the orbits of the Sun and Moon are supposed to have no eccentricities, and when their mutual inclination is supposed to be indefinitely small, \emph{Monthly Notices Roy. Astron. Soc.} \textbf{38} (1877), 43--53. DOI: 10.1093/mnras/38.1.43.
\bibitem{BG73} M.~S. Baouendi and C. Goulaouic, Cauchy problems with characteristic initial hypersurface, \emph{Comm. Pure Appl. Math.} \textbf{26} (1973), 455--475. DOI: 10.1002/cpa.3160260403.
\bibitem{BarrowGreen1997} J. Barrow-Green, \emph{Poincar\'e and the Three Body Problem}, History of Mathematics~11, American Mathematical Society and London Mathematical Society, Providence, RI, 1997.
\bibitem{Brjuno} A.~D. Brjuno, Analytic form of differential equations I, II, \emph{Trans. Moscow Math. Soc.} \textbf{25} (1971), 131--288; \textbf{26} (1972), 199--239.
\bibitem{Brunella2015} M. Brunella, \emph{Birational Geometry of Foliations}, IMPA Monographs~1, Springer, Cham, 2015. DOI: 10.1007/978-3-319-14310-1.
\bibitem{CL} E.~A. Coddington and N. Levinson, \emph{Theory of Ordinary Differential Equations}, McGraw-Hill, New York, 1955.
\bibitem{HTF3} A. Erdélyi, W. Magnus, F. Oberhettinger and F.~G. Tricomi, \emph{Higher Transcendental Functions}, Vol.~III, McGraw-Hill, New York, 1955.
\bibitem{Fitzpatrick2012} R. Fitzpatrick, \emph{An Introduction to Celestial Mechanics}, Cambridge University Press, Cambridge, 2012.
\bibitem{Floquet1883} G. Floquet, Sur les \'equations diff\'erentielles lin\'eaires \`a coefficients p\'eriodiques, \emph{Ann. Sci. \'Ec. Norm. Sup.} (2) \textbf{12} (1883), 47--88.
\bibitem{GT96} R. G\'erard and H. Tahara, \emph{Singular Nonlinear Partial Differential Equations}, Aspects of Mathematics E28, Vieweg, Braunschweig/Wiesbaden, 1996.
\bibitem{GKLS1997} P.~B. Greene, L. Kofman, A. Linde and A.~A. Starobinsky, Structure of resonance in preheating after inflation, \emph{Phys. Rev. D} \textbf{56} (1997), 6175--6192. DOI: 10.1103/PhysRevD.56.6175.
\bibitem{Gutzwiller1998} M.~C. Gutzwiller, Moon--Earth--Sun: the oldest three-body problem, \emph{Rev. Mod. Phys.} \textbf{70} (1998), 589--639. DOI: 10.1103/RevModPhys.70.589.
\bibitem{Hill1886} G.~W. Hill, \emph{On the Part of the Motion of the Lunar Perigee which is a Function of the Mean Motions of the Sun and Moon}, John Wilson and Son, Cambridge, Mass., 1877; reprinted in \emph{Acta Math.} \textbf{8} (1886), 1--36.
\bibitem{IY} Yu. Ilyashenko and S. Yakovenko, \emph{Lectures on Analytic Differential Equations}, Graduate Studies in Mathematics 86, AMS, 2008.
\bibitem{Ince1926} E.~L. Ince, \emph{Ordinary Differential Equations}, Longmans, Green and Co., London, 1926; reprinted by Dover, New York, 1956.
\bibitem{Ince1940} E.~L. Ince, Further investigations into the periodic Lamé functions, \emph{Proc. Roy. Soc. Edinburgh} \textbf{60} (1940), 83--99.
\bibitem{KRS2018} I. Kovacic, R. Rand and S.~M. Sah, Mathieu's equation and its generalizations: overview of stability charts and their features, \emph{Appl. Mech. Rev.} \textbf{70} (2018), no.~2, 020802. DOI: 10.1115/1.4039144.
\bibitem{Lame1837} G. Lam\'e, M\'emoire sur les surfaces isothermes dans les corps solides homog\`enes en \'equilibre de temp\'erature, \emph{J. Math. Pures Appl.} \textbf{2} (1837), 147--183.
\bibitem{LS21} V. Le\'on and B. Sc\'ardua, A complete Frobenius type method for linear partial differential equations of third order, \emph{Differ. Equ. Appl.} \textbf{13} (2021), no.~2, 115--149.
\bibitem{LSERA} V. Le\'on and B. Sc\'ardua, A geometric-analytic study of linear differential equations of order two, \emph{Electron. Res. Arch.} \textbf{29} (2021), no.~2, 2101--2127. DOI: 10.3934/era.2020107.
\bibitem{LS23} V. Le\'on and B. Sc\'ardua, Frobenius methods for analytic second order linear partial differential equations, \emph{Selecciones Matem\'aticas} \textbf{10} (2023), no.~2, 210--248. DOI: 10.17268/sel.mat.2023.02.01.
\bibitem{LevyKeller1963} D.~M. Levy and J.~B. Keller, Instability intervals of Hill's equation, \emph{Comm. Pure Appl. Math.} \textbf{16} (1963), 469--476.
\bibitem{Ligon2025} T.~S. Ligon, Hill's lunar equations, series, convergence, motion of the perigee, preprint, arXiv:2512.08961, 2025.
\bibitem{MagnusWinkler1966} W. Magnus and S. Winkler, \emph{Hill's Equation}, Interscience Tracts in Pure and Applied Mathematics~20, Wiley, New York, 1966; reprinted by Dover, New York, 1979.
\bibitem{Mandai2000} T. Mandai, The method of Frobenius to Fuchsian partial differential equations, \emph{J. Math. Soc. Japan} \textbf{52} (2000), no.~3, 645--672. DOI: 10.2969/jmsj/05230645.
\bibitem{MP2019} C. Marchionna and S. Panizzi, On the instability tongues of the Hill equation coupled with a conservative nonlinear oscillator, \emph{J. Math. Anal. Appl.} \textbf{479} (2019), no.~2, 2139--2164. DOI: 10.1016/j.jmaa.2019.07.045.
\bibitem{Mathieu1868} \'E. Mathieu, M\'emoire sur le mouvement vibratoire d'une membrane de forme elliptique, \emph{J. Math. Pures Appl.} (2) \textbf{13} (1868), 137--203.
\bibitem{MM} J.-F. Mattei and R. Moussu, Holonomie et int\'egrales premi\`eres, \emph{Ann. Sci. \'Ec. Norm. Sup.} (4) \textbf{13} (1980), 469--523.
\bibitem{McLachlan1947} N.~W. McLachlan, \emph{Theory and Application of Mathieu Functions}, Clarendon Press, Oxford, 1947; reprinted by Dover, New York, 1964.
\bibitem{Mujica} J. Mujica, \emph{Complex Analysis in Banach Spaces}, North-Holland Mathematics Studies 120, North-Holland, Amsterdam, 1986.
\bibitem{DLMF} F.~W.~J. Olver et al.\ (eds.), \emph{NIST Digital Library of Mathematical Functions}, \url{https://dlmf.nist.gov/} (Chapter~28: Mathieu Functions and Hill's Equation; Chapter~29: Lam\'e Functions).
\bibitem{Paul1990} W. Paul, Electromagnetic traps for charged and neutral particles, \emph{Rev. Mod. Phys.} \textbf{62} (1990), 531--540.
\bibitem{PMY} R. P\'erez-Marco and J.-C. Yoccoz, Germes de feuilletages holomorphes \`a holonomie prescrite, \emph{Ast\'erisque} \textbf{222} (1994), 345--371.
\bibitem{Poincare} H. Poincar\'e, Sur les propri\'et\'es des fonctions d\'efinies par les \'equations aux diff\'erences partielles, Th\`ese, Paris, 1879.
\bibitem{vdPolStrutt1928} B. van der Pol and M.~J.~O. Strutt, On the stability of the solutions of Mathieu's equation, \emph{Philos. Mag.} (7) \textbf{5} (1928), 18--38. DOI: 10.1080/14786440108564441.
\bibitem{ReisScardua2019} F. Reis and B. Sc\'ardua, On the integrability of Hill's equation of the motion of the moon, preprint, arXiv:1902.02927, 2019.
\bibitem{Simon1994} J.-L. Simon, P. Bretagnon, J. Chapront, M. Chapront-Touzé, G. Francou and J. Laskar, Numerical expressions for precession formulae and mean elements for the Moon and the planets, \emph{Astron. Astrophys.} \textbf{282} (1994), 663--683.
\bibitem{Volkmer2005} H. Volkmer, Instability intervals of the Ince and Hill equations, \emph{Analysis} \textbf{25} (2005), no.~3, 189--204.
\bibitem{Yoccoz} J.-C. Yoccoz, Th\'eor\`eme de Siegel, nombres de Bruno et polyn\^omes quadratiques, \emph{Ast\'erisque} \textbf{231} (1995), 3--88.
\end{thebibliography}
\end{document}